\documentclass[12pt]{article}
\usepackage{a4, amssymb, amsmath, exscale}
\usepackage{makeidx, graphicx, multicol}

\newenvironment{evlist}[2]{
\begin{list}{}{
\setlength{\topsep}{0.5ex plus0.2ex minus0.1ex} 
\setlength{\leftmargin}{#1}
\setlength{\itemsep}{#2 plus0.2ex}
\setlength{\listparindent}{0pt}
\setlength{\parsep}{0ex plus0.2ex} }}
{\end{list}}

\newcommand{\Real}{\mathbb{R}}
\newcommand{\Realpi}{\mathbb{R}^+_\infty}
\newcommand{\Mappi}{\mathrm{M}}

\newcommand{\specM}{\mathsf{M}}
\newcommand{\specN}{\mathsf{N}}

\newcommand{\Prob}[2]{\mathrm{P}(#1,#2)}
\newcommand{\pp}[2]{#1^{#2}_\triangleleft}
\newcommand{\Realpos}{\mathbb{R}^+}

\newcommand{\Rat}{\mathbb{Q}}
\newcommand{\Nat}{\mathbb{N}}
\newcommand{\id}{\mathrm{id}}

\newcommand{\half}{{\textstyle\frac{1}{2}}}

\newcommand{\proj}{\mathrm{p}}
\newcommand{\proq}{\mathrm{q}}
\newcommand{\Elem}{\mathrm{E}}
\newcommand{\definition}[1]{\textit{#1}}

\newcommand{\proof}{{\textit{Proof}\enspace}}

\newcommand{\eop}{\ \vbox{\hrule
                       \hbox{\vrule
                             \hskip 6pt
                             \vrule height 6pt width 0pt
                             \vrule}%
                       \hrule}%
                     \vspace{\medskipamount}
                }

\newsavebox{\ttt}
\sbox{\ttt}{}
\newcommand{\startsection}[1]
    {\section[\,#1]{#1}
    \sbox{\ttt}{\thesection\ \ \textsc{#1}}
    \thispagestyle{plain}
}

\newtheorem{lemma}{Lemma}[section]
\newtheorem{proposition}{Proposition}[section]
\newtheorem{theorem}{Theorem}[section]

\begin{document}

\begin{titlepage}

\begin{center}
\phantom{q}
\vspace{60pt}

{\huge Some Notes on}

\medskip
{\huge Standard Borel and Related Spaces}

\bigskip
{\LARGE Chris Preston}

\bigskip
\end{center}

\bigskip
\bigskip
\bigskip
\bigskip
\begin{quote}
These notes give an elementary approach to parts of the theory of standard Borel and analytic spaces.
\end{quote}

\end{titlepage}

\thispagestyle{empty}

\addtocontents{toc}{\vskip 20pt}
\addtolength{\parskip}{-5pt}
\tableofcontents
\addtolength{\parskip}{5pt}
\startsection{Introduction}
\label{intro}

A \definition{measurable space} is a pair $(X,\mathcal{E})$ consisting of a non-empty set $X$ together with a $\sigma$-algebra $\mathcal{E}$ of subsets of $X$.
If $(X,\mathcal{E})$ and $(Y,\mathcal{F})$ are measurable spaces then a mapping $f : X \to Y$ is said to be \definition{measurable}
if $f^{-1}(\mathcal{F}) \subset \mathcal{E}$. In order to show the dependence on the $\sigma$-algebras $\mathcal{E}$ and $\mathcal{F}$ we then say that
$f : (X,\mathcal{E}) \to (Y,\mathcal{F})$ is measurable. Many of the mappings which occur here have a property which is stronger than just being measurable 
in that $f^{-1}(\mathcal{F}) = \mathcal{E}$ holds, and in this case we say that
$f : (X,\mathcal{E}) \to (Y,\mathcal{F})$ is \definition{exactly measurable}.

A measurable space $(X,\mathcal{E})$ is said to be \definition{standard Borel} if there exists a metric on $X$ which makes
it a complete separable metric space in such a way that $\mathcal{E}$ is then the Borel $\sigma$-algebra (this being the smallest $\sigma$-algebra
containing the open sets).
The name `standard Borel' was given to such spaces by Mackey in \cite{mackey} and they are important because there 
are several very useful results which, although they do not hold in general, are true for standard Borel spaces.
For example, if $(X,\mathcal{E})$ and $(Y,\mathcal{F})$ are standard Borel then any bijective measurable mapping
$f : (X,\mathcal{E}) \to (Y,\mathcal{F})$ is automatically an isomorphism (i.e., the inverse mapping is also measurable).
Other examples involve the existence of conditional probability kernels and generalisations of the classical Kolmogorov
extension theorem.

It is often the case that a standard Borel space $(X,\mathcal{E})$ arises from a topological space $X$, but that
the metric occurring in the definition of being standard Borel is not the `natural' one. This has led to the notion of a 
\definition{Polish space}, which
is a separable topological space whose topology can be given by a complete metric.
Thus, for example, although its usual metric is not complete, the open unit interval $(0,1)$ is a Polish space since it is 
homeomorphic to $\Real$, whose usual metric is complete. 

The theory of standard Borel spaces is usually presented as a spin-off of the theory of Polish spaces.  
In these notes we give an alternative treatment, which essentially only uses a single Polish space.
The space involved here is 
$\specM = \{0,1\}^{\Nat}$ (the space of all sequences $\{z_n\}_{n\ge 0}$ of 0's and 1's), considered as a topological space
as the product of $\Nat$ copies of the discrete space $\{0,1\}$.
Thus $\specM$ is compact and the topology is induced, for example, by the metric
\[ d(\{z_n\}_{n\ge 0},\{z'_n\}_{n\ge 0}) = \sum_{n\ge 0} 2^{-n}|z_n - z'_n|\;. \]
The $\sigma$-algebra of Borel subsets of $\specM$ will be denoted by $\mathcal{B}$.

A measurable space $(X,\mathcal{E})$ is said to be \definition{countably generated} if 
$\mathcal{E} = \sigma(\mathcal{S})$ for some countable subset $\mathcal{S}$ of $\mathcal{E}$ and is said to be
\definition{separable} if $\{x\} \in \mathcal{E}$ for each $x \in X$.
In particular, a standard Borel space is both countably generated and separable. (It is countably generated 
since a separable metric space has a countable base for its topology.)

The starting point for our approach  is a result of Mackey (Theorem~2.1 in \cite{mackey}).
This states that a measurable space $(X,\mathcal{E})$ is countably generated if and only if
there exists an exactly measurable mapping $f : (X,\mathcal{E}) \to (\specM,\mathcal{B})$.
Moreover, it is straightforward to show that if $(X,\mathcal{E})$ is standard Borel then
$f$ can be chosen so that $f(X) \in \mathcal{B}$. (If
$X$ is a complete separable metric space then there is a standard construction producing a continuous 
injective mapping $h : X \to [0,1]^\Nat$ such that $h$ is a homeomorphism from $X$ to $h(X)$ (with the relative topology) and such that
$h(X)$ is a $G_\delta$ subset of $[0,1]^\Nat$, i.e., $h(X)$ is 
the intersection of a sequence of open sets.
The mapping $f$ is obtained by composing $h$ with a suitable mapping $g : [0,1]^\Nat \to \specM$. The details of this construction can be found  in the
proof of Theorem~\ref{type-ba}.1.) 
In fact, it is also well-known (i.e., well-known to those who are interested in such things)
that if $(X,\mathcal{E})$ is standard Borel then $f(X) \in \mathcal{B}$ for every
exactly measurable mapping $f : (X,\mathcal{E}) \to (\specM,\mathcal{B})$. (See Proposition~\ref{type-ba}.1.)

There is one ingredient which we have not yet used (and which might help to make the results stated above more familiar):
If $(X,\mathcal{E})$ is a separable countably generated measurable space then any exactly measurable 
$f : (X,\mathcal{E}) \to (\specM,\mathcal{B})$ is injective; moreover it induces an isomorphism
$f : (X,\mathcal{E}) \to (A,\mathcal{B}_{|A})$, where $A = f(X)$ and $\mathcal{B}_{|A}$ is the trace $\sigma$-algebra of $\mathcal{B}$ on $A$.

If $\mathcal{D}$ is a subset of $\mathcal{P}(\specM)$ containing $\mathcal{B}$ and closed under finite intersections then we call
a countably generated measurable space $(X,\mathcal{E})$ a \definition{type $\mathcal{D}$ space} if there exists an
exactly measurable mapping $f : (X,\mathcal{E}) \to (\specM,\mathcal{B})$ such that $f(X) \in \mathcal{D}$.
In particular, a standard Borel space is a separable type $\mathcal{B}$ space; moreover, the converse also holds: This follows from the fact 
(also well-known and proved in Proposition~\ref{anal-subsets}.8)
that if $B$ is an uncountable element of $\mathcal{B}$ then $(B,\mathcal{B}_{|B})$ is isomorphic to $(\specM,\mathcal{B})$.

Type $\mathcal{B}$ spaces should thus be thought of as standard Borel spaces without the assumption of separability, and
in these notes we actually study type $\mathcal{B}$  and not standard Borel spaces. 
One reason for not requiring separability is that there are situations in which the typical results holding for standard Borel spaces
are needed, but the spaces involved cannot be separable. For example, in the set-up of 
the Kolmogorov extension theorem there is a measurable space $(X,\mathcal{E})$ and an increasing sequence $\{\mathcal{E}_n\}_{n\ge 0}$
of sub-$\sigma$-algebras of $\mathcal{E}$ with $\mathcal{E} = \sigma(\bigcup_{n\ge 0} \mathcal{E}_n)$. Now although $(X,\mathcal{E})$
is usually separable this will not be the case for the measurable spaces $(X,\mathcal{E}_n)$, $n \ge 0$. 
However, these spaces need to be `nice', a requirement which will be met if they are all type $\mathcal{B}$ spaces.
In fact, the notion of a type $\mathcal{B}$ space already occurs implicitly in Parthasarthy's 
proof of the Kolmogorov extension theorem (in Chapter~V of \cite{partha}).

Besides type $\mathcal{B}$ spaces we also need to consider type $\mathcal{A}$ spaces, where $\mathcal{A}$ is the set of analytic subsets of $\specM$.
A subset $A$ of $\specM$ is said to be \definition{analytic} if it is either empty or there exists a continuous mapping
$\tau : \specN \to \specM$ with $\tau(\specN) = A$. 
Here $\specN = \Nat^{\Nat}$ (the space of all sequences $\{m_n\}_{n\ge 0}$ of elements from $\Nat$), considered as a topological space
as the product of $\Nat$ copies of $\Nat$ (with the discrete topology). This topology is also induced by the complete metric
$d : \specN \times \specN \to \Realpos$ given by
\[ d(\{m_n\}_{n\ge 0},\{m'_n\}_{n\ge 0}) = \sum_{n\ge 0} 2^{-n}\delta'(m_n,m'_n) \]
where $\delta'(m,m) = 0$ and $\delta'(m,n) = 1$ whenever $m \ne n$. 
(The statement made earlier which suggested that $\specM$ is the only Polish space to  occur in these notes
is thus not quite correct, since the Polish space $\specN$
is essential for the definition of the analytic sets.)
Type $\mathcal{A}$ spaces should thus be thought of as what are usually called analytic spaces, but again
without the assumption of separability.

In the first part of these notes (Sections \ref{e-meas-maps} to \ref{selectors}) we develop the theory of type $\mathcal{B}$ and type $\mathcal{A}$
spaces and obtain the results corresponding to the standard facts which hold for standard Borel and analytic spaces.
Most of the non-trivial properties of type $\mathcal{B}$ and type $\mathcal{A}$ spaces depend on results about the analytic subsets of $\specM$
and the proofs of these results (in Section~\ref{anal-subsets}) are all based
on the corresponding proofs in Chapter~8 of Cohn \cite{cohn}.

In the second part (Sections \ref{image-measures} to \ref{dynkin}) we show that for certain kinds of applications
it seems to be more natural to work with type $\mathcal{B}$ spaces directly rather than with the usual definition of a standard Borel
space. These applications all involve constructing probability measures and, roughly speaking, our approach to this situation is the following: We have
a type $\mathcal{B}$ (or a type $\mathcal{A}$) space $(X,\mathcal{E})$ and want to construct something out of a set
$S \subset \Prob{X}{\mathcal{E}}$ of probability measures defined on $(X,\mathcal{E})$.
By definition there exists an exactly measurable mapping $f : (X,\mathcal{E}) \to (\specM,\mathcal{B})$  with
$f(X) \in \mathcal{B}$ (or $f(X) \in \mathcal{A}$) and we consider the image measures
$\{ \mu f^{-1} : \mu \in S \}$ which are probability measures defined on $(\specM,\mathcal{B})$. We then carry out
the construction on these measures (exploiting the special properties of the space $(\specM,\mathcal{B})$) and pull the result back to
the measurable space $(X,\mathcal{E})$. It is this last step where the fact that 
$f(X) \in \mathcal{B}$ (or $f(X) \in \mathcal{A}$) plays a crucial role.

The precise formulation of the method is given in Section~\ref{image-measures}.
In Section~\ref{kolmogorov} it is applied to show that the Kolmogorov extension property holds for an inverse limit of type $\mathcal{A}$ spaces 
(a form of the extension theorem given in Chapter V of Parthasarathy \cite{partha}). In Section~\ref{pp} we show that
the space of finite point processes defined on a type $\mathcal{B}$ space is itself a type $\mathcal{B}$ space
(a fact which is equivalent to results in Matthes, Kerstan and Mecke \cite{mkm}, Kallenberg \cite{kall} and Bourbaki \cite{bourb}).
Section~\ref{cond-dist} looks at the existence of conditional distributions and gives a proof of the usual result for
standard Borel spaces to be found, for example, in
Doob \cite{doob},  Parthasarathy \cite{partha} or Dynkin and Yushkevich \cite{dynkinyu}. 
Finally, in Section~\ref{dynkin} we consider the construction of a particular kind of entrance boundary for random fields
due to F\"ollmer \cite{foellmer} (based on ideas in  
Dynkin \cite{dynkin}).


\startsection{Exactly measurable mappings}
\label{e-meas-maps}

Recall that if $(X,\mathcal{E})$ and $(Y,\mathcal{F})$ are measurable spaces then a measurable mapping
$f : (X,\mathcal{E}) \to (Y,\mathcal{F})$ is said to be \definition{exactly measurable} if
$f^{-1}(\mathcal{F}) = \mathcal{E}$.
In this section we look at some general properties of exactly measurable mappings.

\begin{lemma}
(1)\enskip
Let $(X,\mathcal{E})$ be a measurable space. Then the identity mapping $\id_X : (X,\mathcal{E}) \to (X,\mathcal{E})$ is exactly measurable.

(2)\enskip
Let $(X,\mathcal{E})$, $(Y,\mathcal{F})$ and $(Z,\mathcal{G})$ 
be measurable spaces and let $f : (X,\mathcal{E}) \to (Y,\mathcal{F})$ and  $g : (Y,\mathcal{F}) \to (Z,\mathcal{G})$ 
be exactly measurable mappings. Then the composition
$g \circ f : (X,\mathcal{E}) \to (Z,\mathcal{G})$ is also exactly measurable.
\end{lemma}

\proof 
This is clear. \eop

Let $A$ be a non-empty subset of a non-empty set $X$; for each subset
$\mathcal{S}$ of $\mathcal{P}(X)$ denote by $\mathcal{S}_{|A}$ the subset of
$\mathcal{P}(A)$ consisting of all sets having the form $S \cap A$ for some
$S \in \mathcal{S}$. Then $\mathcal{S}_{|A}$ is referred to as the 
\definition{trace of $\mathcal{S}$ on $A$}. 
If $\mathcal{E}$ is a $\sigma$-algebra of subsets of $X$ then $\mathcal{E}_{|A}$ is a $\sigma$-algebra of
subsets of $A$.

\begin{proposition}
Let $f : (X,\mathcal{E}) \to (Y,\mathcal{F})$ be a measurable mapping,
let $A$ be a non-empty subset of $X$ and $f_{|A} : A \to Y$ be the restriction of $f$ to $A$.
Then $f_{|A} : (A,\mathcal{E}_{|A}) \to (Y,\mathcal{F})$ is measurable.
Moreover, if $f$ is exactly measurable then so is $f_{|A}$.
\end{proposition}

\proof 
Note that $f_{|A}^{-1}(F) = f^{-1}(F) \cap A$ for all $F \subset Y$.
Suppose first that $f$ is measurable 
and let $F \in \mathcal{F}$; then $f_{|A}^{-1}(F) = f^{-1}(F) \cap A \in \mathcal{E}_{|A}$, which implies
that $f_{|A}^{-1}(\mathcal{F}) \subset \mathcal{E}_{|A}$, i.e., $f_{|A}$ is measurable.
Suppose now that $f$ is exactly measurable and let $E \in \mathcal{E}_{|A}$; thus $E = E' \cap A$ for some $E' \in \mathcal{E}$ and there exists
$F \in \mathcal{F}$ with $E' = f^{-1}(F)$. Then $f_{|A}^{-1}(F) = f^{-1}(F) \cap A = E' \cap A = E$ and therefore
$f_{|A}^{-1}(\mathcal{F}) = \mathcal{E}_{|A}$, i.e., $f_{|A}$ is exactly measurable. \eop

Let $(X,\mathcal{E})$ be a measurable space and $A$ be a non-empty subset of $X$. Applying Proposition~\ref{e-meas-maps}.1
to the identity mapping $\id_X$ shows that the inclusion mapping
$i_A : A \to X$ gives rise to an exactly measurable mapping
$i_A : (A,\mathcal{E}_{|A}) \to (X,\mathcal{E})$.

A mapping $f : X \to Y$ with $A = f(X)$ will also be considered  as a (surjective) mapping $f : X \to A$.

\begin{proposition}
Let $(X,\mathcal{E})$ and $(Y,\mathcal{F})$ be measurable spaces and $f : X \to Y$ be any mapping; put $A = f(X)$. Then
$f : (X,\mathcal{E}) \to (Y,\mathcal{F})$ is measurable (resp.\ exactly measurable) if and only if
$f : (X,\mathcal{E}) \to (A,\mathcal{F}_{|A})$ is measurable (resp.\ exactly measurable).
\end{proposition}

\proof 
This follows immediately from Lemma~\ref{e-meas-maps}.2~(1) below, since the elements in $\mathcal{F}_{|A}$ are exactly the sets of the 
form $F \cap A$ with $F \in \mathcal{F}$. \eop

Here are some useful simple properties which hold for a general mapping. One of these was needed in the proof of Proposition~\ref{e-meas-maps}.2
and the rest will be needed later.

\begin{lemma}
Let $f : X \to Y$ be an arbitrary mapping and put $A = f(X)$. Then:

(1)\enskip $f^{-1}(F \cap A) = f^{-1}(F)$ for all $F \subset Y$. 

(2)\enskip  $f(f^{-1}(F)) =  F \cap A$  for all $F \subset Y$, and so $f(f^{-1}(F)) = F$ for all $F \subset A$.

(3)\enskip If $E = f^{-1}(F)$ for some $F \subset Y$ then $f(E) = F \cap A$ and $f^{-1}(f(E)) = E$.

(4)\enskip If $f$ is injective then $f^{-1}(f(E)) = E$ for all $E \subset X$.
\end{lemma}

\proof 
(1) and (4) are clear.

(2)\enskip If $y \in f(f^{-1}(F))$ then there exists $x \in f^{-1}(F)$ with $y = f(x)$ and then $y \in F$.
Hence $y \in F \cap A$, i.e., $f(f^{-1}(F)) \subset F \cap A$. On the other hand, if
$y \in F \cap A$ then there exists $x \in X$ with $y = f(x)$, thus $x \in f^{-1}(F)$ and so
$y \in f(f^{-1}(F))$, i.e.,  $F \cap A \subset f(f^{-1}(F))$. 

(3)\enskip 
If $E = f^{-1}(F)$ then by (2) $f(E) = f(f^{-1}(F)) = F \cap A$ and then by (1) it follows that $f^{-1}(F \cap A) = f^{-1}(F) = E$.
\eop

If $f : (X,\mathcal{E}) \to (Y,\mathcal{F})$ is measurable and $A = f(X)$ then
by Lemma~\ref{e-meas-maps}.2~(1) and (2) $f(f^{-1}(F)) = F$ for all $F \in \mathcal{F}_{|A}$.
Note that this has nothing to do with measurability and only depends on the fact that $F \subset A$ for each
$F \in \mathcal{F}_{|A}$. The corresponding statement (that $f^{-1}(f(E)) = E$ for all $E \in \mathcal{E}$) in the following result 
is, however, nowhere near so harmless.

\begin{proposition}
Let $f : (X,\mathcal{E}) \to (Y,\mathcal{F})$ be measurable and put $A = f(X)$. Then $f$ is exactly
measurable if and only if $f(E) \in \mathcal{F}_{|A}$ and $f^{-1}(f(E)) = E$ for all $E \in \mathcal{E}$.
\end{proposition}

\proof 
Assume first $f$ is exactly measurable. Let $E \in \mathcal{E}$; then there exists $F \in \mathcal{F}$ such that $E = f^{-1}(F)$ and therefore  by
Lemma~\ref{e-meas-maps}.2~(3) $f(E) = F \cap A \in \mathcal{F}_{|A}$ and $f^{-1}(f(E)) = E$. 
Suppose conversely that $f(E) \in \mathcal{F}_{|A}$ and $f^{-1}(f(E)) = E$ for all $E \in \mathcal{E}$.
Let $E \in \mathcal{E}$; then $f(E) \in  \mathcal{F}_{|A}$ and so $f(E) = F \cap A$ for some $F \in \mathcal{F}$. Thus by
Lemma~\ref{e-meas-maps}.2~(1) $f^{-1}(F) = f^{-1}(F \cap A) = f^{-1}(f(E)) = E$ and this shows that $f^{-1}(\mathcal{F}) = \mathcal{E}$.
\eop

A measurable mapping $f : (X,\mathcal{E}) \to (Y,\mathcal{F})$ is   \definition{bimeasurable} 
if $f(E) \in \mathcal{F}$ for all $E \in \mathcal{E}$ and is an \definition{isomorphism} if it is 
bimeasurable and bijective.
Proposition~\ref{e-meas-maps}.3 implies that an exactly measurable mapping
$f$ is bimeasurable when considered as a mapping from $(X,\mathcal{E})$ to $(A,\mathcal{F}_{|A})$ with $A = f(X)$. 
If $f$ is an isomorphism then 
$f^{-1} : (Y,\mathcal{F}) \to (X,\mathcal{E})$ is also measurable
(and is an isomorphism), where here $f^{-1} : Y \to X$ is the set-theoretical inverse of $f$. 
By Proposition~\ref{e-meas-maps}.3 a bijective measurable mapping $f$ is an isomorphism if and only if
it is exactly measurable (since $f$ being injective implies that $f^{-1}(f(E)) = E$ for all $E \subset X$).

\begin{lemma}
Let $f : (X,\mathcal{E}) \to (Y,\mathcal{F})$ be exactly  measurable, let $N$ be a finite or
countable infinite set and for each $n \in N$ let $E_n \in \mathcal{E}$. Then
\[ \bigcup_{n \in N} f(E_n) = f\Bigl(\bigcup_{n \in N} E_n\Bigr)
\ \ \mbox{and}\ \   \bigcap_{n \in N} f(E_N) = f\Bigl(\bigcap_{n \in N} E_n\Bigr)\;. \]
Moreover, the sets $f(E_1)$ and $f(E_2)$ are disjoint whenever
$E_1$ and $E_2$ are disjoint sets in $\mathcal{E}$.
\end{lemma}

\proof
Put $A = f(X)$.
By Proposition~\ref{e-meas-maps}.3 $f(E_n) \in \mathcal{F}_{|A}$ for each $n \in N$ and thus also 
$\bigcup_{n \in N} f(E_n) \in \mathcal{F}_{|A}$. Therefore by Lemma~\ref{e-meas-maps}.2~(2) and Proposition~\ref{e-meas-maps}.3
\[ \bigcup_{n \in N} f(E_n)  = f\Bigl( f^{-1}\Bigl( \bigcup_{n \in N} f(E_n)\Bigr)\Bigr) 
= f\Bigl( \bigcup_{n \in N} f^{-1}(f(E_n))\Bigr) = f \Bigl(\bigcup_{n \in N} E_n \Bigr)\;.\]
The analogous statement with $\bigcap$ instead of $\bigcup$ follows in the exactly the same way. 
Finally, if $E_1,\,E_2 \in \mathcal{E}$ are disjoint then
$f(E_1) \cap f(E_2) = f(E_1 \cap E_2) = f(\varnothing) = \varnothing$ and so
$f(E_1)$ and $f(E_2)$ are also disjoint. \eop

A subset $\mathcal{S}$ of $\mathcal{P}(X)$ (with $X$ a non-empty set) is said to \definition{separate the points of $X$}
if for each $x_1,\,x_2 \in X$ with $x_1 \ne x_2$ there exists $A \in \mathcal{S}$ containing exactly one of $x_1$ and $x_2$.

\begin{lemma}
The subset $\mathcal{S}$ separates the points of $X$ if and only if $\sigma(\mathcal{S})$ does.
\end{lemma}

\proof 
For each $x_1,\,x_2 \in X$ with $x_1 \ne x_2$
let $\mathcal{T}_{x_1,x_2}$ denote the set of subsets of $X$ which contain either both or neither of the elements $x_1$ and $x_2$;
clearly $\mathcal{T}_{x_1,x_2}$ is a $\sigma$-algebra.
Suppose $\mathcal{S}$ does not separate the points of $X$; then there
exist $x_1,\,x_2 \in X$ with $x_1 \ne x_2$ such that $\mathcal{S} \subset \mathcal{T}_{x_1,x_2}$.
Hence $\sigma(\mathcal{S}) \subset \mathcal{T}_{x_1,x_2}$, which means that
$\sigma(\mathcal{S})$ does not separate the points of $X$. The converse holds trivially: If
$\mathcal{S}$ separates the points of $X$ then so does $\sigma(\mathcal{S})$, since
$\mathcal{S} \subset \sigma(\mathcal{S})$.
\eop

A measurable space $(X,\mathcal{E})$ is \definition{separated} if $\mathcal{E}$ separates the points of $X$, and is 
\definition{separable} if $\{x\} \in X$ for all $x \in X$. 
In particular, a separable measurable space is separated.

\begin{proposition}
An exactly measurable mapping $f : (X,\mathcal{E}) \to (Y,\mathcal{F})$ with 
a separated measurable space $(X,\mathcal{E})$ is injective.
\end{proposition}

\proof 
Let $x_1,\,x_2 \in X$ with $x_1 \ne x_2$. Since $(X,\mathcal{E})$ is separated there exists $E \in \mathcal{E}$ with
$x_1 \in E$ and $x_2 \in X \setminus E$. Thus by Lemma~\ref{e-meas-maps}.3 the sets $f(E)$ and  $f(X \setminus E)$ are disjoint
and hence $f(x_1) \ne f(x_2)$. This shows that $f$ is injective. \eop

The following simple corollary of the previous results will be needed several times later:

\begin{lemma}
Let $f : (X,\mathcal{E}) \to (Y,\mathcal{F})$ be measurable with 
$(X,\mathcal{E})$ separated and  put $B = f(X)$.
Then the mapping
$f : (X,\mathcal{E}) \to (B,\mathcal{F}_{|B})$ is an isomorphism if and only if
$f$ is exactly measurable.
\end{lemma}

\proof This follows from
Propositions \ref{e-meas-maps}.2,  \ref{e-meas-maps}.3  and \ref{e-meas-maps}.4. \eop

In Propositions \ref{e-meas-maps}.5 and \ref{e-meas-maps}.6
let $(X,\mathcal{E})$ and $(Y,\mathcal{F})$ be arbitrary and $(Z,\mathcal{G})$
be a separable measurable space.

\begin{proposition}
Let $f : (X,\mathcal{E}) \to (Y,\mathcal{F})$ be an exactly measurable mapping
and let $g : (X,\mathcal{E}) \to (Z,\mathcal{G})$ be measurable.
Put $A = f(X)$; then there exists a unique mapping $h : A \to Z$ such that 
$h \circ f = g$, and then $h : (A,\mathcal{F}_{|A}) \to (Z,\mathcal{G})$ is measurable.
Moreover, if $g : (X,\mathcal{E}) \to (Z,\mathcal{G})$ is exactly measurable
then so is $h : (A,\mathcal{F}_{|A}) \to (Z,\mathcal{G})$.
\end{proposition}

\proof
Let $x \in X$; then $\{g(x)\} \in \mathcal{Z}$, since 
$(Z,\mathcal{G})$ is separable, and therefore $E = g^{-1}(\{g(x)\})$ is an element of $\mathcal{E}$ containing $x$.
Hence $\{f(x)\} \subset f(E)$ and so by Proposition~\ref{e-meas-maps}.3 $f^{-1}(\{f(x)\}) \subset f^{-1}(f(E)) = E = g^{-1}(\{g(x)\})$.
This shows that if $y \in A$ then
$f^{-1}(\{y\}) \subset g^{-1}(\{g(x)\})$ for each $x \in X$ with $f(x) = y$ (and note that $f^{-1}(\{y\}) \ne \varnothing$, since $y \in A$).
But the sets
$g^{-1}(\{z_1\})$ and $g^{-1}(\{z_2\})$ are disjoint if $z_1 \ne z_2$, which implies that if $x_1,\,x_2 \in X$ are such that $f(x_1) = y = f(x_2)$
then $g(x_1) = g(x_2)$ (since $g^{-1}(\{g(x_1)\})$ and $g^{-1}(\{g(x_2)\})$ both contain $f^{-1}(\{y\})$ and hence are not disjoint).
There thus exists a unique mapping $h : A \to Z$ such that $h \circ f = g$.
Now let $G \in \mathcal{G}$; then by Lemma~\ref{e-meas-maps}.2~(2)
\[ h^{-1}(G) = f(f^{-1}(h^{-1}(G)) \cap A) = f(f^{-1}(h^{-1}(G))) = f(g^{-1}(G)) \]
and by Proposition~\ref{e-meas-maps}.3 $f(g^{-1}(G)) \in \mathcal{F}_{|A}$, since $g^{-1}(G) \in \mathcal{E}$, i.e., $h^{-1}(G) \in \mathcal{F}_{|A}$.
This shows that
$h^{-1}(\mathcal{G}) \subset  \mathcal{F}_{|A}$, and so
$h : (A,\mathcal{F}_{|A}) \to (Z,\mathcal{G})$ is measurable.

Now suppose that $g : (X,\mathcal{E}) \to (Z,\mathcal{G})$ is exactly measurable.
Let $F \in \mathcal{F}_{|A}$; then $F = F' \cap A$ for some $F' \in \mathcal{F}$ and, since
$f^{-1}(F') \in \mathcal{E}$ and $g$ is exactly measurable,
there exists $G \in \mathcal{G}$ with $g^{-1}(G) = f^{-1}(F')$. Thus by Lemma~\ref{e-meas-maps}.2~(2)
\begin{eqnarray*}
  h^{-1}(G) = f(f^{-1}(h^{-1}(G))) &=& f( (h\circ f)^{-1}(G)) \\
&=& f(g^{-1}(G)) = f(f^{-1}(F')) = F' \cap A = F 
\end{eqnarray*}
and this shows that $h^{-1}(\mathcal{G}) = \mathcal{F}_{|A}$. Hence $h$ is exactly measurable.
\eop

\begin{proposition}
Let $f : (X,\mathcal{E}) \to (Y,\mathcal{F})$ and $g : (X,\mathcal{E}) \to (Z,\mathcal{G})$ be exactly measurable
mappings with $f$ surjective. Then there exists a unique $h : Y \to Z$ such that $h \circ f = g$, 
and the mapping $h : (Y,\mathcal{F}) \to (Z,\mathcal{G})$ 
is exactly measurable.
\end{proposition}

\proof 
This is just a special case of Proposition~\ref{e-meas-maps}.5.
\eop


We finish this section by looking at how exactly measurable mappings behave in relation to the product and disjoint union of
measurable spaces.

In the following result let
$(X,\mathcal{E})$, $(Y,\mathcal{F})$ and $(Z,\mathcal{G})$ be arbitrary measurable spaces.

\begin{proposition}
Let $f : (X,\mathcal{E}) \to (Y,\mathcal{F})$ and $g : (X,\mathcal{E}) \to (Z,\mathcal{G})$ be measurable
mappings and let $h : X \to Y \times Z$ be the mapping with $h(x) = (f(x),g(x))$ for all $x \in X$.
Then $h : (X,\mathcal{E}) \to (Y \times Z,\mathcal{F} \times \mathcal{G})$
is measurable. Moreover, $h$ is exactly measurable provided at least one of $f$ and $g$ is exactly measurable.
\end{proposition}

\proof 
If $F \in \mathcal{F}$ and $G \in \mathcal{G}$ then $h^{-1}(F \times G) = f^{-1}(F) \cap g^{-1}(G) \in \mathcal{E}$, and thus
$h^{-1}(\mathcal{R}) \subset \mathcal{E}$, where $\mathcal{R}$ is the set of measurable rectangles in $\mathcal{F} \times \mathcal{G}$
(i.e., the sets of the form $F \times G$ with $F \in \mathcal{F}$, $G \in \mathcal{G}$).
But by definition $\mathcal{F} \times \mathcal{G} = \sigma(\mathcal{R})$ and so
\[h^{-1}(\mathcal{F} \times \mathcal{G}) = h^{-1}(\sigma(\mathcal{R})) = \sigma (h^{-1}(\mathcal{R})) \subset \sigma(\mathcal{E}) \subset \mathcal{E}\;, \]
which shows that $h : (X,\mathcal{E}) \to (Y \times Z,\mathcal{F} \times \mathcal{G})$ is measurable.
Suppose now that $f$ is exactly measurable and let $E \in \mathcal{E}$; then there exists $F \in \mathcal{F}$ with
$f^{-1}(F) = E$. Hence $F \times Z \in \mathcal{F} \times \mathcal{G}$ and
$h^{-1}(F \times Z) = f^{-1}(F) \cap g^{-1}(Z) = E \cap X = E$. This shows that $h$ is exactly measurable, and the same clearly holds
when $g$ is exactly measurable. \eop
 
The final statement in Proposition~\ref{e-meas-maps}.7 is analogous to the following simple fact:
Let $f : X \to Y$ and $g : X \to Z$ be any mappings and again define $h : X \to Y \times Z$ by $h(x) = (f(x),g(x))$ for all $x \in X$.
Then $h$ is injective provided at least one of $f$ and $g$ is injective.

Let $S$ be a non-empty set. For each $s \in S$ let $(X_s,\mathcal{E}_s)$ and $(Y_s,\mathcal{F}_s)$ 
be measurable spaces and let
$X = \prod_{s \in S} X_s$, $\mathcal{E} = \prod_{s \in S} \mathcal{E}_s$, $Y = \prod_{s \in S} Y_s$ and
$\mathcal{F} = \prod_{s \in S} \mathcal{F}_s$.
Also for each $s \in S$ let $f_s : X_s \to Y_s$ be a mapping; then there is a mapping
$f : X \to Y$ given by $f(\{x_s\}_{s \in S}) = \{f_s(x_s)\}_{s \in S}$ 
for all $\{x_s\}_{s \in S}) \in X$.

\begin{proposition}
If the mapping $f_s : (X_s,\mathcal{E}_s) \to (Y_s,\mathcal{F}_s)$ is measurable (resp.\ exactly measurable) for each $s \in S$ then
$f : (X,\mathcal{E}) \to (Y,\mathcal{F})$ is measurable (resp.\ exactly measurable).
\end{proposition}

\proof 
Suppose first that the mappings $f_s$, $s \in S$, are measurable.
Let $\mathcal{R}_X$ (resp.\ $\mathcal{R}_Y$) be the measurable rectangles in $X$ (resp.\ in $Y$). 
If $R = \prod_{s \in S} F_s \in \mathcal{R}_Y$ then
$f^{-1}(R) = \prod_{s \in S} f_s^{-1}(F_s) \in \mathcal{R}_X$ 
and so $f^{-1}(\mathcal{R}_Y) \subset \mathcal{R}_X$. 
Hence 
\[f^{-1}(\mathcal{F}) = f^{-1}(\sigma(\mathcal{R}_Y)) = \sigma(f^{-1}(\mathcal{R}_Y)) \subset 
\sigma(\mathcal{R}_X) = \mathcal{E}\;,\]
which shows that $f : (X,\mathcal{E}) \to (Y,\mathcal{F})$ is measurable.
Now suppose that the mappings $f_s$, $s \in S$, are exactly measurable.
Let $R = \prod_{s \in S} E_s \in \mathcal{R}_X$; then, since $f_s^{-1}(\mathcal{F}_s) = \mathcal{E}_s$, there exists
$F_s \in \mathcal{F}_s$ with $f_s^{-1}(F_s) = E_s$ and, since
$f_s^{-1}(Y_s) = X_s$, we can choose $F_s = Y_s$ whenever $E_s = X_s$. Thus $R' = \prod_{s \in S} F_s \in \mathcal{R}_Y$
and $f^{-1}(R') =  \prod_{s \in S} f_s^{-1}(F_s) = \prod_{s \in S} E_s = R$ 
and this shows that $f^{-1}(\mathcal{R}_Y) = \mathcal{R}_X$ (since in the first part we established that
$f^{-1}(\mathcal{R}_Y) \subset \mathcal{R}_X$.).
Hence 
\[ f^{-1}(\mathcal{F}) = f^{-1}(\sigma(\mathcal{R}_Y)) 
= \sigma(f^{-1}(\mathcal{R}_Y)) = \sigma(\mathcal{R}_X) = \mathcal{E} \]
and therefore
$f : (X,\mathcal{E}) \to (Y,\mathcal{F})$ is exactly measurable. \eop

Besides the product of measurable spaces there is the dual concept of a disjoint union.
Let $S$ be a non-empty set 
and for each $s \in S$ let $(X_s,\mathcal{E}_s)$ be a measurable space; 
assume that the sets $X_s$, $s \in S$, are disjoint and put $X = \bigcup_{s\in S} X_s$. Let 
\[\bigcup_{s \in S} \mathcal{E}_s = \{ E \subset X : E \cap X_s \in \mathcal{E}_s\ \mbox{for all}\ s \in S \}\;; \]
then $\mathcal{E} = \bigcup_{s \in S} \mathcal{E}_s$
is clearly a $\sigma$-algebra and the measurable space $(X,\mathcal{E})$ is called
the \index{disjoint union}\definition{disjoint union} of the measurable spaces $(X_s,\mathcal{E}_s)$, $s \in S$.

Suppose now $S$ is countable (i.e., finite or countable infinite). 
Let $(X,\mathcal{E})$ be the disjoint union of the measurable spaces $(X_s,\mathcal{E}_s)$, $s \in S$, let
$(Y,\mathcal{F})$ and $(Z,\mathcal{G})$ be measurable spaces with $(Z,\mathcal{G})$ separable.
Also for each $s \in S$ let $f_s : X_s \to Y$ be a mapping and $\gamma : S \to Z$ be an injective mapping. 
Define $f : X \to Y \times Z$ by 
\[f(x) = (f_s(x),\gamma(s))\]
for all $x \in X_s$, $s \in S$.
(The simplest case here is with $(Z,\mathcal{G}) = (S,\mathcal{P}(S))$ and with $\gamma$ the identity mapping.)

\begin{proposition}
If the mapping $f_s : (X_s,\mathcal{E}_s) \to (Y,\mathcal{F})$ is measurable (resp.\ exactly measurable) for each $s \in S$ then
$f : (X,\mathcal{E}) \to (Y \times Z,\mathcal{F} \times \mathcal{G})$ is measurable (resp.\ exactly measurable).
\end{proposition}

\proof
Suppose first the mappings $f_s$,
$s \in S$, are measurable. Let $F \in \mathcal{F} \times \mathcal{G}$; then the
section $F_s = \{ y \in Y : (y,\gamma(s)) \in F \}$ is an element of $\mathcal{F}$ for each $s \in S$
and $f^{-1}(F) \cap X_s = f_s^{-1}(F_s)$. Therefore
$f^{-1}(F) \cap X_s \in \mathcal{E}_s$ for each $s \in S$, which means that 
$f^{-1}(F) \in \mathcal{E}$, i.e., $f : (X,\mathcal{E}) \to (Y \times Z,\mathcal{F} \times \mathcal{G})$ is measurable.
Now suppose that the mappings $f_s$, $s \in S$, are exactly measurable. Let $E \in \mathcal{E}$; then $E \cap X_s \in \mathcal{E}_s$
for each $s \in S$ and so there exists $F_s \in \mathcal{F}$ such that $f_s^{-1}(F_s) = E \cap X_s$.
Hence $F_s \times \{\gamma(s)\}  \in \mathcal{F} \times \mathcal{G}$ for each $s \in S$ and so
$F = \bigcup_{s \in S} (F_s \times \{\gamma(s)\}) \in \mathcal{F} \times \mathcal{G}$, since $S$ is countable.
But $f^{-1}(F) = E$, since $f^{-1}(F) \cap X_s = f_s^{-1}(F_s) = E \cap X_s$ for each $s \in S$, and this implies that 
$f : (X,\mathcal{E}) \to (Y \times Z,\mathcal{F} \times \mathcal{G})$ is exactly measurable. \eop


\startsection{Countably generated measurable spaces}
\label{cg-meas-spaces}

A $\sigma$-algebra $\mathcal{E}$ is \definition{countably generated} if 
$\mathcal{E} = \sigma(\mathcal{S})$ for some countable subset $\mathcal{S}$ of $\mathcal{E}$
and  a measurable space $(X,\mathcal{E})$ is then  \definition{countably generated} if $\mathcal{E}$ is.
Note that by Lemma~\ref{e-meas-maps}.4
a countably generated measurable space $(X,\mathcal{E})$ is separated  if and only if it is countably separated, where
a measurable space $(Y,\mathcal{F})$ is said to be \definition{countably separated} if there exists a countable 
subset $\mathcal{T}$ of $\mathcal{F}$ which separates the points of $Y$.

Countably generated measurable spaces often occur as follows:

\begin{proposition}
Let $X$ be a topological space having a countable base for its topology and let 
$\mathcal{B}_X$ be the $\sigma$-algebra of Borel subsets of $X$. Then $(X,\mathcal{B}_X)$ is countably generated.
\end{proposition}

\proof 
If $\mathcal{O}_X$ is the set of open subsets of $X$ and
$\mathcal{U}$ is a countable base for the topology then each $U \in \mathcal{O}_X$ can be written as a countable
union of elements from $\mathcal{U}$ and thus $\mathcal{O}_X \subset  \sigma(\mathcal{U})$. Hence
$\mathcal{B}_X = \sigma(\mathcal{O}_X) \subset \sigma(\mathcal{U})$ and so
$\mathcal{B}_X = \sigma(\mathcal{U})$. Therefore $(X,\mathcal{B}_X)$ is countably generated. \eop

A topological space is \definition{separable} 
if it possesses a countable dense set, and it is easy to see that a metric space is separable if and only
its topology has a countable base. Thus if $X$ is a separable metric space then
by Proposition~\ref{cg-meas-spaces}.1 the measurable space $(X,\mathcal{B}_X)$ is countably generated.

It is important to note that
in general a sub-$\sigma$-algebra of a countably generated $\sigma$-algebra will 
not be countably generated.

As in the Introduction let $\specM = \{0,1\}^{\Nat}$, considered as a topological space
as the product of $\Nat$ copies of $\{0,1\}$ (with the discrete topology); $\mathcal{B}$ will always denote  the 
$\sigma$-algebra of Borel subsets of $\specM$.
By Proposition~\ref{cg-meas-spaces}.1 is $(\specM,\mathcal{B})$ is countably generated (since a compact metric space is separable). 
Moreover $(\specM,\mathcal{B})$ is separable (as a measurable space), since $\{z\}$ is a closed subset of $\specM$ for each $z \in \specM$.

The following sets (and their denotation) will be used throughout these notes. 
For $m \in \Nat$ and $z_0,\,\ldots,\,z_m \in \{0,1\}$ let 
\[  \specM(z_0,\ldots,z_m) = \{ \{z'_n\}_{n \ge 0} \in \specM :  z'_j = z_j\ \mbox{for}\ j = 0,\,\ldots,\,m \}\;, \]
denote the set of all such subsets of by $\mathcal{C}^o_{\specM}$ and let
$\mathcal{C}_{\specM}$ be the set of  all subsets of $\specM$ which can be written as a finite union of elements from
$\mathcal{C}^o_{\specM}$.
Then $\mathcal{C}^o_{\specM} \subset \mathcal{C}_{\specM}$, 
each element of the countable set $\mathcal{C}_{\specM}$
is both open and closed and both $\mathcal{C}^o_{\specM}$ and $\mathcal{C}_{\specM}$ 
are  bases for the topology on $\specM$.
Moreover, $\mathcal{C}_{\specM}$ is an algebra, which is known as the algebra of
\definition{cylinder sets} in $\specM$. In particular, since
$\mathcal{C}^o_{\specM}$ and $\mathcal{C}_{\specM}$ are countable bases for the topology it follows that
$\mathcal{B} = \sigma(\mathcal{C}_{\specM}) = \sigma(\mathcal{C}^o_{\specM})$.

The next result appears as Theorem~2.1 in Mackey \cite{mackey}.

\begin{proposition}
A measurable space $(X,\mathcal{E})$ is countably generated if and only if
there exists an exactly measurable mapping $f : (X,\mathcal{E}) \to (\specM,\mathcal{B})$.
\end{proposition}

\proof 
For each $m \ge 0$ let
$\Lambda_m = \{ \{z_n\}_{n\ge 0} \in \specM : z_m = 1 \}$. Then $\Lambda_m \in \mathcal{C}_{\specM}$ for each 
$m \ge 0$ and, on the other hand, each element of $\mathcal{C}_{\specM}$ can written as a finite intersection of elements from
the set $\{ \Lambda_m : m \ge 0 \} \cup \{ X \setminus \Lambda_m : m \ge 0 \}$. Hence
$\mathcal{B} = \sigma(\mathcal{C}_{\specM}) = \sigma(\{ \Lambda_m : m \ge 0 \})$.

Suppose now that $(X,\mathcal{E})$ is countably generated; then there exists a sequence $\{E_n\}_{n\ge 0}$
from $\mathcal{E}$ such that $\mathcal{E} = \sigma(\{ E_n : n \ge 0 \})$.
Define a mapping $f : X \to \specM$ by $f(x) = \{ I_{E_n}(x) \}_{n\ge 0}$.
Then $f^{-1}(\Lambda_n) = E_n$ for each $n \ge 0$ and therefore 
\[f^{-1}(\mathcal{B}) = f^{-1}(\sigma(\{ \Lambda_n : n\ge 0 \})) 
= \sigma (\{ E_n : n \ge 0 \}) = \mathcal{E}\;,\]
and thus $f$ is exactly measurable.
Suppose conversely there exists an exactly measurable mapping
$f : (X,\mathcal{E}) \to (\specM,\mathcal{B})$ and for $n \ge 0$ put $E_n = f^{-1}(\Lambda_n)$. 
Then 
\[
\sigma (\{ E_n : n \ge 0 \}) = \sigma(\{ f^{-1}(\Lambda_n) : n \ge 0 \})
= f^{-1}(\sigma(\{ \Lambda_n : n\ge 0 \})) = f^{-1}(\mathcal{B}) = \mathcal{E}
\]
and thus $\mathcal{E}$ is countably generated. \eop

\begin{proposition}
If $(X,\mathcal{E})$ is a countably generated  measurable space then there exists a countable algebra
$\mathcal{G}$ with $\mathcal{E} = \sigma(\mathcal{G})$.
\end{proposition}

\proof
By Proposition~\ref{cg-meas-spaces}.2 there exists an exactly measurable $f : (X,\mathcal{E}) \to (\specM,\mathcal{B})$,
then $\mathcal{G} = f^{-1}(\mathcal{C}_{\specM})$ is a countable algebra and 
\[\sigma(\mathcal{G}) = \sigma(f^{-1}(\mathcal{C}_{\specM})) = f^{-1}(\sigma(\mathcal{C}_{\specM})) = \mathcal{E}\;.\ \eop\]

Here is a property of the space $\specM$ which will play a fundamental role in what follows (where countable means finite or
countably infinite):

\begin{proposition}
If $S$ is a non-empty countable set
then $\specM^S$ (with the product topology) is homeomorphic to $\specM$.
Moreover, if $h : \specM^S \to \specM$ is a homeomorphism then 
$h : (\specM^S,\mathcal{B}^S) \to (\specM,\mathcal{B})$ is exactly measurable
(with $\mathcal{B}^S$ the product $\sigma$-algebra on $\specM^S$).
\end{proposition}

\proof 
The set $S \times \Nat$ is countably infinite, so let $\varphi : \Nat \to S \times \Nat$ be a bijective mapping.
If $\{w_s\}_{s \in S} \in \specM^S$ and $s \in S$ then the element $w_s$ of
$\specM = \{0,1\}^{\Nat}$ will be denoted by $\{w_{s,n}\}_{n\ge 0}$. 
Now define a mapping $g : \specM^S \to \specM$ by letting $g(\{w_s\}_{s \in S}) = \{z_n\}_{n\ge 0}$, where
$z_n = w_{s,k}$ and $(s,k) = \varphi(n)$. Then it is easy to see that $g$ is bijective, and it is continuous, since
$\proj_n \circ g = \proj_k \circ \proj'_s$ for each $n \ge 0$, where again $(s,k) = \varphi(n)$ and
$\proj'_s : \specM^S \to \specM$ is the projection onto the $s\,$th component.
Thus $g$ is a homeomorphism, since $\specM^S$ is compact and compact subsets of $\specM$ are closed.
Now $\mathcal{B}^S$ is the $\sigma$-algebra of Borel subsets of $\specM^S$ and therefore if
$h : \specM^S \to \specM$ is a homeomorphism then 
$h^{-1}(\mathcal{B}) = \mathcal{B}^S$. 
\eop

We next look at constructions involving measurable spaces which preserve the property of being countably generated.

\begin{proposition}
(1)\enskip Let $(X,\mathcal{E})$ be a countably generated measurable space and 
$A$ be a non-empty subset of $X$. Then $(A,\mathcal{E}_{|A})$ is  countably generated
(with $\mathcal{E}_{|A}$ the trace $\sigma$-algebra).

(2)\enskip Let $(X,\mathcal{E})$ and $(Y,\mathcal{F})$  be measurable spaces with $(X,\mathcal{E})$ countably generated.
If there exists an exactly measurable mapping $g : (Y,\mathcal{F}) \to (X,\mathcal{E})$ then
$(Y,\mathcal{F})$ is countably generated.  

In (3) and (4) let $S$ be a non-empty countable set 
and for each $s \in S$ let $(X_s,\mathcal{E}_s)$ be a countably generated measurable space.

(3)\enskip 
The product measurable space $(X,\mathcal{E})$  is countably generated.

(4)\enskip 
Assume the sets $X_s$, $s \in S$, are disjoint.
Then the disjoint union measurable space $(X,\mathcal{E})$  is countably generated.
\end{proposition}

\proof
(1)\enskip
Proposition~\ref{cg-meas-spaces}.2 implies there exists an exactly measurable mapping
$f : (X,\mathcal{E}) \to (\specM,\mathcal{B})$ and then by Proposition~\ref{e-meas-maps}.1 
$f_{|A} : (A,\mathcal{E}_{|A}) \to (\specM,\mathcal{B})$ is exactly measurable (with $f_{|A}$ the restriction of $f$ to $A$).
Thus by Proposition~\ref{cg-meas-spaces}.2 $(A,\mathcal{E}_{|A})$ is  countably generated.

(2)\enskip
By Proposition~\ref{cg-meas-spaces}.2 there exists an exactly measurable
$f : (X,\mathcal{E}) \to (\specM,\mathcal{B})$ and then by Lemma~\ref{e-meas-maps}.1~(2)
$h = f \circ g : (Y,\mathcal{F}) \to (\specM,\mathcal{B})$ is exactly measurable.
Thus by Proposition~\ref{cg-meas-spaces}.2  $(Y,\mathcal{F})$ is countably generated.

As stated above, in (3) and (4)
$S$ is a non-empty countable set 
and $(X_s,\mathcal{E}_s)$ is a countably generated measurable space for each $s \in S$.
By Proposition~\ref{cg-meas-spaces}.2 there exists for each $s \in S$ an exactly measurable mapping $f_s : (X_s,\mathcal{E}_s) \to (\specM,\mathcal{B})$.

(3)\enskip
By Proposition~\ref{e-meas-maps}.8 the mapping $f : (X,\mathcal{E}) \to (\specM^S,\mathcal{B}^S)$ is exactly measurable, where
$f(\{x_s\}_{s \in S}) = \{f_s(x_s)\}_{s \in S}$ for each $\{x_s\}_{s \in S} \in X$.
Now by Proposition~\ref{cg-meas-spaces}.4 there exists a homeomorphism
$h : \specM^S \to \specM$ and then $h : (\specM^S,\mathcal{B}^S) \to (\specM,\mathcal{B})$ is exactly measurable.
Thus by Lemma~\ref{e-meas-maps}.1~(2) $g = h \circ f : (X,\mathcal{E}) \to (\specM,\mathcal{B})$ is exactly measurable and hence by
Proposition~\ref{cg-meas-spaces}.2 $(X,\mathcal{E})$ is countably generated.

(4)\enskip
Choose an injective mapping $\gamma : S \to \specM$ and define  $f : X \to \specM^2$ by 
letting $f(x) = (f_s(x),\gamma(s))$ for each $x \in X_s$, $s \in S$.
Then Proposition~\ref{e-meas-maps}.9  implies that $f : (X,\mathcal{E}) \to (\specM^2,\mathcal{B}^2)$ 
is exactly measurable.
But by Proposition~\ref{cg-meas-spaces}.4 there exists a homeomorphism
$h : \specM^2  \to \specM$ and then $h : (\specM^2,\mathcal{B}^2) \to (\specM,\mathcal{B})$ is exactly measurable.
Therefore by Lemma~\ref{e-meas-maps}.1~(2) $g = h \circ f : (X,\mathcal{E}) \to (\specM,\mathcal{B})$ is exactly measurable and so by
Proposition~\ref{cg-meas-spaces}.2 $(X,\mathcal{E})$ is countably generated. 
\eop

We now look at what are called atoms in a measurable space.
These are important when the spaces are not separable, since they are needed to formulate conditions which
correspond to being injective for separable spaces. 

Let $(X,\mathcal{E})$ be a measurable space and for each $x \in X$ let $\mathrm{a}_x$ be the intersection of
all the elements in $\mathcal{E}$ containing $x$. Thus $x \in \mathrm{a}_x$ and for all $x,\,y \in X$ either
$\mathrm{a}_x = \mathrm{a}_y$ or 
$\mathrm{a}_x$ and $\mathrm{a}_y$ are disjoint.
A subset $A$ of $X$ is an \definition{atom} of $\mathcal{E}$ if $A = \mathrm{a}_x$ for some $x \in X$
and the set of all atoms of $\mathcal{E}$ will be denoted by $\mathrm{A}(\mathcal{E})$.
Thus $\mathrm{A}(\mathcal{E})$ defines a partition of $X$: For each $x \in X$ there is a unique
atom $A \in \mathrm{A}(\mathcal{E})$ with $x \in A$.
Of course, $(X,\mathcal{E})$ is separated if and only if $\mathrm{a}_x = \{x\}$ for each $x \in X$.
In general atoms need not be measurable, i.e., it will not always be the case that 
$\mathrm{A}(\mathcal{E}) \subset \mathcal{E}$.
However, this problem does not arise if $(X,\mathcal{E})$ countably generated:

\begin{lemma}
If $(X,\mathcal{E})$ is countably generated then $\mathrm{A}(\mathcal{E}) \subset \mathcal{E}$. More precisely,
if $f : (X,\mathcal{E}) \to (\specM,\mathcal{B})$  is exactly measurable then
$\mathrm{A}(\mathcal{E}) = \{ f^{-1}(\{z\}) : z \in f(X) \}$.
\end{lemma}

\proof 
Let $z \in f(X)$, put $A = f^{-1}(\{z\})$ and consider $x \in A$ (and so $f(x) = z$). If $E \in \mathcal{E}$ with $x \in E$ then there exists 
$B \in \mathcal{B}$ with $f^{-1}(B) = E$ and therefore $z = f(x) \in B$. Hence
$A = f^{-1}(\{z\}) \subset f^{-1}(B) = E$, and since $A \in \mathcal{E}$ this shows that
$A = \mathrm{a}_x$ for all $x \in A$.
Thus $\mathrm{A}(\mathcal{E}) = \{ f^{-1}(\{z\}) : z \in f(X) \}$, since clearly
$x \in f^{-1}(\{f(x)\})$ for each $x \in X$. \eop

Lemma~\ref{cg-meas-spaces}.1 implies that a countably generated separated measurable space is separable.

If $(X,\mathcal{E})$ and $(Y,\mathcal{F})$ are countably generated measurable spaces
then we say that a  measurable mapping $f : (X,\mathcal{E}) \to (Y,\mathcal{F})$ \definition{respects atoms} if
$f(A)$ is an atom of $\mathcal{F}$ for each atom $A \in \mathrm{A}(\mathcal{E})$ and that $f$ is \definition{injective on atoms}
if $f^{-1}(A)$ is either empty or an element of $\mathrm{A}(\mathcal{E})$ for each atom $A \in \mathrm{A}(\mathcal{F})$.

Note that if $(X,\mathcal{E})$ and $(Y,\mathcal{F})$ are also separable then $f$ always respects atoms
and $f$ is injective on atoms if and only if it is injective.

The measurable spaces occurring in the rest of this section are always assumed to be countably generated.
The next result shows that a measurable mapping which is both surjective and injective on atoms automatically respects atoms.
(For a mapping which is not surjective this need not be the case.)

\begin{proposition}
Let $f : (X,\mathcal{E}) \to (Y,\mathcal{F})$ be injective on atoms. Then:

(1)\enskip There exists an injective mapping 
$f_{\mathrm{A}} : \mathrm{A}(\mathcal{E}) \to \mathrm{A}(\mathcal{F})$ with
$f^{-1}(f_{\mathrm{A}}(A)) = A$ 
and $f(A) =  f_{\mathrm{A}}(A) \cap f(X)$ for each $A \in \mathrm{A}(\mathcal{E})$.

(2)\enskip $f^{-1}(f(E)) = E$ holds for all $E \in \mathcal{E}$.

(3)\enskip The sets $f(E_1)$ and $f(E_2)$ are disjoint whenever $E_1$ and $E_2$ are disjoint sets in $\mathcal{E}$.

(4)\enskip
If $f$ is also surjective then $f$ respects atoms.
\end{proposition}

\proof 
(1)\enskip Let $A \in \mathrm{A}(\mathcal{E})$; if $x \in A$ and $\mathrm{a}_y$ is the atom of $\mathcal{F}$ containing $y = f(x)$
then $f^{-1}(\mathrm{a}_y)$ is an atom of $\mathcal{E}$ containing $x$ and so $f^{-1}(\mathrm{a}_y) = A$. 
In particular, the atom $\mathrm{a}_y$ does not depend on the choice of $x \in A$
(since if $\mathrm{a}_{y} \cap \mathrm{a}_{y'} = \varnothing$ then
$f^{-1}(\mathrm{a}_{y}) \cap f^{-1}(\mathrm{a}_{y'}) = \varnothing$).
Hence there is a mapping $f_{\mathrm{A}} : \mathrm{A}(\mathcal{E}) \to \mathrm{A}(\mathcal{F})$ such that
$f^{-1}(f_{\mathrm{A}}(A)) = A$ for each $A \in \mathrm{A}(\mathcal{E})$. This implies $f_{\mathrm{A}}$ is injective, and
by Lemma~\ref{e-meas-maps}.2~(2) 
$f(A) = f(f^{-1}(f_{\mathrm{A}}(A)) =  f_{\mathrm{A}}(A) \cap f(X)$ for each $A \in \mathrm{A}(\mathcal{E})$.

(2)\enskip
Let $E \in \mathcal{E}$ and consider $x \in f^{-1}(f(E))$ with $f(x) = y$; then $y \in f(E)$
and so there also exists $x' \in E$ with $f(x') = y$. 
Let $A$ be the atom of $\mathcal{E}$ containing $x'$; then $y \in f(A) = f_{\mathrm{A}}(A) \cap f(X)$ and so
$y \in f_{\mathrm{A}}(A)$. Thus $x \in f^{-1}(f_{\mathrm{A}}(A)) = A$, which implies that $x$ and $x'$ lie in the same atom of
$\mathcal{E}$. Therefore $x \in E$, since $x' \in E$ and 
this shows $f^{-1}(f(E)) \subset E$. Hence
$f^{-1}(f(E)) = E$, since $E \subset f^{-1}(f(E))$ holds trivially. 

(3)\enskip 
If $E_1,\, E_2 \in \mathcal{E}$ are disjoint then by (2) the sets
$f^{-1}(f(E_1))$ and $f^{-1}(f(E_2))$ are disjoint, and hence $f(E_1)$ and $f(E_2)$ are disjoint.

(4)\enskip
If $f$ is surjective then
$f(A) = f_{\mathrm{A}}(A)$ for each $A \in \mathrm{A}(\mathcal{E})$, which means that $f$ respects atoms.
\eop

\begin{proposition}
Each exactly measurable mapping is injective on atoms.
\end{proposition}

\proof
Let $f : (X,\mathcal{E}) \to (Y,\mathcal{F})$ be exactly measurable. 
By Proposition~\ref{cg-meas-spaces}.2 there exists an exactly measurable $g : (Y,\mathcal{F}) \to (\specM,\mathcal{B})$ and
therefore by Lemma~\ref{e-meas-maps}.1~(2) $g \circ f : (X,\mathcal{E}) \to (\specM,\mathcal{B})$ is exactly measurable.
Let $A \in \mathrm{A}(\mathcal{F})$, thus by Lemma~\ref{cg-meas-spaces}.1 $A = g^{-1}(\{z\})$ for some $z \in g(Y)$, and then
$f^{-1}(A) = (g \circ f)^{-1}(\{z\})$. But if $z \notin (g \circ f)(X)$ then
$(g \circ f)^{-1}(\{z\})$ is empty and if $z \in (g \circ f)(X)$ then by Lemma~\ref{cg-meas-spaces}.1
$(g \circ f)^{-1}(\{z\})$ is an element of $\mathrm{A}(\mathcal{E})$. 
\eop

\begin{lemma}
Let $f : (X,\mathcal{E}) \to (Y,\mathcal{F})$ be a measurable mapping which respects atoms and let
$h : (Y,\mathcal{F}) \to (\specM,\mathcal{B})$ be an exactly measurable mapping. Then
$f(E) = h^{-1}(h(f(E)))$ for all $E \in \mathcal{E}$.
\end{lemma}

\proof 
For each $Z \subset Y$ let $Z_{\mathrm{a}} = \bigcup_{y \in Z} \mathrm{a}_y$, where $\mathrm{a}_y$ is the atom of $\mathcal{F}$ containing $y$.
Then by Lemma~\ref{cg-meas-spaces}.1 $h^{-1}(h(Z)) = Z_{\mathrm{a}}$ for each $Z \subset Y$. But 
if $E \in \mathcal{E}$ then $f(E) = f(E)_{\mathrm{a}}$, since $f$ respects atoms, and hence
$f(E) = h^{-1}(h(f(E))$ for all $E \in \mathcal{E}$.
\eop

\begin{lemma}
(1)\enskip If $f : (X,\mathcal{E}) \to (Y,\mathcal{F})$ and $g : (Y,\mathcal{F}) \to (Z,\mathcal{G})$ are measurable mappings which 
respect atoms then $g \circ f : (X,\mathcal{E}) \to (Z,\mathcal{G})$ also respects atoms.

(2)\enskip If $f : (X,\mathcal{E}) \to (Y,\mathcal{F})$ and $g : (Y,\mathcal{F}) \to (Z,\mathcal{G})$ are measurable mappings 
which are injective on atoms 
then $g \circ f : (X,\mathcal{E}) \to (Z,\mathcal{G})$ is also injective on atoms.
\end{lemma}

\proof 
(1)\enskip This clear.

(2)\enskip If $A$ is an atom of $\mathcal{G}$ then $g^{-1}(A)$ is either empty or an atom of $\mathcal{F}$.
But if $g^{-1}(A)$ is empty then so is $(g \circ f)^{-1}(A) = f^{-1}(g^{-1}(A))$ and if
$g^{-1}(A)$ is an atom of $\mathcal{F}$ then $(g \circ f)^{-1}(A) = f^{-1}(g^{-1}(A))$ is either empty or an atom of  $\mathcal{E}$.
\eop

Let $f : (X,\mathcal{E}) \to (\specM,\mathcal{B})$ be an exactly measurable 
and $g : (X,\mathcal{E}) \to (\specM,\mathcal{B})$ be a measurable mapping.
Put $A = f(X)$; then by Proposition~\ref{e-meas-maps}.5 there exists a unique mapping $h : A \to \specM$ such that 
$h \circ f = g$, and then $h : (A,\mathcal{F}_{|A}) \to (\specM,\mathcal{B})$ is measurable.

\begin{lemma}
If $g$ is injective on atoms then $h$ is injective. 
\end{lemma}

\proof 
Let $z \in g(X)$; then $g^{-1}(\{z\})$ is an atom of $\mathcal{E}$ and hence by Lemma~\ref{cg-meas-spaces}.1
$g^{-1}(\{z\}) = f^{-1}(\{z'\})$ for some $z' \in f(X) = A$. Thus 
\[f^{-1}(h^{-1}(\{z\}) = g^{-1}(\{z\}) = f^{-1}(\{z'\})\] 
and so by Lemma~\ref{e-meas-maps}.2~(2) $h^{-1}(\{z\}) = \{z'\}$. Therefore $h$ is injective, 
\eop


\startsection{Classifying classes}
\label{cl-classes}

A subset of $\mathcal{P}(\specM)$ will be called a \definition{classifying class} if it contains $\mathcal{B}$ and is closed under finite 
intersections. Let $\mathcal{D}$ be a classifying class; then a countably generated measurable space $(X,\mathcal{E})$
will be called a \definition{type $\mathcal{D}$ space}
if there exists an exactly measurable mapping $f : (X,\mathcal{E}) \to (\specM,\mathcal{B})$
such that $f(X) \in \mathcal{D}$. 
The cases we are mainly interested in are with $\mathcal{D} = \mathcal{B}$ since, as explained in the Introduction, separable
type $\mathcal{B}$ spaces will turn out to be the standard Borel spaces, and with $\mathcal{D} = \mathcal{A}$, the set of analytic subsets of 
$\specM$ (to be introduced in Section~\ref{anal-subsets}).
Of course, type $\mathcal{P}(\specM)$ just means countably generated.

For the whole of the section let $\mathcal{D}$ be a classifying class.
We will introduce various conditions on $\mathcal{D}$ which ensure that type $\mathcal{D}$ spaces have 
properties of a kind associated with standard Borel spaces. To be more precise, let us start by listing the more important properties of type $\mathcal{B}$
and type $\mathcal{A}$ spaces which will eventually be established in Section~\ref{type-ba}.
These are the following:

\begin{evlist}{22pt}{10pt}
\item[(1)] Type $\mathcal{B}$ spaces are closed under forming countable products and countable disjoint unions.

\item[(2)] 
If $(X,\mathcal{E})$ is a type $\mathcal{B}$ space then $f(X) \in \mathcal{B}$ for every
exactly measurable mapping $f : (X,\mathcal{E}) \to (\specM,\mathcal{B})$.

\item[(3)]
If $(X,\mathcal{E})$ is a separable type $\mathcal{B}$ and $(Y,\mathcal{F})$ a separable countably generated measurable space
then any bijective measurable mapping $f : (X,\mathcal{E}) \to (Y,\mathcal{F})$ is an isomorphism and 
$(Y,\mathcal{F})$ is a type $\mathcal{B}$ space.
\item[(4)]
Let $(X,\mathcal{E})$ be a separable type $\mathcal{B}$ and $(Y,\mathcal{F})$ a separable measurable space which is also
countably separated.
If there exists a bijective measurable mapping $f : (X,\mathcal{E}) \to (Y,\mathcal{F})$ then $\mathcal{F}$ is countable generated.
(Thus by (3) $f$ is an isomorphism and $(Y,\mathcal{F})$ is a type $\mathcal{B}$ space.)
\end{evlist}
The properties (1), (2), (3) and (4) also hold for type $\mathcal{A}$ spaces (i.e., when $\mathcal{B}$ is replaced by $\mathcal{A}$).
In addition $\mathcal{A}$ spaces have the following two properties, which do not hold in general for type $\mathcal{B}$ spaces:
\begin{evlist}{22pt}{10pt}
\item[(5)]
If $(X,\mathcal{E})$ is a type $\mathcal{A}$ and $(Y,\mathcal{F})$ a countably generated measurable space
and there exists a surjective measurable mapping $f : (X,\mathcal{E}) \to (Y,\mathcal{F})$ then
$(Y,\mathcal{F})$ is a type $\mathcal{A}$ space.
\item[(6)]
Let $(X,\mathcal{E})$ be a separable type $\mathcal{A}$ and $(Y,\mathcal{F})$ a separable measurable space which is also
countably separated. If there exists a
surjective measurable mapping $f : (X,\mathcal{E}) \to (Y,\mathcal{F})$ then $\mathcal{F}$ is countable generated.
(Thus by (5) $(Y,\mathcal{F})$ is a type $\mathcal{A}$ space.)

\end{evlist}

Properties (4) and (6) may seem somewhat technical, but, for example, (4) can often be applied in the following situation:
We have a set $X$ equipped with a `strong' and a `weak' topology; there are then the corresponding
Borel $\sigma$-algebras $\mathcal{E}$ and $\mathcal{E}'$ and, although the topologies are very different, there is still
the hope that $\mathcal{E}= \mathcal{E}'$. Now the identity mapping $\id_X : (X,\mathcal{E}) \to (X,\mathcal{E}')$
is bijective and measurable (since it is continuous as a mapping between the topological spaces)
and usually both $(X,\mathcal{E})$ and $(X,\mathcal{E}')$ will be separable (because the topologies will be Hausdorff) and $(X,\mathcal{E})$ will
be a type $\mathcal{B}$ space because the `strong' topology is given in terms of a metric.
Thus by (4) it will follow that $\mathcal{E}= \mathcal{E}'$ provided 
$(X,\mathcal{E}')$ is countably separated, which will be the case if there is a countable set of `weakly' open sets  which separate
the points of $X$. 

We begin the analysis of type $\mathcal{D}$ spaces with two simple facts which hold for all classifying classes:

\begin{lemma}
Let $(X,\mathcal{E})$
be a type $\mathcal{D}$ space and let $f : (X,\mathcal{E}) \to (\specM,\mathcal{B})$ be an exactly measurable mapping 
with $f(X) \in \mathcal{D}$. Then $f(E) \in \mathcal{D}$ for all $E \in \mathcal{E}$.
\end{lemma}

\proof
Let $E \in \mathcal{E}$ ; then by Proposition~\ref{e-meas-maps}.3 $f(E) \in \mathcal{B}_{|A}$, where $A = f(X)$, and therefore
$f(E) = B \cap f(X)$ for some $B \in \mathcal{B}$. This implies that $f(E) \in \mathcal{D}$, since $f(X) \in \mathcal{D}$. \eop

\begin{proposition}
Let $(X,\mathcal{E})$ be a type $\mathcal{D}$ and $(Y,\mathcal{F})$ be a countably generated measurable
space and suppose there exists a surjective exactly measurable mapping $f : (X,\mathcal{E}) \to (Y,\mathcal{F})$.
Then $(Y,\mathcal{F})$ is a type $\mathcal{D}$ space.
\end{proposition}

\proof 
There exists an exactly measurable mapping $g : (X,\mathcal{E}) \to (\specM,\mathcal{B})$ with $g(X) \in \mathcal{D}$.
Thus, applying Proposition~\ref{e-meas-maps}.6 (with $(Z,\mathcal{G}) = (\specM,\mathcal{B})$), 
there exists a unique mapping $h : Y \to \specM$ with $h \circ f = g$ and then 
$h : (Y,\mathcal{F}) \to (\specM,\mathcal{B})$ is exactly measurable.
Hence $h(Y) = h(f(X)) = g(X) \in \mathcal{D}$ (since $f$ is surjective) and 
this shows that $(Y,\mathcal{F})$ is a type $\mathcal{D}$ space. \eop

Next we look at conditions which ensure that type $\mathcal{D}$ spaces are closed under 
standard constructions such as forming countable products and countable disjoint unions.
The conditions which are involved here are the following:

\begin{evlist}{22pt}{10pt}
\item[(a)] If $h : \specM \to \specM$ is a homeomorphism then $h(D) \in \mathcal{D}$ for all $D \in \mathcal{D}$.

\item[(b)] Condition (a) holds and
$\mathcal{D}$ is closed under finite products in the sense that if $D_1$ and $D_2$ are elements of $\mathcal{D}$ 
and $h : \specM \times \specM \to \specM$ is a homeomorphism then $h(D_1 \times D_2) \in \mathcal{D}$.
(This is independent of which homeomorphism is used, since (a) holds.)

\item[(c)]
Condition (a) holds and
$\mathcal{D}$ is closed under countable products in the sense that if $\{D_s\}_{s \in S}$ is a non-empty countable family from
$\mathcal{D}$ and $h : \specM^S \to \specM$ is a homeomorphism then $h(\prod_{s \in S} D_s) \in \mathcal{D}$.
(Again, since (a) holds this is independent of which homeomorphism is used.)

\item[(d)] $\mathcal{D}$ is closed under countable unions.
\end{evlist}

We say that the classifying class $\mathcal{D}$ is \definition{closed under finite products} if (b) holds, and is
\definition{closed under countable products} if (c) holds.
In particular, $\mathcal{B}$ is closed under countable products and unions. (Note that (c) holds because
$\mathcal{B}^S$ is the $\sigma$-algebra of Borel subsets of $\specM^S$.)
In Section~\ref{anal-subsets} we will see that $\mathcal{A}$ is also closed under countable products and unions. 

\begin{proposition}
(1)\enskip If $(X,\mathcal{E})$ is a type $\mathcal{D}$ space 
then so is $(E,\mathcal{E}_{|E})$ for each non-empty $E \in \mathcal{E}$.

(2)\enskip Let $(X,\mathcal{E})$ be a type $\mathcal{D}$ and $(Y,\mathcal{F})$  be an arbitrary measurable space. 
If there exists an exactly measurable mapping $g : (Y,\mathcal{F}) \to (X,\mathcal{E})$ with $g(Y) \in \mathcal{E}$ 
then $(Y,\mathcal{F})$ is a type $\mathcal{D}$ space.

In (3) and (4) let $S$ be a non-empty countable set and for each $s \in S$ let $(X_s,\mathcal{E}_s)$ be a type $\mathcal{D}$  space.

(3)\enskip 
If $\mathcal{D}$ is closed under countable products then the product measurable space $(X,\mathcal{E})$ is a type $\mathcal{D}$ space.

(4)\enskip 
Assume the sets $X_s$, $s \in S$, are disjoint. If $\mathcal{D}$ is closed under countable products and unions
the disjoint union measurable space $(X,\mathcal{E})$  is a type $\mathcal{D}$ space.
\end{proposition}

\proof
(1)\enskip
There exists an exactly measurable mapping $f : (X,\mathcal{E}) \to (\specM,\mathcal{B})$ with
$f(X) \in \mathcal{D}$ and by Proposition~\ref{e-meas-maps}.1 
$f_{|E} : (E,\mathcal{E}_{|E}) \to (\specM,\mathcal{B})$ is exactly measurable.
But by Lemma~\ref{cl-classes}.1 $f_{|E}(E) = f(E) \in \mathcal{D}$ and thus 
$(E,\mathcal{E}_{|E})$ is a type $\mathcal{D}$ space.

(2)\enskip
There exists an exactly measurable
$f : (X,\mathcal{E}) \to (\specM,\mathcal{B})$ with $f(X) \in \mathcal{D}$ 
and by Lemma~\ref{e-meas-maps}.1~(2)
$h = f \circ g : (Y,\mathcal{F}) \to (\specM,\mathcal{B})$ is exactly measurable.
Moreover, by Lemma~\ref{cl-classes}.1 
$h(Y) = f(g(Y)) \in \mathcal{D}$, since $g(Y) \in \mathcal{E}$.
Thus $(Y,\mathcal{F})$ is a type $\mathcal{D}$ space.

In (3) and (4)
let $S$ be a non-empty countable set 
and for each $s \in S$ let $(X_s,\mathcal{E}_s)$ be a countably generated measurable space.
For each $s \in S$ there then exists an exactly measurable mapping $f_s : (X_s,\mathcal{E}_s) \to (\specM,\mathcal{B})$
with $f_s(X_s) \in \mathcal{D}$.

In (3) and (4) $S$ is a non-empty countable set and 
$(X_s,\mathcal{E}_s)$ is a countably generated measurable space for each $s \in S$.
For each $s \in S$ there then exists an exactly measurable mapping $f_s : (X_s,\mathcal{E}_s) \to (\specM,\mathcal{B})$
with $f_s(X_s) \in \mathcal{D}$.

(3)\enskip
As in the proof of Proposition~\ref{cg-meas-spaces}.5~(3) the mapping
$g = h \circ f : (X,\mathcal{E}) \to (\specM,\mathcal{B})$ is exactly measurable, where
$h : \specM^S \to \specM$ is a homeomorphism and the mapping $f : (X,\mathcal{E}) \to (\specM^S,\mathcal{B}^S)$ is given by
$f(\{x_s\}_{s \in S}) = \{f_s(x_s)\}_{s \in S}$ for each $\{x_s\}_{s \in S} \in X$.
But $g(X) = h(\prod_{s \in S} f_s(X_s))$
and so $g(X) \in \mathcal{D}$, since $\mathcal{D}$ is closed under countable products. This implies that
$(X,\mathcal{E})$ is a type $\mathcal{D}$ space.

(4)\enskip
As in the proof of Proposition~\ref{cg-meas-spaces}.5~(4) the mapping
$g = h \circ f : (X,\mathcal{E}) \to (\specM,\mathcal{B})$ is exactly measurable, where
$h : \specM^2  \to \specM$ is a homeomorphism and the mapping
$f : X \to \specM^2$ by is given by $f(x) = (f_s(x),\gamma(s))$ for each $x \in X_s$, $s \in S$, with
$\gamma : S \to \specM$ an arbitrary injective mapping. But
\[g(X) = h\Bigl(\bigcup_{s \in S} (f_s(X_s) \times \{\gamma(s)\})\Bigr)
= \bigcup_{s \in S} h(f_s(X_s) \times \{\gamma(s)\})\]
and so $g(X) \in \mathcal{D}$, since $\mathcal{D}$ is closed under countable products and unions. This implies 
$(X,\mathcal{E})$ is a type $\mathcal{D}$ space. \eop

We say that the classifying class $\mathcal{D}$ is \definition{invariant under isomorphisms} if whenever $A \in \mathcal{D}$
is non-empty and $h : (A,\mathcal{B}_{|A}) \to (\specM,\mathcal{B})$ is an exactly measurable mapping then $h(A) \in \mathcal{D}$.
The reason for employing this terminology is that
Lemma~\ref{e-meas-maps}.5 implies the following:
If $h : (A,\mathcal{B}_{|A}) \to (\specM,\mathcal{B})$ is exactly measurable then 
the surjective mapping
$h : (A,\mathcal{B}_{|A}) \to (h(A),\mathcal{B}_{|h(A)})$ is an isomorphism.
Conversely, if $h : (A,\mathcal{B}_{|A}) \to (B,\mathcal{B}_{|B})$ is an isomorphism 
then $h : (A,\mathcal{B}_{|A}) \to (\specM,\mathcal{B})$ is exactly measurable.
Thus $\mathcal{D}$ being invariant under isomorphisms means that whenever $A \in \mathcal{D}$
then $B \in \mathcal{D}$ for each subset $B \subset \specM$ isomorphic to $A$. 
In Section~\ref{anal-subsets} we will see that both $\mathcal{A}$ and $\mathcal{B}$ have this property..

\begin{proposition}
The following are equivalent for the classifying class $\mathcal{D}$:

(1)\enskip $\mathcal{D}$ is closed under isomorphisms.

(2)\enskip If $(X,\mathcal{E})$ is a type $\mathcal{D}$ space then $f(X) \in \mathcal{D}$ holds
for every exactly measurable mapping $f : (X,\mathcal{E}) \to (\specM,\mathcal{B})$.
\end{proposition}

\proof 
(1) $\Rightarrow$ (2):\enskip 
There exists an exactly measurable $g : (X,\mathcal{E}) \to (\specM,\mathcal{B})$ with $A = g(X) \in \mathcal{D}$;
let $f : (X,\mathcal{E}) \to (\specM,\mathcal{B})$ be any exactly measurable mapping.
Then by Proposition~\ref{e-meas-maps}.5 there is a unique mapping $h : A \to \specM$ with
$h \circ g = f$, and $h : (A,\mathcal{B}_{|A}) \to (\specM,\mathcal{B})$ is exactly measurable.
Hence $h(A) \in \mathcal{D}$, since $\mathcal{D}$ is closed under isomorphisms,
which implies that $f(X) = h(g(X)) = h(A) \in \mathcal{D}$.  

(2) $\Rightarrow$ (1):\enskip 
Let $A \in \mathcal{D}$ be non-empty. Then by Lemma~\ref{e-meas-maps}.1~(1) and Proposition~\ref{e-meas-maps}.2
the inclusion mapping $i_A :(A,\mathcal{B}_{|A}) \to (\specM,\mathcal{B})$ is exactly measurable and
thus $(A,\mathcal{B}_{|A})$ is a type $\mathcal{D}$ space, since  $i_A(A) = A \in \mathcal{D}$.
Hence if $h : (A,\mathcal{B}_{|A}) \to (\specM,\mathcal{B})$ is any exactly measurable mapping then $h(A) \in \mathcal{D}$, and this
shows that $\mathcal{D}$ is closed under isomorphisms. \eop

\begin{lemma}
Suppose $\mathcal{D}$ is closed under isomorphisms, let
$(X,\mathcal{E})$ be a type $\mathcal{D}$ space and $f,\,g : (X,\mathcal{E}) \to (\specM,\mathcal{B})$ be exactly
measurable. Then $f^{-1}(\mathcal{D}) = g^{-1}(\mathcal{D})$.
\end{lemma}

\proof 
Put $A = f(X)$, so by Proposition~\ref{cl-classes}.3 $A \in \mathcal{D}$.
By Proposition~\ref{e-meas-maps}.5 there exists a unique mapping $h : A \to \specM$ with
$h \circ f = g$, and $h : (A,\mathcal{B}_{|A}) \to (\specM,\mathcal{B})$ is exactly measurable.
Let $D \in \mathcal{D}$; then $h^{-1}(h(D \cap A)) = D \cap A$, since $h$ is injective and hence
$f^{-1}(D) = f^{-1}(D \cap A) = f^{-1}(h^{-1}(h(D \cap A))) = g^{-1}(h(D \cap A))$. 
But if $D \cap A \ne \varnothing$ then by Proposition~\ref{e-meas-maps}.1
$h_{D \cap A} : (D \cap A,\mathcal{B}_{|D \cap A}) \to (\specM,\mathcal{B})$ is exactly measurable, and $D \cap A \in \mathcal{D}$.
Therefore $h(D \cap A) \in \mathcal{D}$, since $\mathcal{D}$ is closed under isomorphisms.
This shows that $f^{-1}(D) \subset g^{-1}(\mathcal{D})$ for all $D \in \mathcal{D}$, i.e., 
$f^{-1}(\mathcal{D}) \subset g^{-1}(\mathcal{D})$, and in the same way $g^{-1}(\mathcal{D}) \subset f^{-1}(\mathcal{D})$. \eop

If $\mathcal{D}$ is closed under isomorphisms and $(X,\mathcal{E})$ is a type $\mathcal{D}$ space then
Lemma~\ref{cl-classes}.2 implies that the set $f^{-1}(\mathcal{D})$ does not depend on the exactly measurable mapping
$f : (X,\mathcal{E}) \to (\specM,\mathcal{B})$.
This set  will be denoted by $\mathcal{E}_{\mathcal{D}}$; thus $\mathcal{E} \subset \mathcal{E}_{\mathcal{D}}$ and
in particular $\mathcal{E}_{\mathcal{B}} = \mathcal{E}$.

We say that a non-empty subset $A$ of $\specM$ is \definition{regular} if every 
injective measurable mapping $h : (A,\mathcal{B}_{|A}) \to (\specM,\mathcal{B})$ is exactly measurable
(which by Lemma~\ref{e-meas-maps}.5 means that the mapping $h : (A,\mathcal{B}_{|A}) \to (h(A),\mathcal{B}_{|h(A)})$ is an isomorphism).
The empty set is also considered to be regular and the set of regular subsets of $\specM$ will be denoted by $\mathcal{R}_{\specM}$.
In Section~\ref{anal-subsets} we will see that $\mathcal{A} \subset \mathcal{R}_\specM$ and thus also
$\mathcal{B} \subset \mathcal{R}_\specM$.

Recall that if $(X,\mathcal{E})$ and $(Y,\mathcal{F})$ are countably generated measurable spaces then 
by Proposition~\ref{cg-meas-spaces}.7 any exactly measurable $f : (X,\mathcal{E}) \to (Y,\mathcal{F})$ is injective on atoms (meaning 
that $f^{-1}(A)$ is either empty or an atom of $\mathcal{E}$ for each atom $A$ of $\mathcal{F}$).

\begin{lemma}
Suppose that $\mathcal{D} \subset \mathcal{R}_\specM$,
and let $(X,\mathcal{E})$ be a type $\mathcal{D}$ and $(Y,\mathcal{F})$ be a countably generated measurable space.
Let $f : (X,\mathcal{E}) \to (Y,\mathcal{F})$ be a measurable mapping which is injective on atoms and
which respects atoms. Then $f$ is exactly measurable.
Moreover, if $\mathcal{D}$ is also closed under isomorphisms then $f(X) \in \mathcal{F}_{\mathcal{D}}$.
\end{lemma}

\proof 
Let $g : (X,\mathcal{E}) \to (\specM,\mathcal{B})$ and $h : (Y,\mathcal{F}) \to (\specM,\mathcal{B})$  be exactly measurable mappings with
$g(X) \in \mathcal{A}$. 
Now the mapping $q = h \circ f : (X,\mathcal{E}) \to (\specM,\mathcal{B})$ is measurable and by 
Proposition~\ref{cg-meas-spaces}.7 and
Lemma~\ref{cg-meas-spaces}.3~(2) it is injective on atoms. Put $A = g(X)$; by Proposition~\ref{e-meas-maps}.5 there exists a
unique mapping $p : A \to \specM$ such that $p \circ g = q$ and then $p : (A,\mathcal{B}_{|A}) \to (\specM,\mathcal{B})$ is measurable.
Moreover, by Lemma~\ref{cg-meas-spaces}.4 $p$ is injective
and therefore it is exactly measurable, since $\mathcal{D} \subset \mathcal{R}_\specM$. It thus follows from
Lemma~\ref{e-meas-maps}.5 that the  mapping $p : (A,\mathcal{B}_{|A}) \to (D,\mathcal{B}_{|D})$ is an isomorphism, where
$D = p(A) = p(g(X)) = q(X) = h(f(X)) = h(C)$.

Let $E \in \mathcal{E}$; then by Proposition~\ref{e-meas-maps}.3 $g(E) \in \mathcal{B}_{|A}$ and so $h(f(E)) = p(g(E)) \in \mathcal{B}_{|D}$.
There thus exists $B \in \mathcal{B}$ such that $h(f(E)) = B \cap D$.
But $f$ respects atoms and so by Lemma~\ref{cg-meas-spaces}.2 $C = f(X) = h^{-1}(h(f(X))) = h^{-1}(D)$
and then
\begin{eqnarray*}
f(E) &=& h^{-1}(h(f(E)))\\ &=& h^{-1}(p(g(E))) = h^{-1}(B \cap D) = h^{-1}(B) \cap h^{-1}(D) = h^{-1}(B) \cap C\;.
\end{eqnarray*}
This shows that $f(E) \in \mathcal{F}_{|C}$ for each $E \in \mathcal{E}$. Moreover, by Proposition~\ref{cg-meas-spaces}.6~(2)
$f^{-1}(f(E)) = E$ holds for all $E \in \mathcal{E}$ and therefore by
Proposition~\ref{e-meas-maps}.3 $f$ is exactly measurable. 
Finally, $f(X) = C = h^{-1}(D)$ and if $\mathcal{D}$ is closed under isomorphisms then $D = p(A) \in \mathcal{D}$. Thus in this case
$f(X) \in h^{-1}(\mathcal{D}) = \mathcal{F}_{\mathcal{D}}$.
\eop

\begin{proposition}
Suppose that $\mathcal{D} \subset \mathcal{R}_\specM$,
and let $(X,\mathcal{E})$ be a type $\mathcal{D}$ and $(Y,\mathcal{F})$ be a countably generated measurable space.
Then any surjective measurable mapping 
$f : (X,\mathcal{E}) \to (Y,\mathcal{F})$ which is injective on atoms is exactly measurable.
Moreover, if such a mapping exists then $(Y,\mathcal{F})$ is also a type $\mathcal{D}$ space.
\end{proposition}

\proof 
By Proposition~\ref{cg-meas-spaces}.6~(4) $f$ respects atoms and thus by Lemma~\ref{cl-classes}.3 $f$ is exactly measurable.
Moreover, if $f : (X,\mathcal{E}) \to (Y,\mathcal{F})$ is a surjective exactly measurable mapping then by
Proposition~\ref{cl-classes}.1 $(Y,\mathcal{F})$ is a type $\mathcal{D}$ space. \eop

\begin{proposition}
Suppose that $\mathcal{D} \subset \mathcal{R}_\specM$,
let $(X,\mathcal{E})$ be a separable type $\mathcal{D}$ and $(Y,\mathcal{F})$ be a separable countably generated measurable space.
Then any bijective measurable mapping $f : (X,\mathcal{E}) \to (Y,\mathcal{F})$ is an isomorphism
(and in this case $(Y,\mathcal{F})$ is also a type $\mathcal{D}$ space).
\end{proposition}

\proof 
This is just a special case of Proposition~\ref{cl-classes}.4. \eop

Let $\Delta : \specM \to \specM$ be the mapping given by
\[\Delta(\{z_n\}_{n\ge 0}) = \{z'_n\}_{n\ge 0}\;,\]
where $z'_n = z_{2n+1}$ for each $n \ge 0$; thus $\Delta$ is continuous and surjective.
Moreover, let $\Theta : \specM \times \specM \to \specM$ be the homeomorphism given by
\[\Theta(\{z_n\}_{n\ge 0},\{z'_n\}_{n\ge 0}) = \{w_n\}_{n\ge 0}\;,\]
where $w_{2n} = z_n$ and $w_{2n+1} = z'_n$ for each $n \ge 0$.
Note that $\Delta \circ \Theta = \pi_2$,
with $\pi_2 : \specM \times \specM \to \specM$ the projection onto the second component.

The proof of the following result is based on an idea which occurs
in Theorem~2.2 of Mackey \cite{mackey} and is also used in Lemma~4.1 in Chapter~V of 
Parthasarathy \cite{partha}.

\begin{proposition}
Let $(X,\mathcal{E})$ be a countably generated measurable space and
let $f,\, g : (X,\mathcal{E}) \to (\specM,\mathcal{B})$ be measurable mappings with 
$f$ exactly measurable.
Then there exists an exactly measurable $q : (X,\mathcal{E}) \to (\specM,\mathcal{B})$ 
with $g = \Delta \circ q$ such that $q(X)$ has the form $C \cap \Theta(A \times \specM)$, where
$C \in \mathcal{B}$ and $A = f(X)$. 
\end{proposition}

\proof 
Let $s : X \to \specM \times \specM$ be the mapping with $s(x) = (f(x),g(x))$ for all $x \in X$.
Then by Proposition~\ref{e-meas-maps}.6
$s : (X,\mathcal{E}) \to (\specM \times \specM, \mathcal{B}\times \mathcal{B})$
is exactly measurable (since $f$ is exactly measurable).
Now by Proposition~\ref{e-meas-maps}.5 there  exists a unique mapping $h : A \to \specM$ such that $h \circ f = g$, and then
$h : (A,\mathcal{B}_{|A}) \to (\specM,\mathcal{B})$ is measurable. 
Consider the measurable mapping
$r : (A \times \specM,\mathcal{B}_{|A} \times \mathcal{B})\to (\specM \times \specM,\mathcal{B}\times \mathcal{B})$ given
by letting $r(z_1,z_2) = (h(z_1),z_2)$
for all $z_1 \in A$, $z_2 \in \specM$; then
\begin{eqnarray*}
 s(X) &=& \{ (f(x),g(x)) : x \in X \} 
 = \{ (f(x),h(f(x)) : x \in X \}\\ 
&=& \{ (z,h(z)) : z \in A \} = \{ (z_1,z_2) \in A \times \specM : h(z_1) = z_2 \} \\
 &=& \{ (z_1,z_2) \in A \times \specM : r(z_1,z_2) \in D \} = r^{-1}(D)\;,
\end{eqnarray*}
where $D = \{ (z,z) : z \in \specM \}$ is the diagonal in $\specM \times \specM$. But $D$ is closed and $\mathcal{B} \times \mathcal{B}$ is the 
$\sigma$-algebra of Borel subsets of $\specM \times \specM$, hence
$D \in \mathcal{B} \times \mathcal{B}$ which implies that $s(X) \in \mathcal{B}_{|A} \times \mathcal{B}$. 
Now 
$\mathcal{B}_{|A} \times \mathcal{B} = (\mathcal{B} \times \mathcal{B})_{|A \times \specM}$
(and for those not familiar with this fact a proof is given in Lemma~\ref{cl-classes}.6 at the end of the section) 
and so there exists $C' \in \mathcal{B} \times \mathcal{B}$ such that $s(X) = C' \cap (A \times \specM)$.

Put $q = \Theta \circ s$. Then 
$q^{-1}(\mathcal{B}) = s^{-1}(\Theta^{-1}(\mathcal{B} \times \mathcal{B})) = s^{-1}(\mathcal{B} \times \mathcal{B}) = \mathcal{E}$, which means
that $q : (X,\mathcal{E}) \to (\specM,\mathcal{B})$ is exactly measurable and
\[ q(X) = \Theta(s(X)) = \Theta(C' \cap (A \times \specM)) = \Theta(C') \cap \Theta(A \times \specM) = C \cap \Theta(A \times \specM)\;, \]
where $C = \Theta(C') \in \mathcal{B}$, since $\Theta$ is a homeomorphism. Finally, 
\[g = \pi_2 \circ (f,g) = \Delta \circ \Theta \circ (f,g) = \Delta \circ \Theta \circ s = \Delta \circ q\;.\ \eop \]

\begin{proposition}
Let $\mathcal{D}$ be closed under finite products and
let $(X,\mathcal{E})$ be a type $\mathcal{D}$ and $(Y,\mathcal{F})$ a countably generated measurable space.
Let $f : (X,\mathcal{E}) \to (Y,\mathcal{F})$ be measurable
and $h : (Y,\mathcal{F}) \to (\specM,\mathcal{B})$ exactly measurable.
Then there exists an exactly measurable mapping $q : (X,\mathcal{E}) \to (\specM,\mathcal{B})$
with $q(X) \in \mathcal{D}$ such that $h \circ f = \Delta \circ q$.
\end{proposition}

\proof 
There exists an exactly measurable mapping $g : (X,\mathcal{E}) \to (\specM,\mathcal{B})$ with 
$A = g(X) \in \mathcal{D}$. Thus, applying Proposition~\ref{cl-classes}.6 to the mappings $g$ and $h \circ f$,
there exists an exactly measurable mapping $q : (X,\mathcal{E}) \to (\specM,\mathcal{B})$
such that $h \circ f = \Delta \circ q$ and such that $q(X)$ has the form $C \cap \Theta(A \times \specM)$, where
$C \in \mathcal{B}$. Therefore $q(X) \in \mathcal{D}$, since $\mathcal{D}$ is closed under finite products.
\eop

Proposition~\ref{cl-classes}.7 is useful for dealing with the case when we have some kind of inverse limit of
type $\mathcal{D}$ spaces.
Let $\mathcal{D}$ be closed under finite products, for each $n \ge 0$ let $(X_n,\mathcal{E}_n)$ be a type $\mathcal{D}$ space 
and suppose for each $n \ge 0$ there is a measurable mapping $i_n : (X_{n+1},\mathcal{E}_{n+1}) \to (X_n,\mathcal{E}_n)$. 
Then by
Proposition~\ref{cl-classes}.7 and induction it follows that
for each $n \ge 0$ there exists an exactly measurable mapping $f_n : (X_n,\mathcal{E}_n) \to (\specM,\mathcal{B})$ 
with $f_n(X_n) \in \mathcal{D}$ such that $f_n \circ i_n = \Delta \circ f_{n+1}$ for all $n \ge 0$. 
(See Proposition~\ref{kolmogorov}.2.)

We say that the classifying class $\mathcal{D}$ is \definition{closed under continuous images} if whenever $f : \specM \to \specM$
is a continuous mapping then $f(D) \in \mathcal{D}$ for each $D \in \mathcal{D}$. 
This is the first property for which $\mathcal{A}$ and $\mathcal{B}$ differ: It will follow directly from its definition that
$\mathcal{A}$ is closed under continuous images, but $\mathcal{B}$ does not have this property.

\begin{proposition}
Let $\mathcal{D}$ be closed under finite products and continuous images and
let $(X,\mathcal{E})$ be a type $\mathcal{D}$ space. Then $f(X) \in \mathcal{D}$
for every measurable mapping $f : (X,\mathcal{E}) \to (\specM,\mathcal{B})$.
In particular, Proposition~\ref{cl-classes}.3 implies that $\mathcal{D}$ is closed under isomorphisms.
\end{proposition}

\proof 
There exists an exactly measurable mapping $g : (X,\mathcal{E}) \to (\specM,\mathcal{B})$ with 
$A = g(X) \in \mathcal{D}$. Thus by Proposition~\ref{cl-classes}.6
there exists an exactly measurable mapping $q : (X,\mathcal{E}) \to (\specM,\mathcal{B})$
such that $f = \Delta \circ q$ and such that $q(X)$ has the form $C \cap \Theta(A \times \specM)$, where
$C \in \mathcal{B}$. Hence $q(X) \in \mathcal{D}$, since $\mathcal{D}$ is closed under finite products.
Therefore $f(X) = \Delta(q(X)) \in \mathcal{D}$, since $\mathcal{D}$ is also closed under continuous images.
\eop

\begin{proposition}
Let $\mathcal{D}$ be closed under finite products and continuous images,
let $(X,\mathcal{E})$ be a type $\mathcal{D}$ and $(Y,\mathcal{F})$ a countably generated measurable
space. If  there exists a surjective measurable mapping $f : (X,\mathcal{E}) \to (Y,\mathcal{F})$ then 
$(Y,\mathcal{F})$ is also a type $\mathcal{D}$ space.
\end{proposition}

\proof 
There exists an exactly measurable mapping $g : (X,\mathcal{E}) \to (\specM,\mathcal{B})$ with 
$A = g(X) \in \mathcal{D}$ and also an exactly measurable mapping $h : (Y,\mathcal{F}) \to (\specM,\mathcal{B})$.
Thus, applying Proposition~\ref{cl-classes}.6 to the mappings $g$ and $h \circ f$,
there exists an exactly measurable mapping $q : (X,\mathcal{E}) \to (\specM,\mathcal{B})$
such that $h \circ f = \Delta \circ q$ and such that $q(X)$ has the form $C \cap \Theta(A \times \specM)$, where
$C \in \mathcal{B}$. Therefore $q(X) \in \mathcal{D}$, since $\mathcal{D}$ is closed under finite products.
Hence $h(Y) = h(f(X)) = \Delta(q(X)) \in \mathcal{D}$, since $\mathcal{D}$ is also closed under continuous images and
$f$ is surjective, and this shows that $(Y,\mathcal{F})$ is a type $\mathcal{D}$ space. \eop

Let $\mathcal{D}$ be closed under finite products and continuous images,
let $(X,\mathcal{E})$ be a type $\mathcal{D}$ space and let 
$\mathcal{E}_0$ be a countably generated sub-$\sigma$-algebra of $\mathcal{E}$.
Then, since $\id_X : (X,\mathcal{E}) \to (X,\mathcal{E}_0)$ is a surjective measurable mapping, it follows from
Proposition~\ref{cl-classes}.9 that $(X,\mathcal{E}_0)$ is also a type $\mathcal{D}$ space.

Let $(X,\mathcal{E})$ be an arbitrary measurable space and let $\mathcal{E}_0$ be a sub-$\sigma$-algebra of $\mathcal{E}$.
Then for each $x \in X$ the atom of $\mathcal{E}$ containing $x$ is a subset of the atom of $\mathcal{E}_0$ containing $x$.
Thus $\mathrm{A}(\mathcal{E}_0) = \mathrm{A}(\mathcal{E})$ (i.e., $\mathcal{E}_0$ and $\mathcal{E}$ have the same atoms) if and only
if each atom of $\mathcal{E}_0$ is an atom of $\mathcal{E}$. This can also be expressed as follows: 
The points of $X$ which can be separated by an element of $\mathcal{E}$ can be separated by an element of $\mathcal{E}_0$ (i.e., if
$x_1$ and $x_2$ are such that $x_1 \in E$ and $x_2 \in X \setminus E$ for some $E \in \mathcal{E}$ then there exists
$E_0 \in \mathcal{E}_0$ such that $x_1 \in E_0$ and $x_2 \in X \setminus E_0$).

Note that $\mathrm{A}(\mathcal{E}_0) = \mathrm{A}(\mathcal{E})$ always holds when
$(X,\mathcal{E}_0)$ is separable.

\begin{lemma}
Suppose that $\mathcal{D} \subset \mathcal{R}_\specM$,
let $(X,\mathcal{E})$ be a type $\mathcal{D}$ space and  $\mathcal{E}_0$ be a countably generated sub-$\sigma$-algebra
of $\mathcal{E}$. Then $\mathcal{E}_0 = \mathcal{E}$ if and only if $\mathrm{A}(\mathcal{E}_0) = \mathrm{A}(\mathcal{E})$.
\end{lemma}

\proof 
The identity mapping $\id_X : X \to X$ results in a surjective measurable mapping
$\id_X : (X,\mathcal{E}) \to (X,\mathcal{E}_0)$ and $(X,\mathcal{E}_0)$ is countably generated.
Suppose that $\mathrm{A}(\mathcal{E}_0) = \mathrm{A}(\mathcal{E})$; then $\id_X$ is injective on atoms
and therefore by Proposition~\ref{cl-classes}.4 $\id_X$ is exactly measurable. But this just means that
$\mathcal{E}_0 = \mathcal{E}$. The converse holds, of course, trivially.
\eop

We say that a measurable space $(Y,\mathcal{F})$ is \definition{quasi-countably separated} if $\mathcal{F}$ contains
a countably generated sub-$\sigma$-algebra $\mathcal{F}_0$ such that $\mathrm{A}(\mathcal{F}_0) = \mathrm{A}(\mathcal{F})$.
Note that by Lemma~\ref{e-meas-maps}.4 a separable measurable space is quasi-countably separated if and only if it is countably separated.

\begin{proposition}
Let $\mathcal{D}$ be closed under finite products and continuous
images and with $\mathcal{D} \subset \mathcal{R}_\specM$.
Let $(X,\mathcal{E})$ be a type $\mathcal{D}$ and 
$(Y,\mathcal{F})$ be a quasi-countably separated measurable space and
suppose that there exists a surjective measurable mapping $f : (X,\mathcal{E}) \to (Y,\mathcal{F})$. 
Then $(Y,\mathcal{F})$ is countably generated (and so by Proposition~\ref{cl-classes}.9
it is a type $\mathcal{D}$ space).
\end{proposition}

\proof
Let $\mathcal{F}_0$ be a countably generated sub-$\sigma$-algebra of $\mathcal{F}$ 
with $\mathrm{A}(\mathcal{F}_0) = \mathrm{A}(\mathcal{F})$.
Let $F \in \mathcal{F}$ and put $\mathcal{F}_1 = \sigma(\mathcal{F}_0 \cup \{F\})$. Then $\mathcal{F}_1$ is countably generated
and the mapping $f : (X,\mathcal{E}) \to (Y,\mathcal{F}_1)$ is still measurable and surjective 
and so by Proposition~\ref{cl-classes}.9 $(Y,\mathcal{F}_1)$ is a type $\mathcal{D}$ space. But
$\mathrm{A}(\mathcal{F}_0) = \mathrm{A}(\mathcal{F}_1)$ and therefore by Lemma~\ref{cl-classes}.4 
$\mathcal{F}_0 = \mathcal{F}_1$, i.e., $F \in \mathcal{F}_0$. Since this holds for all $F \in \mathcal{F}$ it follows that
$\mathcal{F}_0 = \mathcal{F}$.
\eop

We end the section by proving the result (Lemma~\ref{cl-classes}.6) used in the proof of Proposition~\ref{cl-classes}.5.
First we need the following:

\begin{lemma}
Let $A$ be a non-empty subset of a set $X$. Then
$\sigma(\mathcal{S}_{|A}) = \sigma(\mathcal{S})_{|A}$ for each $\mathcal{S} \subset \mathcal{P}(X)$.
\end{lemma}

\proof Let $\mathcal{F}$ denote the subset of $\mathcal{P}(X)$ consisting
of all sets having the form $F \cup (G \setminus A)$ with
$F \in \sigma(\mathcal{S}_{|A})$ and $G \in \sigma(\mathcal{S})$. Then it is clear that
$\mathcal{F}_{|A} = \sigma(\mathcal{S}_{|A})$ and
it is easily checked that $\mathcal{F}$ is a $\sigma$-algebra. Moreover,
$\mathcal{S} \subset \mathcal{F}$, since if $S \in \mathcal{S}$ then
$S = (S \cap A) \cup (S \setminus A)$ and
$S \cap A \in \mathcal{S}_{|A} \subset \sigma(\mathcal{S}_{|A})$.
Thus $\sigma(\mathcal{S}) \subset \mathcal{F}$, which implies that
$\sigma(\mathcal{S})_{|A} \subset \mathcal{F}_{|A} = \sigma(\mathcal{S}_{|A})$.
On the other hand, $\sigma(\mathcal{S})_{|A}$ is a $\sigma$-algebra
containing $\mathcal{S}_{|A}$, and hence also
$\sigma(\mathcal{S}_{|A}) \subset \sigma(\mathcal{S})_{|A}$. \eop 

\begin{lemma}
Let $(X,\mathcal{E})$ and $(Y,\mathcal{F})$ be measurable spaces, let $A$ be a non-empty subset of $X$ and 
$B$ a non-empty subset of $Y$. Then $\mathcal{E}_{|A} \times \mathcal{F}_{|B} = (\mathcal{E} \times \mathcal{F})_{|A \times B}$.
\end{lemma}

\proof 
Let $\mathcal{R}$ be the set of all subsets of $X \times Y$ having the form $E \times F$ with $E \in \mathcal{E}$ and $F \in \mathcal{F}$;
thus $\mathcal{E} \times \mathcal{F} = \sigma(\mathcal{R})$.
Then $\mathcal{R}_{|A\times B}$ consists of all subsets of $A \times B$ having the form 
$(E \times F) \cap (A \times B)$ with $E \in \mathcal{E}$ and $F \in \mathcal{F}$.
But $(E \times F) \cap (A \times B) = (E \cap A)\times (F \cap B)$ and hence
$\mathcal{R}_{|A\times B}$ consists of all subsets of $A \times B$ having the form 
$E' \times F' $ with $E' \in \mathcal{E}_{|A}$ and $F' \in \mathcal{F}_{|B}$.
It thus follows that $\sigma(\mathcal{R}_{|A\times B}) = \mathcal{E}_{|A} \times \mathcal{F}_{|B}$.
But by Lemma~\ref{cl-classes}.5
$\sigma(\mathcal{R}_{|A\times B}) = \sigma(\mathcal{R})_{|A\times B} = (\mathcal{E} \times \mathcal{F})_{|A \times B}$
and so $\mathcal{E}_{|A} \times \mathcal{F}_{|B} = (\mathcal{E} \times \mathcal{F})_{|A \times B}$.
\eop


\startsection{Analytic subsets of $\,\specM$}
\label{anal-subsets}

The two most important classifying classes are $\mathcal{B}$, the set of Borel subsets of $\specM$, and $\mathcal{A}$, the
set of analytic subsets of $\specM$. In the present section we establish the basic properties of this latter class.
The analytic subsets are not only important in their own right; they also play a crucial role in establishing most of the
non-elementary properties of the Borel subsets of $\specM$.

Practically all the proofs given here are based on the corresponding proofs in Chapter~8 of Cohn \cite{cohn}.

Let $\specN = \Nat^{\Nat}$ (the space of all sequences $\{m_n\}_{n\ge 0}$ of elements from $\Nat$), considered as a topological space
as the product of $\Nat$ copies of $\Nat$ (with the discrete topology). This topology is also induced by the complete metric
$d : \specN \times \specN \to \Realpos$ given by
\[ d(\{m_n\}_{n\ge 0},\{m'_n\}_{n\ge 0}) = \sum_{n\ge 0} 2^{-n}\delta'(m_n,m'_n) \]
where $\delta'(m,m) = 0$ and $\delta'(m,n) = 1$ whenever $m \ne n$. 

A subset $A$ of $\specM$ is said to be \definition{analytic} if it is either empty or there exists a continuous mapping
$\tau : \specN \to \specM$ with $\tau(\specN) = A$. The set of analytic subsets of $\specM$ will be denoted by $\mathcal{A}$.

The space $\specN$ shares the property enjoyed by the space $\specM$:

\begin{lemma}
If $S$ is a non-empty countable set then the product topological space $\specN^S$ is homeomorphic to $\specN$.
\end{lemma}

\proof 
This is the same as the proof of Proposition~\ref{cg-meas-spaces}.4. \eop

The following subsets of $\specN$ will be needed in the proof of Lemma~\ref{anal-subsets}.2 and also later.
For $n,\, m_0,\,\ldots,\,m_n \in \Nat$ let 
\[  \specN(m_0,\ldots,m_n) = \{ \{m'_p\}_{p \ge 0} \in \specN :  m'_j = m_j\ \mbox{for}\ j = 0,\,\ldots,\,n \} \]
and denote the set of all such subsets of $\specN$ by $\mathcal{C}^o_{\specN}$. Then each element of $\mathcal{C}^o_{\specN}$
is both open and closed and $\mathcal{C}^o_{\specN}$ is a base for the topology on $\specN$.
(The notation employed here corresponds to that used for the analogous subsets of $\specM$.)
Note that $\specN(m_0,\ldots,m_n) = \bigcup_{p \ge 0} \specN(m_0,\ldots,m_n,p)$ for all
$n,\, m_0,\,\ldots,\,m_n \in \Nat$.

\begin{lemma}
If $D$ is a non-empty closed subset of $\specN$ then there is a continuous mapping $f : \specN \to \specN$
with $f(\specN) = D$ such that $f(\mathsf{m}) = \mathsf{m}$ for all $\mathsf{m} \in D$.
\end{lemma}

\proof 
For each $n,\, m_0,\,\ldots,\,m_n \in \Nat$ put $D(m_0,\ldots,m_n) = D \cap \specN(m_0,\ldots,m_n)$;
thus if $\mathsf{m} = \{m_n\}_{n \ge 0} \in \specN$ then $\mathsf{m} \in D$ if and only if
$D(m_0,\ldots,m_n) \ne \varnothing$ for all $n \ge 0$, since $D$ is closed.
Note also that if $D(m_0,\ldots,m_n) = \varnothing$ for some $n \ge 0$ then
$D(m_0,\ldots,m_p) = \varnothing$ for all $p \ge n$.

For each $n,\, m_0,\,\ldots,\,m_n \in \Nat$ 
such that $D(m_0,\ldots,m_n) \ne \varnothing$ choose an element $d(m_0,\ldots,m_n) \in D(m_0,\ldots,m_n)$; also
choose an element $d \in D$.
Now define a mapping $f : \specN \to \specN$ as follows: If $\mathsf{m} \in D$ then put $f(\mathsf{m}) = \mathsf{m}$.
Thus consider $\mathsf{m} = \{m_n\}_{n \ge 0} \in \specN \setminus D$; then either
$D(m_0) = \varnothing$, in which case put $f(\mathsf{m}) = d$, or $D(m_0) \ne \varnothing$ but
$D(m_0,\ldots,m_n) = \varnothing$ for some $n \ge 1$  and in this case
put $f(\mathsf{m}) = d(m_0,\ldots,m_p)$, where $p = \max \{ n \ge 0 : D(m_0,\ldots,m_n) \ne \varnothing \}$.
Therefore $f(\specN) \subset D$ and 
$f(\mathsf{m}) = \mathsf{m}$ for all $\mathsf{m} \in D$, which implies that $f(\specN) = D$.
Moreover, $f$ is continuous, since if
$\mathsf{m} =  \{m_n\}_{n \ge 0}$ and $\mathsf{m}' =  \{m'_n\}_{n \ge 0}$ with
$f(\mathsf{m}) =  \{k_n\}_{n \ge 0}$ and $f(\mathsf{m}') =  \{k'_n\}_{n \ge 0}$, and
$m_n = m'_n$ for $n = 0,\,\ldots,\,\ell$ then also $k_n = k'_n$ for $n = 0,\,\ldots,\,\ell$.
\eop

\begin{lemma}
Let $D$ be a non-empty closed subset of $\specN$ and $g : D \to \specM$ be a continuous mapping. Then $g(D)$ is analytic.
\end{lemma}

\proof 
By Lemma~\ref{anal-subsets}.1 there is a continuous mapping $f : \specN \to \specN$ with $f(\specN) = D$ and then
$h = g \circ f : \specN \to \specM$ is continuous with $h(\specN) = g(f(\specN)) = g(D)$. Hence 
$g(D)$ is analytic.
\eop

\begin{lemma}
(1)\enskip 
For each $n \ge 0$ let $A_n \in \mathcal{A}$. Then
$\bigcup_{n \ge 0} A_n$ and $\bigcap_{n \ge 0} A_n$ are both analytic.

(2)\enskip Each open and each closed subset of $\specM$ is analytic.
\end{lemma}

\proof
(1)\enskip
Consider first the countable union. Put $A = \bigcup_{n \ge 0} A_n$; if $A  = \varnothing$ then there is nothing to prove
and so we can assume that $A_m \ne \varnothing$ for some $m \ge 0$. Thus, replacing the $A_n$'s which are empty by this $A_m$,
we can in fact assume that $A_n \ne \varnothing$ for all $n \ge 0$.
Hence for each $n \ge 0$ there is a continuous mapping $f_n : \specN \to \specM$ with $f_n(\specN) = A_n$, and so
the mapping $f : \Nat \times \specN \to \specM$ given by $f(n,\mathsf{m}) = f_n(\mathsf{m})$ is continuous and 
$f(\Nat \times \specN) = A$. Moreover, the mapping
$g : \specN \to \Nat \times \specN$ given by
$g(\{m_n\}_{n\ge 0}) = (m_0,\{m'_n\}_{n\ge 0})$, where $m'_n = m_{n+1}$, is clearly a homeomorphism
and therefore $h = f \circ g : \specN \to \specM$ is continuous and $h(\specN) = f(\Nat \times \specN) =  \bigcup_{n \ge 0} A_n$. This shows 
$A$ is analytic.

Now for the countable intersection. Put $A = \bigcap_{n \ge 0} A_n$; we can assume here that $A \ne \varnothing$, which means that
$A_n \ne \varnothing$ for each $n \ge 0$. Again for each $n \ge 0$ there is a continuous mapping $f_n : \specN \to \specM$ with $f_n(\specN) = A_n$ 
and let $f : \specN^\Nat \to \specM^\Nat$ be the continuous mapping given by $f(\{\mathsf{m}_n\}_{n\ge 0}) = \{f_n(\mathsf{m}_n)\}_{n\ge 0}$.
Let 
\[ D = \{ \{z_n\}_{n \ge 0} \in \specM^\Nat : \mbox{$z_n$ is independent of $n$} \} \]
be the diagonal in $\specM^\Nat$ and let $\delta : D \to \specM$ be defined by $\delta(\{z_n\}_{n \ge 0}) = z_0$; thus $D$ is closed
and $\delta$ is continuous. Now put $C = f^{-1}(D)$; then $C$ is closed and non-empty: Let $z \in A$; for each $n \ge 0$
there exists $\mathsf{m}_n \in \specN$ with $f_n(\mathsf{m}_n) = z$ and then $\{\mathsf{m}_n\}_{n \ge 0} \in C$.
This shows that if $z \in A$ then there exists $\{\mathsf{m}_n\}_{n \ge 0} \in C$ with
$(\delta \circ f)(\{\mathsf{m}_n\}_{n \ge 0}) = z$, and in fact
$(\delta \circ f)(C) = A$, since $(\delta \circ f)(\{\mathsf{m}_n\}_{n \ge 0}) \in A$ for all $\{\mathsf{m}_n\}_{n \ge 0} \in C$.
By Lemma~\ref{anal-subsets}.1 there exists a homeomorphism $g : \specN \to \specN^\Nat$ and then $D = g^{-1}(C)$ is a non-empty closed 
subset of $\specN$. Moreover, the mapping $q =   \delta  \circ f \circ  g : D \to \specM$ is continuous 
and $q(D) = A$ and therefore by Lemma~\ref{anal-subsets}.3 $\bigcap_{n \ge 0} A_n$ is analytic.

(2)\enskip
First note that
each element of $\mathcal{C}^o_{\specM}$ is analytic:
For each $k = 0,\,\ldots,\,m$ let $\tau_k : \Nat \to \{0,1\}$ be the constant mapping with value $z_k$ and for each
$n > m$ let $\tau_n : \Nat \to \{0,1\}$ be any surjective mapping. Then the mapping
$\tau : \specN \to \specM$ given by $\tau(\{m_n\}_{n \ge 0}) = \{\tau_n(m_n)\}_{n \ge 0}$
is continuous and $\tau(\specN) = \specM(z_0,\ldots,z_m)$, i.e., $\specM(z_0,\ldots,z_m)$ is analytic.
But $\mathcal{C}^o_{\specM}$ is a countable base for the topology on $\specM$ and so each open set can be written as a countable
union of elements from $\mathcal{C}^o_{\specM}$. Thus by (1) each open subset of $\specM$ is analytic.

Finally, each closed subset of $\specM$ can be written as a countable intersection of open sets. (This is true in any metric space.)
Therefore, by (1) each closed subset of $\specM$ is analytic.
\eop

\begin{proposition}
Each Borel subset of $\specM$ is analytic, i.e., $\mathcal{B} \subset \mathcal{A}$.
\end{proposition}

\proof
If $\mathcal{G}$ is a subset of $\mathcal{P}(\specM)$ containing the open and closed sets and which is closed under countable 
intersections and countable unions then $\mathcal{B} \subset \mathcal{G}$. (A proof of this standard fact can be found, for example
in Cohn \cite{cohn}, Lemma~8.2.4.) Therefore by Lemma~\ref{anal-subsets}.4 $\mathcal{B} \subset \mathcal{A}$.
\eop

There exist analytic subsets of $\specM$ which are not Borel, see, for example, Cohn \cite{cohn}, Corollary~8.2.17 and Exercise 6 which follows it.

\begin{lemma}
(1)\enskip If $f : \specM \to \specM$ is continuous then $f(A) \in \mathcal{A}$ for all $A \in \mathcal{A}$.

(2)\enskip Let $S$ be a non-empty countable set, for each $s \in S$ let $A_s \in \mathcal{A}$ and let $f : \specM^S \to \specM$
be a continuous mapping. Then $f(\prod_{s \in S} A_s)$ is analytic.
\end{lemma}

\proof 
(1)\enskip
There exists a continuous mapping $g : \specN \to \specM$ with $g(\specN) = A$ and then
$h = f \circ g : \specN \to \specM$ is a continuous mapping with $h(\specN) = f(g(\specN)) = f(A)$.
Thus $f(A) \in \mathcal{A}$.

(2)\enskip
For each $s \in S$ there exists a continuous mapping
$h_s : \specN \to \specM$ such that  $h_s(\specN) = A_s$ and hence there is a continuous mapping
$h :  \specN^S \to \specM^S$ given by $h(\{z_s\}_{s \in S}) =  \{ h_s(z_s)\}_{s \in S}$
with $h(\specN^S) = \prod_{s \in S} A_s$. But by Lemma~\ref{anal-subsets}.1 there exists a homeomorphism
$g : \specN \to \specN^S$ and therefore $s = f \circ h \circ g : \specN \to \specM$ is a continuous mapping with
$s(\specN) = f(\prod_{s \in S} A_s)$. 
This shows that $f(\prod_{s \in S} A_s)$ is analytic. \eop

\begin{proposition}
The set of analytic subsets $\mathcal{A}$ is a classifying class which is closed under countable products, countable unions and under
continuous images.
\end{proposition}

\proof 
This follows from Lemmas \ref{anal-subsets}.4 and \ref{anal-subsets}.5 and Proposition~\ref{anal-subsets}.1.
\eop

Recall that $\mathcal{B}$ is a classifying class which is closed under countable products and countable unions. However, it
is not closed under continuous images.

\begin{proposition}
Let $A \in \mathcal{A}$ be non-empty and let $f : (A,\mathcal{B}_{|A}) \to (\specM,\mathcal{B})$ be any measurable mapping.
Then $f(A) \in \mathcal{A}$. In particular, this shows (together with Proposition~\ref{cl-classes}.3) that $\mathcal{A}$ is closed under
isomorphisms.
\end{proposition}

\proof
This is really just a special case of Proposition~\ref{cl-classes}.7: By Proposition~\ref{anal-subsets}.2 $\mathcal{A}$ is closed under
finite products and continuous images and 
$(A,\mathcal{B}_{|A})$ is a type $\mathcal{A}$ space, since by Lemma~\ref{e-meas-maps}.1~(1) and Proposition~\ref{e-meas-maps}.2 the inclusion mapping
$i_A : (A,\mathcal{B}_{|A}) \to (\specM,\mathcal{B})$ is exactly measurable and $i_A(A) = A \in \mathcal{A}$. 
Thus $f(A) \in \mathcal{A}$ for any measurable mapping $f : (A,\mathcal{B}_{|A}) \to (\specM,\mathcal{B})$. 
\eop

Let $A \in \mathcal{A}$ be non-empty and $f : (A,\mathcal{B}_{|A}) \to (\specM,\mathcal{B})$ be a measurable mapping.
Then Proposition~\ref{anal-subsets}.3 implies that $f(A \cap C) \in \mathcal{A}$ for all $C \in \mathcal{A}$. (This holds trivially
if $A \cap C = \varnothing$. If $A \cap C \ne \varnothing$ then, since $A \cap C \in \mathcal{A}$,  
Proposition~\ref{anal-subsets}.3 applied to the
measurable mapping $f_{|A \cap C} : (A \cap C,\mathcal{B}_{|A \cap C}) \to (\specM,\mathcal{B})$ 
shows that $f(A \cap C)  = f_{|A \cap C}(A \cap C) \in \mathcal{A}$.)

The following result is known as the Separation Theorem:

\begin{proposition}
If $A_1$ and $A_2$ are disjoint analytic subsets of $\specM$ then there exist disjoint Borel sets $B_1$ and $B_2$ with
$A_1 \subset B_1$ and $A_2 \subset B_2$. 
\end{proposition}

\proof
By definition disjoint subsets $C_1$ and $C_2$ of $\specM$ can be \definition{separated by Borel sets}
if there exist disjoint $B_1,\,B_2 \in \mathcal{B}$ with $C_1 \subset B_1$ and $C_2 \subset B_2$. 

\begin{lemma}
Let $\{C_n\}_{n \ge 0}$ and $\{D_n\}_{n \ge 0}$ be two sequences of subsets of $\specM$ such that for all $m,\,n \ge 0$ the sets
$C_m$ and $D_n$ can be separated by Borel sets. Then 
$\bigcup_{n \ge 0} C_n$ and $\bigcup_{n \ge 0} D_n$ can be separated by Borel sets.
\end{lemma}

\proof
For all $m,\,n \ge 0$ there exist disjoint sets $B^C_{m,n},\, B^D_{m,n} \in \mathcal{B}$ with
$C_m \subset  B^C_{m,n}$ and $D_n \subset B^D_{m,n}$. Now for each $m \ge 0$ the
sets $\bigcup_{n \ge 0} B^C_{m,n}$ and $\bigcap_{n \ge 0} B^D_{m,n}$ are disjoint and contain $\bigcup_{n \ge 0} C_n$ and $D_m$ respectively
and from this follows that
$\bigcap_{m \ge 0 }\bigcup_{n \ge 0} B^C_{m,n}$ and $\bigcup_{m \ge 0}\bigcap_{n \ge 0} B^D_{m,n}$ are disjoint Borel sets containing 
$\bigcup_{n \ge 0} C_n$ and $\bigcup_{m\ge 0} D_m$ respectively
\eop

Now to the proof of Proposition~\ref{anal-subsets}.4. 
We can assume that both of $A_1$ and $A_2$ are non-empty
(since otherwise either $\varnothing$ and $\specM$ or $\specM$ and $\varnothing$ separate $A_1$ and $A_2$). 
There thus exist continuous mappings $f_1,\,f_2 : \specN \to \specM$ with $f_1(\specN) = A_1$ and $f_2(\specN) = A_2$.
Since $\specN(m_0,\ldots,m_n) = \bigcup_{p \ge 0} \specN(m_0,\ldots,m_n,p)$ 
it follows that
\[h(\specN(m_0,\ldots,m_n)) = \bigcup_{p \ge 0} h(\specN(m_0,\ldots,m_n,p))\] 
for each mapping $h : \specN \to \specM$ and all $n,\, m_0,\,\ldots,\,m_n \in \Nat$.
Thus if $f_1(\specN) = A_1$ and $f_2(\specN) = A_2$ cannot be separated by Borel sets then by Lemma~\ref{anal-subsets}.6
there exist elements $\mathsf{k} = \{k_n\}_{n\ge 0}$ and $\mathsf{m} = \{m_n\}_{n\ge 0}$ from $\specN$ such that for each $n \ge 0$
the sets $f_1(\specN(k_0,\ldots,k_n))$ and $f_2(\specN(m_0,\ldots,m_n))$ cannot be separated by Borel sets.
Now $f_1(\mathsf{k}) \in A_1$ and $f_2(\mathsf{m}) \in A_2$ and $A_1 \cap A_2 = \varnothing$ and so $f_1(\mathsf{k}) \ne f_2(\mathsf{m})$. 
There thus exist disjoint open subsets $U_1$ and $U_2$ of $\specM$ containing $f_1(\mathsf{k})$ and $f_2(\mathsf{m})$ respectively and then
$f_1^{-1}(U_1)$ and $f_2^{-1}(U_2)$ are open sets in $\specN$ containing $\mathsf{k}$ and $\mathsf{m}$ respectively.
However, this implies  $\specN(k_0,\ldots,k_n) \subset f_1^{-1}(U_1)$ and $\specN(m_0,\ldots,m_n) \subset f_2^{-1}(U_2)$ 
and hence that
$f_1(\specN(k_0,\ldots,k_n)) \subset U_1$ and $f_2(\specN(m_0,\ldots,m_n)) \subset U_2$ 
for all large enough $n$. But $U_1$ and $U_2$ are disjoint Borel sets, and therefore
$f_1(\specN(k_0,\ldots,k_n))$ and $f_2(\specN(m_0,\ldots,m_n))$ can be  separated by Borel sets for all large enough $n$.
This contradiction shows that $A_1$ and $A_2$ can be separated by Borel sets. \eop

Proposition~\ref{anal-subsets}.4 implies that the only subsets of $\specM$ which are both analytic and have an analytic complement 
are the Borel subsets.

\begin{lemma}
If $\{A_n\}_{n \ge 0}$ is disjoint sequence of analytic subsets of $\specM$ then there exists a disjoint sequence
$\{B_n\}_{n\ge 0}$ of Borel sets such that $A_n \subset B_n$ for all $n \ge 0$.
\end{lemma}

\proof 
Let $n \ge 0$; then by Lemma~\ref{anal-subsets}.4~(1) $A'_n = \bigcup_{m\ne n} A_m$ is analytic and it is disjoint from $A_n$.
Thus by Proposition~\ref{anal-subsets}.4  there exist disjoint Borel sets $B'_n$ and $B''_n$ with
$A_n \subset B'_n$ and $A'_n \subset B''_n$, and hence also $A'_n \subset \specM \setminus B'_n$.
For each $n \ge 0$ put $B_n = B'_n \setminus \bigcup_{m\ne n} B'_m$; then 
$\{B_n\}_{n\ge 0}$ is a disjoint sequence of Borel sets with $A_n \subset B_n$ for all $n \ge 0$.
\eop

\begin{proposition}
Let $A \in \mathcal{A}$ be non-empty and let  $f : (A,\mathcal{B}_{|A}) \to (\specM,\mathcal{B})$ be an injective measurable mapping.
Then $f$ is exactly measurable. Thus $\mathcal{A} \subset \mathcal{R}_{\specM}$.
\end{proposition}

\proof
Let $B \in \mathcal{B}$ and put $B' = \specM \setminus B$;
by Proposition~\ref{anal-subsets}.3 the sets $f(B \cap A)$ and $f(B' \cap  A)$ are analytic and they are disjoint,
since $f$ is injective.  Thus by Proposition~\ref{anal-subsets}.4 there exist disjoint Borel sets $C$ and $C'$ with $f(B \cap A) \subset C$ and
$f(B' \cap A) \subset C'$, and hence also $f(B \cap A) \subset C \cap f(A)$ and $f(B' \cap A) \subset C' \cap f(A)$.
But the sets $f(B \cap A)$ and $f(B' \cap A)$ are disjoint and $f(B \cap A) \cup f(B' \cap A) = f(A)$, which implies that
$f(B \cap A) = C \cap f(A)$ and $f(B' \cap  A) = C' \cap f(A)$. In particular, $f(B \cap A) = C \cap f(A) \in \mathcal{B}_{|f(A)}$.
\eop

Since $\mathcal{B} \subset \mathcal{A}$ it follows immediately from Proposition~\ref{anal-subsets}.5 that also
$\mathcal{B} \subset \mathcal{R}_{\specM}$.

\begin{proposition}
Let $B \in \mathcal{B}$ be non-empty and let  $f : (B,\mathcal{B}_{|B}) \to (\specM,\mathcal{B})$ be an injective measurable mapping.
Then $f(B) \in \mathcal{B}$. In particular, this shows (together with Proposition~\ref{cl-classes}.3) that $\mathcal{B}$ is closed under
isomorphisms.
\end{proposition}

\proof
By Lemma~\ref{anal-subsets}.8 below there is a measurable mapping $g : (\specM,\mathcal{B}) \to (\specM,\mathcal{B})$ with $g(\specM) \subset B$ 
such that $g(f(z)) = z$ for all $z \in B$. If $z \in f(B)$ then $z = f(z')$ for some $z' \in B$ and hence
$f(g(z)) = f(g(f(z'))) = f(z') = z$. On the other hand, if $z \in \specM$ with $f(g(z)) = z$ then $z \in f(B)$, since $g(\specM) \subset B$.
This shows that $f(B) = \{ z \in \specM : g(f(z)) = z \}$. Let $h = g \circ f$; then $h : (\specM,\mathcal{B}) \to (\specM,\mathcal{B})$  
is measurable and $f(B) = \{ z \in \specM : h(z) = z \}$. But
\[  \{ z \in \specM : h(z) = z \}
   = \bigcap_{n \ge 0} \bigcup_{z_0,\ldots,z_n}  \specM(z_0,\ldots,z_n) \cap h^{-1}(\specM(z_0,\ldots,z_n)) \in \mathcal{B} \] 
and thus $f(B) \in \mathcal{B}$. 
\eop

\begin{lemma}
Let $B \in \mathcal{B}$ be non-empty and $f : (B,\mathcal{B}_{|B}) \to (\specM,\mathcal{B})$ be injective and measurable.
Then there exists a measurable mapping $g : (\specM,\mathcal{B}) \to (\specM,\mathcal{B})$ with $g(\specM) \subset B$ such that
$g(f(z)) = z$ for all $z \in B$.
\end{lemma}

\proof
Fix some element $w \in B$.
For each $n \ge 0$ and all $z_0,\,\ldots,\,z_n \in \{0,1\}$ choose an element $[z_0,\ldots,z_n]$ in
$B \cap \specM(z_0,\ldots,z_n)$ if this set is non-empty and put $[z_0,\ldots,z_n] = w$ otherwise.
Let $n \ge 0$; then by Proposition~\ref{anal-subsets}.3 the sets $f(B  \cap \specM(z_0,\ldots,z_n))$,
$z_0,\,\ldots,\,z_n \in \{0,1\}$ are analytic, and they are disjoint, since $f$ is injective.
Thus by Lemma~\ref{anal-subsets}.7 there exist disjoint elements 
$R(z_0,\ldots,z_n)$, $z_0,\,\ldots,\,z_n \in \{0,1\}$, from $\mathcal{B}$ such that
$f(B \cap  \specM(z_0,\ldots,z_n)) \subset R(z_0,\ldots,z_n)$
for all $z_0,\,\ldots,\,z_n \in \{0,1\}$.
Define a mapping $g_n : \specM \to \specM$ by letting
\[ g_n(z) = \left\{   \begin{array}{cl}
                       [z_0,\ldots,z_0] & \ \mbox{if $z \in R(z_0,\ldots,z_n)$ for some $z_0,\,\ldots,\,z_n$}\;,\\
                              w          & \ \mbox{otherwise}\;.
                      \end{array}  \right.
\]
Then $g_n : (\specM,\mathcal{B}) \to (\specM,\mathcal{B})$ is measurable, since $g_n(\specM)$ is finite and
$g_n^{-1}(\{z\}) \in \mathcal{B}$ for each $z \in \specM$, and $g_n(\specM) \subset B$. Moreover,
if $z = \{z_n\}_{n \ge 0}\in B$ then $g_n(f(z))$ and $z$ both lie in $\specM(z_0,\ldots,z_n)$.
Now define $g : \specM \to \specM$ by letting
\[ g(z) = \left\{   \begin{array}{cl}
                       \lim\limits_{n\to\infty}g_n(z) & \ \mbox{if the limit exists and lies in $B$}\;,\\
                              w          & \ \mbox{otherwise}\;.
                      \end{array}  \right.
\]
Then $g : (\specM,\mathcal{B}) \to (\specM,\mathcal{B})$ is measurable, $g(\specM) \subset B$ and 
$g(f(z)) = z$ for all $z \in B$.
\eop

To end the section we show that every uncountable element of $\mathcal{B}$ is isomorphic to $\specM$. More precisely,
if $B \in \mathcal{B}$ is uncountable then there exists an isomorphism $f : (B,\mathcal{B}_{|B}) \to (\specM,\mathcal{B})$.
Surprisingly, this powerful result is needed less often than might be expected.
We start with what can be considered as a measurable version of the Cantor-Bernstein theorem.

\begin{lemma}
Let $(X,\mathcal{E})$ be a measurable space, let $E \in \mathcal{E}$ and suppose there exists an isomorphism
$h :(X,\mathcal{E}) \to (E,\mathcal{E}_{|E})$. Then for every $F \in \mathcal{E}$ with $E \subset F$ there exists an
isomorphism $q : (F,\mathcal{E}_{|F}) \to (X,\mathcal{E})$. 
\end{lemma}

\proof
This is more-or-less identical with the usual modern proof of the Cantor-Bernstein theorem.
Let $F' = \bigcup_{n\ge 0} h^n(F \setminus E)$. Then $F' \subset F$, $F \setminus F' \subset E$ and $h(F') \subset F'$.
Now define a mapping $q : F \to X$ by
\[ q(x) = \left\{  \begin{array}{cl} 
                       x      &\ \mbox{if $x \in F'$}\;,\\
                       h^{-1}(x) &\ \mbox{if $x \in F \setminus F'$}\;.
           \end{array} \right. \]
It is then easily checked that $q$ is bijective, and in fact $q : (F,\mathcal{E}_{|F}) \to (X,\mathcal{E})$ is an isomorphism, since
$h$ is. \eop

\begin{proposition}
Let $A \in \mathcal{A}$ be uncountable. Then there exists a continuous injective mapping $f : \specM \to \specM$ with
$f(\specM) \subset A$ (which means that any uncountable analytic set has the power of the continuum).
\end{proposition}

\proof 
There exists a continuous mapping $\tau : \specN \to \specM$ with $\tau(\specN) = A$ and so it is enough
show there exists a continuous mapping $g : \specM \to \specN$ such that $f = \tau \circ g$ is injective. 
Choose a subset $D$ of $\specN$ so that the restriction $\tau_{|D}$ of $\tau$ to $D$ maps $D$ bijectively
onto $A$. (In Section~\ref{selectors} we will see that $D$ can actually be defined without using the axiom of
choice.) Let $S$ be the set of elements in $D$ which are condensation points of $D$ (so each neighbourhood of an 
element in $S$ contains uncountably many elements of $D$). Then $D\setminus S$ is countable (this is part of
the Cantor-Bendixon theorem) and hence $S$ is uncountable. The required mapping $g  : \specM \to \specN$ 
is now defined so that $g(\mathsf{m})$ is the limit of a suitably chosen sequence from $S$ for each
$\mathsf{m} \in \specM$; the details are left to the reader.
\eop

\begin{proposition}
Let $B \in \mathcal{B}$ be uncountable. Then there exists an isomorphism
$h : (B,\mathcal{B}_{|B}) \to (\specM,\mathcal{B})$.
\end{proposition}

\proof 
By Proposition~\ref{anal-subsets}.1 $B \in \mathcal{A}$ and therefore by Proposition~\ref{anal-subsets}.7 there exists a continuous
injective mapping $f : \specM \to \specM$ with $C = f(\specM) \subset B$. In particular, this means 
$f : (\specM,\mathcal{B}) \to (\specM,\mathcal{B})$ is measurable and injective, thus by Proposition~\ref{anal-subsets}.6
$C \in \mathcal{B}$ and by Proposition~\ref{anal-subsets}.5
$f : (\specM,\mathcal{B}) \to (C,\mathcal{B}_{|C})$ is an isomorphism.
Hence by Lemma~\ref{anal-subsets}.9 there exists an isomorphism
$h : (B,\mathcal{B}_{|B}) \to (\specM,\mathcal{B})$.
\eop


\startsection{Type $\mathcal{B}$ and type $\mathcal{A}$ spaces}
\label{type-ba}

Recall from Section~\ref{cl-classes} that a classifying class is a subset of $\mathcal{P}(\specM)$ containing $\mathcal{B}$ and closed under
finite intersections. The main examples of such classes are $\mathcal{B}$, the set of Borel subsets of $\specM$, and $\mathcal{A}$, the set of analytic 
subsets of $\specM$ (and note that by Lemma~\ref{anal-subsets}.4~(1) and Proposition~\ref{anal-subsets}.1 $\mathcal{A}$ is 
a classifying class).

In this section we look at the properties of type $\mathcal{B}$ and type $\mathcal{A}$ spaces. This essentially just involves
collecting together the results from Sections \ref{cl-classes} and \ref{anal-subsets}. The one new fact is the the first part of the following, 
which implies that a standard Borel space is a type $\mathcal{B}$ space.

\begin{theorem}
(1)\enskip Let $X$ be a complete separable metric space and 
$\mathcal{B}_X$ be the $\sigma$-algebra of Borel subsets of $X$. Then $(X,\mathcal{B}_X)$ is a type $\mathcal{B}$ space.
This implies that any standard Borel space is a separable type $\mathcal{B}$ space.

(2)\enskip
Each separable type $\mathcal{B}$ space is a standard Borel space.
\end{theorem}

\proof This is given at the end of the section. \eop

Let $(X,\mathcal{E})$ be a countably generated measurable space $(X,\mathcal{E})$. By definition $(X,\mathcal{E})$ 
is then a  type $\mathcal{B}$ space (resp.\ a type $\mathcal{A}$ space) if there exists an
exactly measurable mapping $f : (X,\mathcal{E}) \to (\specM,\mathcal{B})$ such that
$f(X) \in \mathcal{B}$ (resp.\ such that $f(X) \in \mathcal{A}$). In fact, this then holds for every such mapping.

\begin{proposition}
Let $(X,\mathcal{E})$ be a type $\mathcal{B}$ (resp.\ a type $\mathcal{A}$ space).
Then $f(X) \in \mathcal{B}$ (resp.\ $f(X) \in \mathcal{A}$) 
for every exactly measurable mapping $f : (X,\mathcal{E}) \to (\specM,\mathcal{B})$.
\end{proposition}

\proof 
This follows immediately from Proposition~\ref{cl-classes}.3 since by Proposition~\ref{anal-subsets}.6 (resp.\ Proposition~\ref{anal-subsets}.3)
$\mathcal{B}$ (resp.\ $\mathcal{A}$) is closed under isomorphisms.
\eop

\begin{proposition}
(1)\enskip 
If $(X,\mathcal{E})$ is a type $\mathcal{B}$ space (resp.\ a type $\mathcal{A}$ space) then
so is $(E,\mathcal{E}_{|E})$ for each non-empty $E \in \mathcal{E}$.

(2)\enskip Let $(X,\mathcal{E})$ be a type $\mathcal{B}$ (resp.\ a type $\mathcal{A}$)
and $(Y,\mathcal{F})$  be an arbitrary measurable space.
If there exists an exactly measurable mapping $g : (Y,\mathcal{F}) \to (X,\mathcal{E})$ with $g(Y) \in \mathcal{E}$ 
then $(Y,\mathcal{F})$ is a type $\mathcal{B}$ (resp.\ a type $\mathcal{A}$ space).

In (3) and (4) let $S$ be a non-empty countable set 
and for each $s \in S$ let $(X_s,\mathcal{E}_s)$ be a type $\mathcal{B}$ (resp.\ a type $\mathcal{A}$) space.

(3)\enskip 
The product measurable space is a type $\mathcal{B}$ (resp.\ a type $\mathcal{A}$ space).

(4)\enskip 
If the sets $X_s$, $s \in S$, are disjoint then the disjoint union measurable space is a type $\mathcal{B}$ (resp. a type $\mathcal{A}$) space.
\end{proposition}

\proof 
This is special case of Proposition~\ref{cl-classes}.2, since $\mathcal{B}$ and $\mathcal{A}$ are both closed and countable products
and countable unions. (For $\mathcal{B}$ this is clear and for $\mathcal{A}$ 
it follows from Proposition~\ref{anal-subsets}.2.) \eop

\begin{theorem}
Let $(X,\mathcal{E})$ be a type $\mathcal{A}$ space and $(Y,\mathcal{F})$ be countably generated.
If there exists a surjective measurable mapping $f : (X,\mathcal{E}) \to (Y,\mathcal{F})$ then $(Y,\mathcal{F})$ is a type $\mathcal{A}$
space.
\end{theorem}

\proof 
This follows from Proposition~\ref{cl-classes}.9, since by Proposition~\ref{anal-subsets}.2 $\mathcal{A}$ is closed under finite products and
continuous images. \eop

\begin{theorem}
Let $(X,\mathcal{E})$ be a type $\mathcal{B}$ space and $(Y,\mathcal{F})$ be countably generated.
If there exists a surjective exactly measurable mapping 
$f : (X,\mathcal{E}) \to (Y,\mathcal{F})$ then $(Y,\mathcal{F})$ is a type $\mathcal{B}$ space.
In particular, this is the case if $f$ is an isomorphism.
\end{theorem}

\proof 
This is a special case of Proposition~\ref{cl-classes}.1. \eop

\begin{theorem}
Let $(X,\mathcal{E})$ be a separable type $\mathcal{B}$ (resp.\ separable type $\mathcal{A}$)
and $(Y,\mathcal{F})$ be a separable countably generated measurable space.
Then any bijective measurable mapping $f : (X,\mathcal{E}) \to (Y,\mathcal{F})$ is an isomorphism
(and in this case $(Y,\mathcal{F})$ is also a type $\mathcal{B}$ (resp.\ a type $\mathcal{A}$) space).
\end{theorem}

\proof 
This follows from Proposition~\ref{cl-classes}.4, since by Proposition~\ref{anal-subsets}.5  $\mathcal{A} \subset \mathcal{R}_{\specM}$
(and thus also $\mathcal{B} \subset \mathcal{R}_{\specM}$).  \eop

\begin{theorem}
Let $f : (X,\mathcal{E}) \to (Y,\mathcal{F})$ be an injective measurable mapping with $(X,\mathcal{E})$ and $(Y,\mathcal{F})$ separable type $\mathcal{B}$ 
spaces. Then $f(E) \in \mathcal{F}$ for all $E \in \mathcal{E}$ and in particular $f(X) \in \mathcal{F}$.
\end{theorem}

\proof 
By Proposition~\ref{anal-subsets}.5  $\mathcal{B} \subset \mathcal{R}_{\specM}$ and by Proposition~\ref{anal-subsets}.6
$\mathcal{B}$ is closed under isomorphisms; moreover, $\mathcal{F}_\mathcal{B} = \mathcal{B}$.
Thus by Lemma~\ref{cl-classes}.3 $f$ is exactly measurable and $A = f(X) \in \mathcal{F}$, and so by Proposition~\ref{e-meas-maps}.3
$f(E) \in \mathcal{F}_{|A} \subset \mathcal{F}$ for all $E \in \mathcal{E}$.
\eop

\begin{theorem}
Let $(X,\mathcal{E})$ be a separable type $\mathcal{A}$ and $(Y,\mathcal{F})$ a separable measurable space which is countably
separated. If there exists a
surjective measurable mapping $f : (X,\mathcal{E}) \to (Y,\mathcal{F})$ then $\mathcal{F}$ is countably generated
(and thus by Theorem~\ref{type-ba}.3 $(Y,\mathcal{F})$ is a type $\mathcal{A}$ space.)
\end{theorem}

\proof 
This follows from Proposition~\ref{cl-classes}.10, since by Proposition~\ref{anal-subsets}.2 $\mathcal{A}$ is closed under
finite products and continuous images and by Proposition~\ref{anal-subsets}.5 $\mathcal{A} \subset \mathcal{R}_\specM$.
\eop

\begin{theorem}
Let $(X,\mathcal{E})$ be a separable type $\mathcal{B}$ and $(Y,\mathcal{F})$ a separable measurable space which is countably
separated. If there exists a
bijective measurable mapping $f : (X,\mathcal{E}) \to (Y,\mathcal{F})$ then $\mathcal{F}$ is countably generated,
$(Y,\mathcal{F})$ is a type $\mathcal{B}$ space and $f$ is an isomorphism.
\end{theorem}

\proof 
This is just a special case of Theorem~\ref{type-ba}.6 (together with Theorem~\ref{type-ba}.4). \eop

\textit{Proof of Theorem~\ref{type-ba}.1}\enskip
(2)\enskip
Let $(X,\mathcal{E})$ be a separable type $\mathcal{B}$-space; there thus exists an exactly measurable mapping
$f : (X,\mathcal{E}) \to (\specM,\mathcal{B})$ with $B = f(X) \in \mathcal{B}$. Moreover, since $(X,\mathcal{E})$ is separable,
Proposition~\ref{e-meas-maps}.5 implies $f$ is injective and hence by Proposition~\ref{e-meas-maps}.3  the mapping
$f : (X,\mathcal{E}) \to (B,\mathcal{B}_{|B})$ is an isomorphism.

Suppose first that $B$ is uncountable; then by Proposition~\ref{anal-subsets}.8 there exists an isomorphism
$h : (B,\mathcal{B}_{|B}) \to (\specM,\mathcal{B})$ and therefore $h \circ f : (X,\mathcal{E}) \to (\specM,\mathcal{B})$ is an isomorphism,
i.e., $(X,\mathcal{E})$ is isomorphic to the standard Borel space $(\specM,\mathcal{B})$. But clearly any space isomorphic to a standard Borel
space is itself standard Borel, and thus $(X,\mathcal{E})$ is a standard Borel space. 

Suppose now that $B$ is countable; then $\mathcal{B}_{|B} = \mathcal{P}(B)$ and it follows that $\mathcal{B}_{|B}$ is the Borel $\sigma$-algebra
of $B$ considered with the discrete topology. But $B$ with this topology is a Polish space. Again
$(X,\mathcal{E})$ is isomorphic to a standard Borel space which implies that it is itself standard Borel.

(1)\enskip
This will be proved as follows:
Let $I$ be the closed interval $[0,1]$ and $\mathcal{B}_I$ be the $\sigma$-algebra of Borel subsets of $I$. 
We first show that $(I,\mathcal{B}_I)$ is a type $\mathcal{B}$ space.
It then follows from Proposition~\ref{type-ba}.2~(3) that
$(I^\Nat,\mathcal{B}_I^\Nat)$ is a type $\mathcal{B}$ space.
Thus by Proposition~\ref{type-ba}.2~(2) it is enough to construct
an exactly measurable mapping $h : (X,\mathcal{B}_X) \to (I^\Nat,\mathcal{B}_I^\Nat)$ with $h(X) \in \mathcal{B}_I^\Nat$.

Here are the details: Let $b : \specM \to I$ be the mapping with
\[b(\{z_n\}_{n\ge 0}) = \sum_{n\ge 0} 2^{-n-1} z_n\;;\] 
then $b$ is continuous and hence  $b^{-1}(\mathcal{B}_I) \subset \mathcal{B}$.
Now for each $C \in \mathcal{C}_{\specM}$ there is a dyadic interval $J$ such that $b^{-1}(J) = C$ and hence
$b^{-1}(\mathcal{B}_I) \supset \mathcal{C}_{\specM}$, which implies that $b^{-1}(\mathcal{B}_I) \supset \sigma(\mathcal{C}_{\specM}) = \mathcal{B}$,
i.e., $b^{-1}(\mathcal{B}_I) = \mathcal{B}$. Put
\[ N = \bigl\{ \{z_n\}_{n\ge 0} \in \specM : z_0 = 0\ \mbox{and}\ z_n = 1\ \mbox{for all}\ n \ge m\ \mbox{for some}\ m \ge 1 \bigr\}\;; \]
then $N$ is countable and $b$ maps $\specM_0 = \specM \setminus N$ bijectively onto $I$. Let $v : I \to \specM$
be the unique mapping with $v(b(z)) = z$ for all $z \in \specM_0$; then $v(I) = \specM_0$ and so in particular
$v(I) \in \mathcal{B}$. Let $B \in \mathcal{B}$; then $B \cap \specM_0 \in \mathcal{B}$ and thus there exists $A \in \mathcal{B}_I$
with $b^{-1}(A) = B \cap \specM_0$, which implies  $v^{-1}(B) = v^{-1}(B \cap M_0) = A \in \mathcal{B}_I$. This
shows $v^{-1}(\mathcal{B}) \subset \mathcal{B}_I$. But if $A \in \mathcal{B}_I$ then $b^{-1}(A) \in \mathcal{B}$ and 
$v^{-1}(b^{-1}(A)) = A$, and hence $v^{-1}(\mathcal{B}) = \mathcal{B}_I$. 
Therefore $v : (I,\mathcal{B}_I) \to (\specM,\mathcal{B})$ is an exactly measurable mapping 
with $v(I) \in \mathcal{B}$, and hence $(I,\mathcal{B}_I)$ is a type $\mathcal{B}$ space.
As mentioned above, it now follows from Proposition~\ref{type-ba}.2~(3) that $(I^\Nat,\mathcal{B}_I^\Nat)$ is a type $\mathcal{B}$ space.
Note that the product $\sigma$-algebra $\mathcal{B}_I^\Nat$
is the Borel $\sigma$-algebra of $I^\Nat$ with the product topology (since $\Nat$ is countable and $I$ has a countable base
for its topology).

Finally, let $(X,d)$ be a complete separable metric space. Then there is a standard construction (given below) producing a continuous 
injective mapping $h : X \to I^\Nat$ such that $h$ is a homeomorphism from $X$ to $h(X)$ (with the relative topology) and such that
$h(X)$ is the intersection of a sequence of open subsets of $I^\Nat$. In particular
with $h(X) \in \mathcal{B}_I^\Nat$, and $h^{-1}(\mathcal{B}_I^\Nat) \subset \mathcal{B}_X$, since $h$ is continuous.
On the other hand, let $U \subset X$ be open; since $h : X \to h(X)$ is a homeomorphism there exists an open subset $V$ of
$I^\Nat$ with $h^{-1}(h(X) \cap V) = U$. But then $h^{-1}(V) = U$, and this shows that
$\mathcal{O}_X \subset h^{-1}(\mathcal{B}_I^\Nat)$, where $\mathcal{O}_X$ is the set of open subsets of $X$. Hence
$\mathcal{B}_X = \sigma(\mathcal{O}_X) 
\subset \sigma(h^{-1}(\mathcal{B}_I^\Nat)) = h^{-1}(\sigma(\mathcal{B}_I^\Nat)) = h^{-1}(\mathcal{B}_I^\Nat)$ and thus
$h : (X,\mathcal{B}_X) \to (I^\Nat,\mathcal{B}_I^\Nat)$ is exactly measurable and $h(X) \in \mathcal{B}_I^\Nat$.
Therefore by Proposition~\ref{type-ba}.2~(2) $(X,\mathcal{B}_X)$ is also a type $\mathcal{B}$ space.

Here is how the mapping $h$ can be constructed:
Choose a dense sequence of elements $\{x_n\}_{n\ge 0}$ from $X$
and for each $n \ge 0$ let $h_n : X \to I$ be the continuous mapping given by $h_n(x) = \min \{ d(x,x_n),1\}$
for each $x \in X$. If $x,\,y \in X$ with $x \ne y$ then there exists $n \ge 0$ such that $h_n(x) \ne h_n(y)$ 
(since if $n \ge 0$ is such that $d(x,x_n) < \varepsilon$, where 
$\varepsilon = \half \min \{ d(x,y),1\}$, then $h_n(x) < \varepsilon <  h_n(y)$).
Define $h : X \to I^\Nat$ by letting
$h(x) = \{h_n(x)\}_{n\ge 0}$ for each $x \in X$.
Then $h$ is continuous (since $\proj_n \circ h = h_n$ is continuous for each $n \ge 0$, with 
$\proj_n : I^\Nat \to I$ the projection onto the $n$ th component) and injective.
Moreover, for each $x \in X$ and each $0 < \varepsilon < 1$ there exists $n \ge 0$ so that
$|h_n(y) - h_n(x)| > \varepsilon/2$ for all $y \in X$ with $d(y,x) > \varepsilon$. (Just take $n \ge 0$ so that
$d(x,x_n) < \varepsilon/4$.) This implies that the bijective mapping $h : X \to h(X)$ is a homeomorphism
from $X$ to $h(X)$ with the relative topology.
(Note that the completeness of $X$ was not needed here.)

The topological space $I^\Nat$ is metrisable (since $\Nat$ is countable and $I$ is metrisable).
Let $\delta$ be any metric generating the topology on $I^\Nat$, and for each $n \ge 0$ let 
\[U_n = \{ y \in I^\Nat : \delta(y,y') < 2^{-n}\ \mbox{for some}\ y' \in h(X) \}\;.\]
Then $\{U_n\}_{n\ge 0}$ is a decreasing sequence of open subsets of $I^\Nat$ with $h(X) \subset U_n$ for each
$n \ge 0$. In fact $h(X) = \bigcap_{n\ge 0} U_n$:
Let $y \in \bigcap_{n\ge 0} U_n$; then for each $n \ge 0$ there exists $y_n \in h(X)$ with $\delta(y,y_n) < 2^{-n}$
and so $\{y_n\}_{n\ge 0}$ is a Cauchy sequence in $h(X)$. Thus $\{x_n\}_{n\ge 0}$ is a Cauchy sequence in $X$,
where $x_n$ is the unique element with $h(x_n) = y_n$. Since $X$ is complete the sequence
$\{x_n\}_{n\ge 0}$ has a limit $x \in X$ and then $h(x) = y$, i.e., $y \in h(X)$.
This shows that $h(X)$ is the intersection of a sequence of open subsets of $I^\Nat$. \eop


\startsection{Universal measurability}
\label{uni-meas}

In Section~\ref{anal-subsets}\ we saw (in Proposition~\ref{anal-subsets}.1) that each Borel subset of $\specM$ is analytic, i.e., 
$\mathcal{B} \subset \mathcal{A}$. We also noted that there exist analytic sets which are not Borel.
However, when dealing with finite measures this fact usually doesn't cause a problem: In this section we will see that if
$\mu$ is a finite measure on $\mathcal{B}$ and $A \in \mathcal{A}$ then there exist $B^-,\,B^+ \in \mathcal{B}$ 
with $B^- \subset A \subset B^+$ such that $\mu(B^+ \setminus B^-) = 0$, and so the difference between $A$ and a Borel set
is `negligible'.
As is usual we treat this topic using the notion of what is called universal measurability.

A \definition{measure space} is a triple $(X,\mathcal{E},\mu)$, where $(X,\mathcal{E})$ is a measurable space and $\mu$ is 
a measure on $\mathcal{E}$. If the measure is finite then we say that the measure space is \definition{finite}.
The measure space is said to be \definition{complete} if $N \in \mathcal{E}$ whenever
$N \subset N'$ for some $N' \in \mathcal{E}$ with $\mu(N') = 0$.

Let $(X,\mathcal{E},\mu)$ be a measure space and let $\mathcal{E}_\mu$ consist of  those subsets $E$ of $X$ for which there
exist $E^-,\,E^+ \in \mathcal{E}$ with $E^- \subset E \subset E^+$ and $\mu(E^+ \setminus E^-) = 0$ (and note that then
$\mu(E^-) = \mu(E^+)$). If $E_1^-,\,E_1^+,\,E_2^-,\,E_2^+ \in \mathcal{E}$ with $E_1^- \subset E \subset E_1^+$, $E_2^- \subset E \subset E_2^+$ and 
$\mu(E_1^+ \setminus E_1^-) = \mu(E_2^+ \setminus E_2^-) = 0$ then $\mu(E_1^-) = \mu(E_2^-)$ and hence there is a unique mapping
$\bar{\mu} : \mathcal{E}_\mu \to \Realpi$ such that 
$\bar{\mu}(E) = \mu(E^-)$ whenever $E^- \subset E \subset E^+$ with $\mu(E^+ \setminus E^-) = 0$.
Then $\mathcal{E}_\mu$ is a $\sigma$-algebra which contains $\mathcal{E}$, $\bar{\mu}$ is  the unique measure on $\mathcal{E}_\mu$ which
extends $\mu$ and the measure space $(X,\mathcal{E}_\mu,\bar{\mu})$ is complete. 
This measure space is called the \definition{completion} of $(X,\mathcal{E},\mu)$.

Now for each measurable space $(X,\mathcal{E})$ we denote by
$\mathcal{E}_*$  the intersection of all the $\sigma$-algebras $\mathcal{E}_\mu$, the intersection
being taken over all finite measures $\mu$ (or, equivalently, over all probability measures) defined on $\mathcal{E}$.
Then $\mathcal{E}_*$ is a $\sigma$-algebra with $\mathcal{E} \subset \mathcal{E}_*$ 
which is called the
\definition{$\sigma$-algebra of universally measurable sets} (with respect to $\mathcal{E}$).
Each finite measure $\mu$ on $\mathcal{E}$ has a extension to a measure on 
$\mathcal{E}_*$, namely the restriction of $\bar{\mu}$ to $\mathcal{E}_*$, and this measure will also be denoted
by $\bar{\mu}$. Again, $\bar{\mu}$ is the unique extension of $\mu$ to a measure on $\mathcal{E}_*$.

Here is the main main result in this section:

\begin{theorem}
$\mathcal{A} \subset \mathcal{B}_*$.
\end{theorem}

\proof
This will follow immediately from Proposition~\ref{uni-meas}.1 below, which states that $\mathcal{A} \subset \mathcal{B}_\mu$
for each finite measure $\mu$ on $\mathcal{B}$. \eop

The proof of Proposition~\ref{uni-meas}.1 is taken from Section~8.6 of Cohn \cite{cohn} and involves outer and inner measures. 
In what follows let $(X,\mathcal{E},\mu)$ be a finite measure space.
The outer measure $\mu^* : \mathcal{P}(X) \to \Realpos$ associated with $\mu$ is defined by
\[ \mu^*(B) = \inf \{\, \mu(E) : \mbox{$E \in \mathcal{E}$ with $B \subset E$} \,\}\]
and the inner measure $\mu_* : \mathcal{P}(X) \to \Realpos$ by
\[ \mu_*(B) = \sup \{\, \mu(E) : \mbox{$E \in \mathcal{E}$ with $E \subset B$} \,\} \]
for all $F \subset X$. 
Thus $\mu^*$ and $\mu_*$ are both increasing (meaning  $\mu^*(B) \le \mu^*(B')$ and $\mu_*(B) \le \mu_*(B')$ hold whenever
$B \subset B'$); moreover,
$\mu_*(B) \le \mu^*(B)$ holds for all $B \subset X$ and $\mu_*(E) = \mu(E) = \mu^*(E)$ for all $E \in \mathcal{E}$.

\begin{lemma}
$\mathcal{E}_\mu = \{ B \in \mathcal{P}(X) : \mu_*(B) = \mu^*(B) \}$ and
$\mu_*(E) = \bar{\mu}(E) = \mu^*(E)$ for all $E \in \mathcal{E}_\mu$.
\end{lemma}

\proof 
If $E \in  \mathcal{E}_\mu$ then there exist $E^-,\,E^+ \in \mathcal{E}$ with $E^- \subset E \subset E^+$ such that
$\mu(E^+ \setminus E^-) = 0$. Therefore
$\mu_*(E) \ge \mu(E^-) = \mu(E^+) \ge \mu^*(E)$
and it follows that $\mu_*(E) = \bar{\mu}(E) = \mu^*(E)$.
Suppose conversely that $B \subset X$ with $\mu_*(B) = \mu^*(B)$.
For each $n \ge 1$ there then exist $E^-_n,\,E^+_n \in \mathcal{B}$ with $E^-_n \subset B \subset E^+_n$ such that
$\mu(E^-_n) > \mu_*(B) - 1/n$ and $\mu(E^+_n) < \mu_*(B) + 1/n$. Put $E^- = \bigcup_{n\ge 1} E^-_n$ and
$E^+ = \bigcap_{n\ge 1} E^+_n$; then $E^-,\,E^+ \in \mathcal{E}$, $E^- \subset B \subset E^+$ and
\[\mu(E^+ \setminus E^-) \le \mu(E^+_n \setminus E^-_n) = \mu(E^+_n) - \mu(E^-_n) < 2/n\] 
for all $n \ge 1$, i.e., $\mu(E^+ \setminus E^-) = 0$. Therefore $B \in \mathcal{E}_\mu$. \eop

\begin{lemma}
Let $\{B_n\}_{n\ge 0}$ be any increasing sequence of subsets of $X$. Then 
$\mu^*\bigl(\bigcup_{n \ge 0} B_n\bigr) = \lim_{n} \mu^*(B_n)$.
\end{lemma}

\proof 
Put $B = \bigcup_{n \ge 0} B_n$.
The sequence $\{\mu^*(B_n)\}_{n \ge 0}$ is increasing and bounded above by $\mu^*(B)$; thus
it converges  and $\lim_{n} \mu^*(B_n) \le \mu^*(B)$.
Let $\varepsilon > 0$ and for each $n \ge 0$ let $E'_n \in \mathcal{E}$ be such that 
$B_n \subset E'_n$ and $\mu(E'_n) \le \mu^*(B_n) + \varepsilon$.
Put $E_n = \bigcup_{k \ge n} E'_k$ for each $n \ge 0$; then $\{E_n\}_{n \ge 0}$ is an increasing sequence
from $\mathcal{E}$ with $B_n \subset E_n$ and $\mu(E_n) < \mu^*(B_n) + \varepsilon$ for all $n \ge 0$.
Let $E = \bigcup_{n\ge 0} E_n$; then $E \in \mathcal{E}$, $B \subset E$ and
$\lim_{n} \mu(E_n) = \mu(E)$. Therefore
\[ \mu^*(B) \le \mu(E) = \lim_{n \to \infty} \mu(E_n) \le  \lim_{n \to \infty} \mu^*(B_n) + \varepsilon \] 
and, since $\varepsilon > 0$ is arbitrary, it follows that $\mu^*(B) \le \lim_{n} \mu^*(B_n) + \varepsilon$. This shows that
$\mu^*(B) = \lim_{n} \mu^*(B_n)$.
\eop

\begin{proposition}
If $\mu$ is a finite measure on $\mathcal{B}$ then $\mathcal{A} \subset \mathcal{B}_\mu$.
\end{proposition}

\proof 
For $n,\, m_0,\,\ldots,\,m_n \in \Nat$ let 
\[  \specN^*(m_0,\ldots,m_n) = \{ \{m'_p\}_{p \ge 0} \in \specN :  m'_j \le m_j\ \mbox{for}\ j = 0,\,\ldots,\,n \}\;. \]
Note that $\specN^*(m_0,\ldots,m_n,p) \subset \specN^*(m_0,\ldots,m_n,p+1)$ for all $p \in \Nat$
and that $\specN^*(m_0,\ldots,m_n) = \bigcup_{p \ge 0} \specN^*(m_0,\ldots,m_n,p)$ for all
$n,\, m_0,\,\ldots,\,m_n \in \Nat$.

Let $A \in \mathcal{A}$ and since $\varnothing \in \mathcal{B}_\mu$ we can assume that $A \ne \varnothing$. Thus there exists
a continuous mapping $f : \specN \to \specM$ with $f(\specN) = A$. Let $\varepsilon > 0$; then 
$\{f(\specN^*(p)\}_{p \ge 0}$ is an increasing sequence of subsets of $\specM$ with
$\bigcup_{p\ge 0} f(\specN^*(p)) = f(\specN) = A$, and so by Lemma~\ref{uni-meas}.2 there exists $m_0 \in \Nat$, so that
$\mu^*(f(\specN^*(m_0))) > \mu^*(A) - \varepsilon$. In the same way
$\{f(\specN^*(m_0,p)\}_{p \ge 0}$ is an increasing sequence of subsets of $\specM$ with
$\bigcup_{p\ge 0} f(\specN^*(m_0,p)) = f(\specN^*(m_0))$, and so by Lemma~\ref{uni-meas}.2 there exists $m_1 \in \Nat$, so that
$\mu^*(f(\specN^*(m_0,m_1))) > \mu^*(A) - \varepsilon$. Iterating the process results in an element
$\mathsf{m} = \{m_n\}_{n \ge 0}$ of $\specN$ with
$\mu^*(f(\specN^*(m_0,\ldots,m_n))) > \mu^*(A) - \varepsilon$ for each $n \in \Nat$.
Put 
\[\specN^*(\mathsf{m}) = \bigcap_{n \ge 0} \specN^*(m_0,\ldots,m_n)
= \{ \{m'_n\}_{n \ge 0} \in \specN :  m'_n \le m_n\ \mbox{for all}\ n \in \Nat \}\;;\]
then $\specN^*(\mathsf{m})$ is a compact subset of $\specN$ (since it is the product of finite, and hence compact, subsets of $\Nat$), and therefore
$K = f(\specN^*(\mathsf{m}))$ is a compact subset of $\specM$ with, of course, $K \subset A$.
Put $K_n = \overline{f(\specN^*(m_0,\ldots,m_n))}$ for each $n \in \Nat$
(with $\overline{B}$ denoting the closure of the set $B$).

\begin{lemma}
$K = \bigcap_{n \ge 0} K_n$.
\end{lemma}

\proof 
For each $n \in \Nat$ let $\mathsf{m}_n \in \specN^*(m_0,\ldots,m_n)$; then for each $p \in \Nat$ the $p$ th components of the
elements in the sequence $\{\mathsf{m}_n\}_{n \ge 0}$ are bounded. Hence by the usual diagonal argument there exists a 
subsequence $\{n_k\}_{k \ge 0}$ so that the sequence
$\{\mathsf{m}_{n_k}\}_{k \ge 0}$ converges. Moreover, the limit lies in
$\specN^*(m_0,\ldots,m_n)$ for each $n \in \Nat$ and thus in $\specN^*(\mathsf{m})$. 

Let $z \in \bigcap_{n\ge 0} K_n$.
Then for each $n \in \Nat$ there exists $\mathsf{m}_n \in \specN^*(m_0,\ldots,m_n)$
such that $\{f(\mathsf{m}_n)\}_{n \ge 0}$ converges to $z$ in $\specM$.
Let $\{n_k\}_{k \ge 0}$ be a subsequence such that the sequence
$\{\mathsf{m}_{n_k}\}_{k \ge 0}$ converges and let $\mathsf{p} \in \specN^*(\mathsf{m})$ be the limit.
Then $\{f(\mathsf{m}_{n_k})\}_{k \ge 0}$ still converges to $z$ and hence
$f(\mathsf{p}) = z$, since $f$ is continuous. This shows that $\bigcap_{n\ge 0} K_n \subset K$. But $K \subset \bigcap_{n\ge 0} K_n$ holds 
trivially, since $K = \overline{K} \subset K_n$ for each $n \in \Nat$. Therefore $K = \bigcap_{n\ge 0} K_n$. \eop

Now $\{K_n\}_{n\ge 0}$ is a decreasing sequence of closed subsets of $\specM$ (and in particular
of elements from $\mathcal{B}$) and hence by Lemma~\ref{uni-meas}.3
$\mu(K) = \lim_{n} \mu(K_n)$.
since $\mu$ is a finite measure. But for each $n \ge 0$
\[
 \mu^*(A) - \varepsilon < \mu^*(f(\specN^*(m_0,\ldots,m_n))) \le \mu^*(K_n) = \mu(K_n) \]
and therefore $\mu(K) \ge \mu^*(A) - \varepsilon$. Thus
$\mu_*(A) \ge \mu^*(A) - \varepsilon$, since $\mu(K) \le \mu_*(A)$, which implies that $\mu_*(A) = \mu^*(A)$.
It follows from Lemma~\ref{uni-meas}.1
that $A \in \mathcal{B}_\mu$. This shows that $\mathcal{A} \subset \mathcal{B}_\mu$ and completes the proof of
Proposition~\ref{uni-meas}.1.
\eop

Now let $(X,\mathcal{E})$ and $(Y,\mathcal{F})$ be measurable spaces and 
$f : (X,\mathcal{E}) \to (Y,\mathcal{F})$ be a measurable mapping. Then
for each measure $\mu$ on $\mathcal{E}$ there is a measure $\mu f^{-1}$
on $\mathcal{F}$ called the \definition{image of $\mu$ under $f$}
and defined by $(\mu f^{-1})(F) = \mu(f^{-1}(F))$ for each $F \in \mathcal{F}$.

\begin{lemma}
Let $\mu$ be a finite measure on $\mathcal{E}$ and 
let $\nu = \mu f^{-1}$ be the image measure.
Then $f^{-1}(\mathcal{F}_\nu) \subset \mathcal{E}_\mu$ and $\bar{\nu} = \bar{\mu} f^{-1}$.
\end{lemma}

\proof 
Let $F \in \mathcal{F}_\nu$; there thus exist $F^-,\,F^+ \in \mathcal{F}$ with $F^- \subset F \subset F^+$ 
such that $\nu(F^+ \setminus F^-) = 0$. Put $E^- = f^{-1}(F^-)$ and $E^+ = f^{-1}(F^+)$; then 
$E^-,\,E^+ \in \mathcal{E}$ with $E^- \subset f^{-1}(E) \subset E^+$ and
$\mu(E^+ \setminus E^-) = \nu(f^{-1}(E^+ \setminus E^-)) = \nu(F^+ \setminus F^-) = 0$.
Hence $f^{-1}(F) \in \mathcal{E}_\mu$, which shows that $f^{-1}(\mathcal{F}_\nu) \subset \mathcal{E}_\mu$.
Finally, $\bar{\mu} f^{-1}$ is a measure on $\mathcal{F}_\nu$ which is an extension of the measure
$\nu = \mu  f^{-1}$ on $\mathcal{F}$, and so $\bar{\nu} = \bar{\mu} \circ f^{-1}$ by the uniqueness of 
$\bar{\nu}$.  \eop

\begin{lemma}
$f^{-1}(\mathcal{E}_*) \subset \mathcal{F}_*$.
Moreover, if $\mu$ is a finite measure on $\mathcal{E}$ and $\nu = \mu f^{-1}$ then
$\bar{\nu} = \bar{\mu} f^{-1}$, here with $\bar{\mu}$ and $\bar{\nu}$ the measures on 
$\mathcal{E}_*$ and $\mathcal{F}_*$ respectively.
\end{lemma}

\proof 
This follows from Lemma~\ref{uni-meas}.4. \eop


\startsection{Measurable selectors}
\label{selectors}

Let $(X,\mathcal{E})$ and $(Y,\mathcal{F})$ be measurable spaces 
and let $f : (X,\mathcal{E}) \to (Y,\mathcal{F})$ be a surjective measurable mapping.
By the axiom of choice there then exists a \definition{selector} for $f$, i.e., a mapping $g : Y \to X$ such that 
$f \circ g = \id_Y$. Unfortunately, it is not always possible to choose $g$ to be a measurable mapping from
$(X,\mathcal{E})$ to $(Y,\mathcal{F})$, even when 
$(X,\mathcal{E})$ and $(Y,\mathcal{F})$ are separable type $\mathcal{B}$ spaces: Let $I = [0,1]$ and let
$\proj_1 : I \times I \to I$ be the projection onto the first component. Then there exists a Borel subset $A$ of
$I \times I$ with $\proj_1(A) = I$ for which there does not exist a Borel measurable mapping $g : I \to I \times I$ with
$g(I) \subset A$ such that $\proj_1(g(x)) = x$ for all $x \in I$. (See, for example, Blackwell \cite{blackwell}.) 

However, the following result due to Yankov \cite{yankov} and von Neumann \cite{vonneumann} 
shows that universally measurable selectors exist.
Recall that if $(Y,\mathcal{F})$ is a measurable space then, as introduced in the previous section, $\mathcal{F}_*$ denotes the corresponding 
$\sigma$-algebra of universally measurable sets.

\begin{theorem}
Let $f : (X,\mathcal{E}) \to (Y,\mathcal{F})$ be a surjective measurable mapping
with $(X,\mathcal{E})$ and $(Y,\mathcal{F})$ separable type $\mathcal{A}$ spaces.
Then there is a measurable mapping
$g : (Y,\mathcal{F}_*) \to (X,\mathcal{E})$  such that $f(g(y)) = y$ for all $y \in Y$.
\end{theorem}

\proof 
This is given below. The main step in the proof is to establish an analogous result 
(Proposition~\ref{selectors}.1) for a continuous mapping $f : \specM \to \specM$.
The proof of Proposition~\ref{selectors}.1 is taken from Section~8.5 of Cohn \cite{cohn}.
\eop

Denote the lexicographical order on $\specN$ by $\preceq$; thus $\mathsf{m} \preceq \mathsf{m}$ for all $\mathsf{m} \in \specN$ 
and if $\mathsf{m} = \{m_n\}_{n \ge 0},\, \mathsf{m}' = \{m'_n\}_{n\ge 0}$  with
$\mathsf{m} \ne \mathsf{m}'$ then $\mathsf{m} \preceq \mathsf{m}'$ if and only if
$m_p < m'_p$, where $p = \min \{ n \ge 0 : m_n \ne m'_n \}$. Clearly $\preceq$ defines a total order on $\specN$.

\begin{lemma}
Each non-empty closed subset of $\specN$ possesses a least element (with respect to $\preceq$).
\end{lemma}

\proof 
Let $D$ be a non-empty closed subset of $\specN$. As in the proof of Lemma~\ref{anal-subsets}.2 let
$D(m_0,\ldots,m_n) = D \cap \specN(m_0,\ldots,m_n)$ for all $n,\,m_0,\,\ldots,\,m_n \in \Nat$. Note that
$D(p) \ne \varnothing$ for some $p \in \Nat$ and if
$D(m_0,\ldots,m_n) \ne \varnothing$ then $D(m_0,\ldots,m_n,p) \ne \varnothing$ for at least one $p \in \Nat$.
Define $\mathsf{m} = \{m_n\}_{n\ge 0} \in \specN$ with $D(m_0,\ldots,m_n) \ne \varnothing$ for all $n \in \Nat$
inductively as follows: Put $m_0 = \min \{ k \in \Nat : D(k) \ne \varnothing \}$ and for $n \ge 0$ let
$m_{n+1} = \min \{ k \in \Nat : D(m_0,\ldots,m_n,k) \ne \varnothing \}$.
Then $\mathsf{m} \in D$, since $D$ is closed. Moreover, if $\mathsf{m}' = \{m'_n\}_{n \ge 0}$ with
$\mathsf{m}' \ne \mathsf{m}$ and $p = \min \{ n \ge 0 : m'_n \ne m_n \}$ then $m'_p < m_p$ (by the definition of
$m_p$) and hence $\mathsf{m}' \prec \mathsf{m}$. This implies that $\mathsf{m}$ is the least element of $D$ 
(with respect to $\preceq$). \eop

For each $\mathsf{m} \in \specN$ let $U_{\mathsf{m}} = \{ \mathsf{n} \in \specN : \mathsf{n} \prec \mathsf{m} \}$.
Denote the $\sigma$-algebra of Borel subsets of $\specN$ by $\mathcal{B}_{\specN}$.

\begin{lemma}
The set $U_{\mathsf{m}}$ is open for each $\mathsf{m} \in \specN$. Moreover,
$\mathcal{B}_{\specN}$ is the smallest $\sigma$-algebra containing the sets $U_{\mathsf{m}}$, $\mathsf{m} \in \specN$.
\end{lemma}

\proof 
Let $\mathsf{m} = \{m_n\}_{n\ge 0} \in \specN$; then
\[   U_\mathsf{m} = \bigcup_{k < m_0} \specN(k) \cup 
\bigcup_{p \ge 1} \bigcup_{k < m_p}  \specN(m_0,\ldots,m_{p-1},k)\;, \]
which is a union of open sets and thus $U_\mathsf{m}$ is open.
Now let $\mathcal{B}'_{\specN}$ be the smallest $\sigma$-algebra containing the sets $U_{\mathsf{m}}$, $\mathsf{m} \in \specN$; 
thus $\mathcal{B}'_{\specN} \subset \mathcal{B}_{\specN}$, since 
$U_{\mathsf{m}}$ is open for each $\mathsf{m} \in \specN$. 
Let $p,\,m_0,\,\ldots,\,m_p \in \Nat$, and let
$\mathsf{k} = \{k_n\}_{n\ge 0}$ and $\mathsf{k}' = \{k'_n\}_{n\ge 0}$ be the elements of $\specN$ with 
$k_j = k'_j = m_j$ for $j = 0,\,\ldots,\,p-1$, $k_p = m_p + 1$, $k'_p = m_p$ and $k_j = k'_j = 0$ for all $j > p$. 
Then $\specN(m_0,\ldots,m_p) = U_\mathsf{k} \setminus U_{\mathsf{k}'}$, and this implies that
$\mathcal{C}^o_\specN \subset \mathcal{B}'_{\specN}$.  Therefore
$\mathcal{B}_\specN = \sigma(\mathcal{C}^o_\specN) \subset \mathcal{B}'_{\specN}$,
i.e., $\mathcal{B}_\specN = \mathcal{B}'_{\specN}$.  \eop

\begin{lemma}
Let $f : \specN \to \specM$ be continuous; then $f(U) \in \mathcal{A}$ for each open subset $U$ of $\specN$.
\end{lemma}

\proof 
We can assume that $U \ne \varnothing$; thus, since $\mathcal{C}^o_\specN$ is a base for the topology on $\specN$ there exists a
sequence $\{C_n\}_{n \ge 0}$ from $\mathcal{C}^o_\specN$ with $U = \bigcup_{n \ge 0} C_n$ and then
$f(U) = \bigcup_{n \ge 0} f(C_n)$. But the elements of $\mathcal{C}^o_\specN$ are also closed and so by Lemma~\ref{anal-subsets}.3
$f(C_n) \in \mathcal{A}$ for each $n \ge 0$. Therefore by Lemma~\ref{anal-subsets}.4~(1) $f(U) \in \mathcal{A}$. 
\eop

Put $\mathcal{A}^\sigma = \sigma(\mathcal{A})$, i.e., $\mathcal{A}^\sigma$ is the smallest $\sigma$-algebra containing the
analytic subsets of $\specM$.

\begin{proposition}
Let $f : \specM \to \specM$ be a continuous mapping and let $A$ be a non-empty analytic subset of $\specM$; put $C = f(A)$.
Then there exists a measurable mapping $g : (C,\mathcal{A}^{\sigma}_{|C}) \to (\specM,\mathcal{B})$ with $g(C) \subset A$ and such that
$f(g(z)) = z$ for all $z \in C$.
\end{proposition}

\proof
There exists a continuous mapping $\tau : \specN \to \specM$ with $\tau(\specN) = A$ and then the mapping
$p = f \circ \tau : \specN \to \specM$ is continuous with $p(\specN) = C$.
Now $p^{-1}(\{z\})$ is a non-empty closed subset of $\specN$ for each $z \in C$ and so by Lemma~\ref{selectors}.1 we can define a mapping
$q : C \to \specN$ by letting $q(z)$ be the least element in $p^{-1}(\{z\})$ for each $z \in C$.
Put $g = \tau \circ q$; then $g : C \to \specM$ with $g(C) \subset \tau(\specN) = A$ and
$f(g(z)) = f(\tau(q(z))) = p(q(z)) = z$ for all $z \in C$.
Now
$\tau^{-1}(\mathcal{B}) \subset \mathcal{B}_{\specN}$, since $\tau$ is continuous, and $g^{-1}(B) = q^{-1}(\tau^{-1}(B))$ for
each $B \in \mathcal{B}$ and so it is now enough to show that
$q^{-1}(\mathcal{B}_{\specN}) \subset \mathcal{A}^{\sigma}_{|C}$.
But $q^{-1}(U_{\mathsf{m}}) = p(U_{\mathsf{m}})$
for all $\mathsf{m} \in \specN$. (If $z \in q^{-1}(U_{\mathsf{m}})$ then $q(z) \in U_{\mathsf{m}}$ and so
$z = p(q(z)) \in p(U_{\mathsf{m}})$, i.e.,  $q^{-1}(U_{\mathsf{m}}) \subset p(U_{\mathsf{m}})$.
On the other hand, if $z \in p(U_{\mathsf{m}})$ with $z = p(\mathsf{n})$ and $\mathsf{n} \in U_{\mathsf{m}}$ then
$q(z) \preceq \mathsf{n}$, since
$\mathsf{n} \in p^{-1}(\{z\})$ and $q(z)$ is the smallest element in this set.
Thus $q(z) \preceq \mathsf{n} \prec \mathsf{m}$, which means that $q(z) \in U_{\mathsf{m}}$, i.e., 
$p(U_{\mathsf{m}}) \subset q^{-1}(U_{\mathsf{m}})$. Hence $q^{-1}(U_{\mathsf{m}}) = p(U_{\mathsf{m}})$.)
Moreover, by Lemma~\ref{selectors}.3 the set $p(U_{\mathsf{m}})$ is analytic, and  $p(U_{\mathsf{m}}) = p(U_{\mathsf{m}}) \cap C$ and  
so $q^{-1}(U_{\mathsf{m}}) \in \mathcal{A}^{\sigma}_{|C}$ for each  $\mathsf{m} \in \specN$. 
Therefore by Lemma~\ref{selectors}.2 $q^{-1}(\mathcal{B}_{\specN}) \subset \mathcal{A}^{\sigma})$.
\eop

\textit{Proof of Theorem~\ref{selectors}.1:}
As usual let $\Delta : \specM \to \specM$ be the continuous surjective mapping given by
$\Delta(\{z_n\}_{n\ge 0}) = \{z'_n\}_{n\ge 0}$,
where $z'_n = z_{2n+1}$ for each $n \ge 0$.
By Proposition~\ref{cl-classes}.7 there exist exactly measurable mappings $p : (X,\mathcal{E}) \to (\specM,\mathcal{B})$
and $q : (Y,\mathcal{F}) \to (\specM,\mathcal{B})$ with $p(X),\,q(Y) \in \mathcal{A}$ such that $q \circ f = \Delta \circ p$.
Put $A = p(X)$ and $C = q(Y)$, and hence $\Delta(A) = \Delta(p(X)) = q(f(X)) = q(Y) = C$.
Then by Proposition~\ref{selectors}.1 there
exists a measurable mapping $h : (C,\mathcal{A}^{\sigma}_{|C}) \to (\specM,\mathcal{B})$ with $h(C) \subset A$ and such that
$\Delta(h(z)) = z$ for all $z \in C$.
Now by Proposition~\ref{e-meas-maps}.4 the mappings $p$ and $q$ are injective and so let 
$p^{-1} : A \to X$ be the inverse of the bijective mapping $p : X \to A$ and $q^{-1} : C \to Y$ be the inverse of the bijective mapping 
$q : Y \to C$. Define $g : Y \to X$ by $g = p^{-1} \circ h \circ q$; then 
\[ q(f(g(y))) =   \Delta(p(g(y)) = \Delta(h(q(y))) = q(y)\]
and hence $f(g(y)) = q^{-1}(q(f(g(y))) = q^{-1}(q(y)) = y$ for all $y \in Y$.
Let $E \in \mathcal{E}$; there thus exists $B \in \mathcal{B}$ with $E = p^{-1}(B)$ and then
\begin{eqnarray*}
g^{-1}(E) &=& (p^{-1} \circ h \circ q)^{-1}(E)\\ 
&=& q^{-1}(h^{-1}((p^{-1})^{-1}(E))) = q^{-1}(h^{-1}(B \cap A)) = q^{-1}(h^{-1}(B))\;.
\end{eqnarray*}
But $h^{-1}(B) \in \mathcal{A}^\sigma_{|C}$ and so $h^{-1}(B) = D \cap C$ for some $D \in \mathcal{A}^\sigma$, which implies
that $g^{-1}(E) = q^{-1}(h^{-1}(B)) = q^{-1}(D \cap C) = q^{-1}(D)$. Moreover, by Theorem~\ref{uni-meas}.1
$\mathcal{A}^\sigma \subset \mathcal{B}_*$ and by Lemma~\ref{uni-meas}.5 $q^{-1}(\mathcal{B}_*) \subset \mathcal{F}_*$,
and therefore $g^{-1}(E) = q^{-1}(D) \in \mathcal{F}_*$. This shows that $g^{-1}(\mathcal{E}) \subset \mathcal{F}_*$, i.e.,
$g : (Y,\mathcal{F}_*) \to (X,\mathcal{E})$ is measurable. \eop


\startsection{A method for constructing measures}
\label{image-measures}

Recall that
if $(X,\mathcal{E})$ and $(Y,\mathcal{F})$ are measurable spaces and 
$f : (X,\mathcal{E}) \to (Y,\mathcal{F})$ is a measurable mapping then
for each measure $\mu$ on $\mathcal{E}$ there is a measure $\mu f^{-1}$
on $\mathcal{F}$ called the \definition{image of $\mu$ under $f$}
and defined by $(\mu f^{-1})(F) = \mu(f^{-1}(F))$ for each $F \in \mathcal{F}$.

In the following sections one of the basic task involves, in some form or another, constructing an element 
$\mu \in \Prob{X}{\mathcal{E}}$ out of a given sequence $\{\mu_n\}_{n\ge 0}$ of elements of $\Prob{X}{\mathcal{E}}$,
where $\Prob{X}{\mathcal{E}}$ denotes the set of probability measures on $(X,\mathcal{E})$.
If $(X,\mathcal{E})$ is countably generated (and in almost all cases $(X,\mathcal{E})$ will be a type $\mathcal{A}$ or a type $\mathcal{B}$
space) then we can try the following:

\begin{evlist}{25pt}{5pt}
\item[(1)] Choose an exactly measurable mapping $f : (X,\mathcal{E}) \to (\specM,\mathcal{B})$.
\item[(2)] Work with the sequence of images
$\{\mu_nf^{-1}\}_{n\ge 0}$ in $\Prob{\specM}{\mathcal{B}}$ 
and exploit the properties of the space $(\specM,\mathcal{B})$ to produce a measure $\nu \in \Prob{\specM}{\mathcal{B}}$.
\item[(3)] Pull the measure $\nu$ back to an element $\mu \in \Prob{X}{\mathcal{E}}$ with $\nu = \mu f^{-1}$.
\end{evlist}
Step (2) is, of course, somewhat vague, but
step (3) can be made more precise: The problem is whether there exists a measure
$\mu \in \Prob{X}{\mathcal{E}}$ with $\nu = \mu f^{-1}$, and this
means that the measure $\nu \in \Prob{\specM}{\mathcal{B}}$ constructed 
in (2) should be such that Proposition~\ref{image-measures}.1 can be applied. First, however, a definition.
Let $(Y,\mathcal{F})$ be a measurable space and 
let $\nu$ be a measure on $\mathcal{F}$; then a subset $B \subset Y$ is said to be 
\definition{thick} with respect to
$\nu$ if $\nu(F) = 0$ for all $F \in \mathcal{F}$ with $F \cap B = \varnothing$.
Of course, an element $B \in \mathcal{F}$ is thick with respect to $\nu$ if and only if $\nu(Y \setminus B) = 0$.

In Proposition~\ref{image-measures}.1 let
$(X,\mathcal{E})$ and $(Y,\mathcal{F})$ be measurable spaces.

\begin{proposition}
Let $f : (X,\mathcal{E}) \to (Y,\mathcal{F})$ be exactly measurable
and let $\nu$ be a measure on $\mathcal{F}$.
Then there exists a measure $\mu$ on $\mathcal{F}$ 
such that $\nu = \mu f^{-1}$ if and only if $f(X)$ is thick with respect to $\nu$.
Moreover, if the measure $\mu$ exists then it is unique.
\end{proposition}

\proof 
The condition is clearly necessary, since 
if $\nu = \mu f^{-1}$ and $F \in \mathcal{F}$ with $F \cap f(X) = \varnothing$ then $f^{-1}(F) = \varnothing$ and hence
$\nu(F) = \mu(f^{-1}(F)) = \mu(\varnothing) = 0$. Thus
suppose conversely that $\nu(F) = 0$ for all $F \in \mathcal{F}$ with $F \cap f(X) = \varnothing$.
Let $F_1,\, F_2 \in \mathcal{F}$ with $f^{-1}(F_1) = f^{-1}(F_2)$. Then 
$(F_1 \bigtriangleup F_2) \cap f(X) = \varnothing$
(with $F_1 \bigtriangleup F_2 = (F_1 \setminus F_2)\cup (F_2 \setminus F_1)$  the symmetric difference of $F_1$ and $F_2$),
hence $\nu(F_1 \bigtriangleup F_2) = 0$
and so $\nu(F_1) = \nu(F_2)$. Therefore, since $f^{-1}(\mathcal{F}) = \mathcal{E}$,
there exists a unique mapping $\mu : \mathcal{E} \to \Realpi$ such that
$\mu(f^{-1}(F)) = \nu(F)$ for all $F \in \mathcal{F}$, and it only remains to show that
$\mu$ is $\sigma$-additive.
Let $\{E_n\}_{n\ge 0}$ be a disjoint sequence from $\mathcal{E}$ and put $E = \bigcup_{n\ge 0} E_n$.
For each $n \ge 0$ let $F_n \in \mathcal{F}$ with $f^{-1}(F_n) = E_n$. 
Let $F'_0 = F_0$ and for each $n \ge 1$ put $F'_n = F_n \setminus \bigcup_{k=0}^{n-1} F_k$.
Then 
\[f^{-1}(F'_n) = f^{-1}\Bigl(F_n \setminus \bigcup_{k=0}^{n-1} F_k\Bigr)
   = f^{-1}(F_n) \setminus \bigcup_{k=0}^{n-1} f^{-1}(F_k) = E_n \setminus \bigcup_{k=0}^{n-1} E_k = E_n\]
and so $\nu(F'_n) = \mu(E_n) = \nu(F_n)$ for each $n \ge 0$.
But the sequence $\{F'_n\}_{n\ge 0}$ is disjoint and thus
\begin{eqnarray*} \sum_{n\ge 0} \mu(E_n) &=& \sum_{n\ge 0} \mu(f^{-1}(F_n)) = \sum_{n\ge 0} \nu(F_n) 
= \sum_{n\ge 0} \nu(F'_n) = \nu\Bigl(\bigcup_{n\ge 0} F'_n\Bigr)\\ 
&=& \nu\Bigl(\bigcup_{n\ge 0} F_n\Bigr) = \mu\Bigl(f^{-1}\Bigl(\bigcup_{n\ge 0} F_n\Bigr)\Bigr) 
= \mu\Bigl(\bigcup_{n\ge 0} f^{-1}(F_n)\Bigr) = \mu\Bigl(\bigcup_{n\ge 0} E_n\Bigr)\;. 
\end{eqnarray*}
Finally, the uniqueness follows because $\mu(f^{-1}(F)) = \nu(F)$ for all $F \in \mathcal{F}$ and $f^{-1}(\mathcal{F}) = \mathcal{E}$. \eop

Proposition~\ref{image-measures}.1 is most often (but not here) applied to the case in which
$f$ is a surjective exactly measurable mapping. For each measure $\nu$ on $\mathcal{E}$ there then exists a unique measure 
$\mu$ on $\mathcal{F}$ with $\nu = \mu f^{-1}$. 

If $\nu \in \Prob{Y}{\mathcal{F}}$ and $F \in \mathcal{F}$ then clearly $F$ is thick with respect to $\nu$ if and only if
$\nu(F) = 1$. The following generalisation of this fact is useful when dealing with type $\mathcal{A}$ spaces:

\begin{lemma}
Let $\nu \in \Prob{Y}{\mathcal{F}}$ and $F \in \mathcal{F}_*$. Then $F$ is thick with respect to $\nu$ if and only if
$\bar{\nu}(F) = 1$, where $\bar{\nu}$ is the unique extension of $\nu$ to $\mathcal{F}_*$.
\end{lemma}

\proof 
Clearly $F$ is thick with respect to $\nu$ if and only if $\mu_*(Y \setminus F) = 0$ and by Lemma~\ref{uni-meas}.1
$\bar{\nu}(F) = 1 - \bar{\nu}(Y \setminus F) = 1 - \mu_*(Y \setminus F)$. \eop

One of the main reasons for working with image measures in the space $(\specM,\mathcal{B})$ is the 
following wonderful property of finite measures on $(\specM,\mathcal{B})$:

\begin{proposition}
Any additive mapping $\mu : \mathcal{C}_{\specM} \to \Realpos$ defined on the algebra $\mathcal{C}_{\specM}$ of cylinder sets
is automatically  $\sigma$-additive and thus has a unique
extension to a measure on $\mathcal{B}$.
(Note that the mapping $\mu$ here is bounded since
$\mu(C) \le \mu(\specM) \in \Realpos$ for each $C \in \mathcal{C}_{\specM}$.)
\end{proposition}

\proof 
If $\{C_n\}_{n\ge 0}$ is a decreasing sequence from $\mathcal{C}_{\specM}$ with $\bigcap_{n\ge 0} C_n = \varnothing$ then,
since the elements of $\mathcal{C}$ are compact, there exists $m \ge 0$ so that $C_n = \varnothing$
for all $n \ge m$. Thus $\mu(C_n) = 0$ for all $n \ge m$ and hence $\lim_{n\to\infty} \mu(C_n) = 0$.
Therefore $\mu$ is $\sigma$-additive. \eop

Let us now illustrate our method by applying it to give a proof of part of the Dunford-Pettis theorem.
If $(X,\mathcal{E})$ is a measurable space then a subset $Q$ of $\Prob{X}{\mathcal{E}}$ is \definition{equicontinuous} if for each 
decreasing sequence $\{E_n\}_{n\ge 1}$ from $\mathcal{E}$
with $\bigcap_{n\ge 1} E_n = \varnothing$ and each $\varepsilon > 0$ there exists $p \ge 1$ so that 
$\mu(E_p) < \varepsilon$ for all $\mu \in Q$.

The following is the elementary (but more useful) half of the the Dunford-Pettis theorem.
(The proof of the converse can be found in Dunford and Schwartz \cite{dunfordsch}, Chapter IV.9.)

\begin{proposition}
Let $(X,\mathcal{E})$ be a measurable space and $Q$ be an 
equicontinuous subset of $\Prob{X}{\mathcal{E}}$.
Then for each sequence $\{\mu_n\}_{n\ge 1}$ from $Q$ there is a subsequence
$\{n_j\}_{j \ge 1}$ and a measure $\mu \in \Prob{X}{\mathcal{E}}$ such that $\mu(E) = \lim_{j} \mu_{n_j}(E)$
for all $E \in \mathcal{E}$. 
\end{proposition}

\proof
We apply the method outlined above to show that Proposition~\ref{image-measures}.3 holds for countably generated
measurable spaces and then use a standard technique to reduce the general case to the countably generated one.

Let us 
say that a measurable space $(Y,\mathcal{F})$ has the
\definition{weak sequential compactness property}
if whenever $Q$ is an equicontinuous subset of $\Prob{Y}{\mathcal{F}}$
then for each sequence $\{\mu_n\}_{n\ge 1}$ from $Q$ there exists a subsequence
$\{n_j\}_{j \ge 1}$ and  a measure $\mu \in \Prob{Y}{\mathcal{F}}$ such that 
$\mu(F) = \lim_{j} \mu_{n_j}(F)$
for all $F \in \mathcal{F}$.

\begin{lemma}
The space $(\specM,\mathcal{B})$ has the weak sequential compactness property.
\end{lemma}

\proof
Let $Q \subset \Prob{\specM}{\mathcal{B}}$ be equicontinuous   and $\{\mu_n\}_{n\ge 1}$ be a sequence from 
$Q$. Then, since the algebra $\mathcal{C}_{\specM}$ of cylinder sets is countable and the values $\mu_n(C)$ all lie in the compact interval $[0,1]$
the usual diagonal argument implies there exists a subsequence
$\{n_j\}_{j\ge 1}$ and $\nu : \mathcal{C}_{\specM} \to \Realpos$ so that
$\nu(C) = \lim_{j} \mu_{n_j}(C)$ for all $C \in \mathcal{C}_{\specM}$.
But $\nu$ is clearly additive and $\nu(M) = 1$ and hence by Proposition~\ref{image-measures}.2 there exists $\mu \in \Prob{\specM}{\mathcal{B}}$
with $\mu(C) = \nu(C)$ for all $C \in \mathcal{C}_{\specM}$; thus $\mu(C) = \lim_{j} \mu_{n_j}(C)$ for all 
$C \in \mathcal{C}_{\specM}$. Now let
\[\mathcal{K} = \bigl\{ B \in \mathcal{B} : \mu(B) = \lim_{j\to\infty} \mu_{n_j}(B) \bigr\}\;;\]
then $\mathcal{C} \subset \mathcal{K}$ and $\sigma(\mathcal{C}) = \mathcal{B}$, and so by the monotone class theorem
it is enough to show that $\mathcal{K}$ is a monotone class.
Let $\{B_n\}_{n\ge 1}$ be an increasing sequence from $\mathcal{K}$ and put $B = \bigcup_{n\ge 1} B_n$.
For each $p \ge 1$ let $A_p = B \setminus  B_p$; then $\{A_p\}_{p \ge 1}$ is a decreasing sequence from $\mathcal{B}$ with 
$\bigcap_{p \ge 1} A_p = \varnothing$.
Let $\varepsilon > 0$; there thus exists $p \ge 1$ so that $\mu(A_p) < \varepsilon/3$ and
so that $\omega(A_p) < \varepsilon/3$ for all $\omega \in Q$. 
Moreover, since $B_p \in \mathcal{K}$, there exists $m \ge 1$
so that $|\mu(B_p) - \mu_{n_j}(B_p)| < \varepsilon/3$ for all $j \ge m$. Hence
\[ |\mu(B) - \mu_{n_j}(B)| \le |\mu(B_p) - \mu_{n_j}(B_p)| + \mu(A_p) + \mu_{n_j}(A_p) < \varepsilon\]
for all $j \ge m$, and so $\mu(B) = \lim_{j} \mu_{n_j}(B)$, i.e., $B \in \mathcal{K}$.
The case of a decreasing sequence from $\mathcal{K}$ is almost exactly the same. \eop

\begin{lemma}
Let $(X,\mathcal{E})$, $(Y,\mathcal{F})$ be measurable spaces and
$f : (X,\mathcal{E}) \to (Y,\mathcal{F})$ be an exactly measurable mapping.
If $(Y,\mathcal{F})$ has the weak sequential compactness property then so does $(X,\mathcal{E})$. 
\end{lemma}

\proof 
For each $\mu \in \Prob{X}{\mathcal{E}}$ denote the image measure $\mu f^{-1}  \in \Prob{Y}{\mathcal{F}}$
by $\mu'$. Let $Q \subset \Prob{X}{\mathcal{E}}$ be equicontinuous;
then the subset $Q' = \{ \mu' : \mu \in Q \}$ of $\Prob{Y}{\mathcal{F}}$ is also equicontinuous:
If $\{F_n\}_{n\ge 1}$ is a decreasing sequence from $\mathcal{F}$ with $\bigcap_{n\ge 1} F_n = \varnothing$ 
and $E_n = f^{-1}(F_n)$ for each $n \ge 1$ then $\{E_n\}_{n\ge 1}$ is a decreasing sequence 
from $\mathcal{E}$ with $\bigcap_{n\ge 1} E_n = \varnothing$;
thus, given $\varepsilon > 0$, there exists $p \ge 1$ so that $\mu(E_p) < \varepsilon$ for all $\mu \in Q$
and hence $\mu'(F_p) = \mu(E_p) < \varepsilon$ for all $\mu' \in Q'$.
Let $\{\nu\}_{n\ge 1}$ be a sequence from $Q$; thus $\{\nu'\}_{n\ge 1}$ is a sequence from $Q'$, and so if
$(Y,\mathcal{F})$ has the weak sequential compactness property then there exists a subsequence $\{n_j\}_{j\ge 1}$
and $\nu \in \Prob{Y}{\mathcal{F}}$ such that $\nu(F) = \lim_{j} \mu'_{n_j}(F)$ for all
$F \in \mathcal{F}$. In particular, if $F \in \mathcal{F}$ with $F \cap f(X) = \varnothing$ then
\[ \nu(F) = \lim_{j\to \infty} \mu'_{n_j}(F) = \lim_{j\to \infty} \mu_{n_j}(f^{-1}(F))
= \lim_{j\to \infty} \mu_{n_j}(\varnothing) = 0\] 
and so by Proposition~\ref{image-measures}.1 there exists a  $\mu \in \Prob{X}{\mathcal{E}}$ with
$\nu = \mu f^{-1}$. Let $E \in \mathcal{E}$; then $E = f^{-1}(F)$ for some $F \in \mathcal{F}$ and therefore
\[ \mu(E) = \mu(f^{-1}(F)) = \nu(F) = \lim_{j\to\infty} \mu'_{n_j}(F) 
= \lim_{j\to\infty} \mu_{n_j}(f^{-1}(F)) = \lim_{j\to\infty} \mu_{n_j}(E)\;.\]
This shows that $(X,\mathcal{E})$ has the weak sequential compactness property. \eop

\begin{lemma}
Each countably generated measurable space has the weak sequential compactness property.
\end{lemma}

\proof
This follows immediately from Proposition~\ref{cg-meas-spaces}.2 together with
Lemmas \ref{image-measures}.2 and \ref{image-measures}.3.
\eop

We turn to the general case, so now let $(X,\mathcal{E})$ be an arbitrary measurable space.
Things will be reduced to the countably generated case by modifying the proof of a similar reduction to be found in
Dunford and Schwartz \cite{dunfordsch}, Chapter IV.9.

\begin{lemma}
If $\{\mu_n\}_{n\ge 1}$ is a sequence from $\Prob{X}{\mathcal{E}}$ then there is a $\nu \in \Prob{X}{\mathcal{E}}$
such that $\mu_n \ll \nu$ for all $n \ge 1$.
\end{lemma}

\proof
Just take, for example, $\nu = \sum_{n\ge 1} 2^{-n} \mu_n$.
\eop

For each sub-$\sigma$-algebra $\mathcal{F}$ of $\mathcal{E}$ let
$\Mappi(\mathcal{F})$ denote the set of measurable mappings $g : (X,\mathcal{F}) \to (\Realpos,\mathcal{B}_{\Realpos})$ 
with $\mathcal{B}_{\Realpos}$ the $\sigma$-algebra of Borel subsets of $\Realpos$.

Let $Q$ be an equicontinuous subset of $\Prob{X}{\mathcal{E}}$ and let $\{\mu_n\}_{n\ge 1}$ be a sequence from $Q$.
By Lemma~\ref{image-measures}.5 there exists $\nu \in \Prob{X}{\mathcal{E}}$ such that $\mu_n \ll \nu$ for each $n \ge 1$;
by the Radon-Nikodym theorem there  then exists $h_n \in \Mappi(\mathcal{E})$ such that 
$\mu_n = \int_E h_n\, d\mu$ for all $E \in \mathcal{E}$.
Now consider any countably generated $\sigma$-algebra $\mathcal{F} \subset \mathcal{E}$.
For each $\mu \in \Prob{X}{\mathcal{E}}$ let $\mu' \in \Prob{Y}{\mathcal{F}}$ be the restriction of $\mu$ to
$\mathcal{F}$. Then $Q' = \{ \mu' : \mu \in Q \}$ is an equicontinuous subset of $\Prob{Y}{\mathcal{F}}$
and $\{\mu'_n\}_{n\ge 1}$ is a sequence from $Q'$.
Thus, applying Lemma~\ref{image-measures}.4 to $(X,\mathcal{F})$, there exists a subsequence
$\{n_j\}_{j \ge 1}$ and $\omega \in \Prob{Y}{\mathcal{F}}$ with $\omega(F) = \lim_{j} \mu'_{n_j}(F)$
for all $F \in \mathcal{F}$. 
But if $F \in \mathcal{F}$ with $\nu'(F) = 0$ then $\nu(F) = 0$, thus
$\omega(F)  = \lim_j \mu'_{n_j}(F) = 0$ and hence $\omega \ll \nu'$.
Therefore by the Radon-Nikodym theorem there exists $h \in \Mappi(\mathcal{F})$ such that
$\omega(F) = \int_F h\,d\nu'$ for all $F \in \mathcal{F}$.
Now define $\mu \in \Prob{X}{\mathcal{E}}$ by letting $\mu(E) = \int_E h\,d\nu$ for all $E \in \mathcal{E}$.
In particular $\mu(F) = \lim_{j} \mu_{n_j}(F)$ for all $F \in \mathcal{F}$ and hence also
$\int g \,d\mu = \lim_{j} \int g\,d\mu_{n_j}$ for all bounded elements of $\Mappi(\mathcal{F})$.

\begin{lemma}
Suppose $h_n \in \Mappi(\mathcal{F})$ for each $n \ge 1$. Then
$\mu(E) = \lim_{j} \mu_{n_j}(E)$ for all $E \in \mathcal{E}$.
\end{lemma}

\proof 
Let $E \in \mathcal{E}$; then $I_E \le 1$ and therefore there exists $g \in \Mappi(\mathcal{F})$
with $g \le 1$ such that  $\int gf\,d\nu = \int_E f\,d\nu$ for all $f \in \Mappi(\mathcal{F})$ (i.e.,
$g$ is the conditional expectation of $I_E$ with respect to $\mathcal{F}$).
Then 
\[\mu(E) =  \int_E h\,d\nu = \int gh\,d\nu  = \int g\,d\mu\] 
and $\mu_n(E) = \int_E h_n\,d\nu = \int gh_n\,d\nu = \int g\,d\mu_n$ 
for each $n \ge 1$ and therefore 
\[ \mu(E) = \int g\,d\,\mu = \lim_{j\to\infty} \int g\,d\mu_{n_j} = \lim_{j\to\infty} \mu_{n_j}(E)\;.\ \eop \]

The next result completes the proof of Proposition~\ref{image-measures}.3.

\begin{lemma}
There exists a countably generated sub-$\sigma$-algebra $\mathcal{F}$ of $\mathcal{E}$ such that
$h_n \in \Mappi(\mathcal{F})$ for each $n \ge 1$. 
\end{lemma}

\proof
Let $\mathcal{J}$ be the countable set consisting  of all elements of $\mathcal{E}$ of the form
$\{ x \in X : h_n(x) > r \}$ with $n \ge 1$ and $r \in \Rat^+$, and put
$\mathcal{F} = \sigma(\mathcal{J})$. Then $\mathcal{F} \subset \mathcal{E}$ and $\mathcal{F}$ is  countably
generated. But if $a \in \Realpos$ then there is a decreasing sequence $\{r_m\}_{m\ge 1}$ from $\Rat^+$ with
$\lim_{m} r_m = a$ and thus
\[ \{ x \in X : h_n(x) > a \} = \bigcap_{m\ge 1} \{ x \in X : h_n(x) > r_m \} \in \mathcal{F}\;. \]
This implies that $h_n \in \Mappi(\mathcal{F})$ for each $n \ge 1$. \eop


\startsection{The Kolmogorov extension property}
\label{kolmogorov}

We consider the following set-up: 

\begin{evlist}{25pt}{5pt}
\item[(1)] For each $n \ge 0$ there is a countably generated measurable space $(X_n,\mathcal{E}_n)$ and a measurable mapping
$i_n : (X_{n+1},\mathcal{E}_{n+1}) \to (X_n,\mathcal{E}_n)$.

\item[(2)] There is also a measurable space $(X,\mathcal{E})$ and for each $n \ge 0$ a surjective measurable mapping 
$\tau_n : (X,\mathcal{E}) \to (X_n,\mathcal{E}_n)$ with $i_n \circ \tau_{n+1} = \tau_n$ for all $n \ge 0$.

\item[(3)] Since $i_n \circ \tau_{n+1} = \tau_n$ it follows that 
$\tau_n^{-1}(\mathcal{E}_n) = \tau_{n+1}^{-1}(i_n^{-1}(\mathcal{E}_n)) \subset \tau_{n+1}^{-1}(\mathcal{E}_{n+1})$ 
for each $n \ge 0$ and so $\{\tau_n^{-1}(\mathcal{E}_n)\}_{n\ge 0}$ is an increasing sequence of sub-$\sigma$-algebras of $\mathcal{E}$.
We assume that $\mathcal{E} = \sigma\big(\bigcup_{n \ge 0} \tau_n^{-1}(\mathcal{E}_n)\bigr)$.

\item[(4)] Finally, we also assume the following: For each sequence $\{A_n\}_{n \ge 0}$ of atoms with
$A_n \in \mathsf{A}(\mathcal{E}_n)$ such that $A_{n+1} \subset i_n^{-1}(A_n)$ for each $n \ge 0$ there exists
an element $x \in X$ with $\tau_n(x) \in A_n$ for all $n \ge 0$.
\end{evlist}

Note that (3) implies $(X,\mathcal{E})$ is also countably generated.

\begin{proposition}
Let $\mathcal{D}$ be a classifying class closed under countable products and suppose that $(X_n,\mathcal{E}_n)$ is a type $\mathcal{D}$
space for each $n \ge 0$. Then $(X,\mathcal{E})$ is also a type $\mathcal{D}$ space.
\end{proposition}

\begin{theorem}
Let $(X_n,\mathcal{E}_n)$ be  a type $\mathcal{A}$ space for each $n \ge 0$ (and so $(X,\mathcal{E})$
is also a type $\mathcal{A}$ space). 
For each $n \ge 0$ let $\mu_n \in \Prob{X_n}{\mathcal{E}_n}$ and suppose the sequence of measures $\{\mu_n\}_{n\ge 0}$
is consistent in that $\mu_{n+1}i_n^{-1} = \mu_n$ for each $n \ge 0$. 
Then there exists a unique measure $\mu \in \Prob{X}{\mathcal{E}}$ such that $\mu_n \tau_n^{-1} = \mu$ for all $n \ge 0$.
\end{theorem}

Before beginning the proofs of these two results we look at the usual form in which they are applied.
Let $(Y,\mathcal{F})$ be a measurable space and $\{\mathcal{F}_n\}_{n\ge 0}$ be an increasing sequence 
of countably generated sub-$\sigma$-algebras of $\mathcal{F}$ with $\mathcal{F} = \sigma(\bigcup_{n\ge 0} \mathcal{F}_n)$.
A sequence of measures $\{\mu_n\}_{n\ge 0}$ with
$\mu_n \in \Prob{Y}{\mathcal{F}_n}$ for each $n \ge 0$ 
is  \definition{consistent} 
if $\mu_n(F) = \mu_{n+1}(F)$ for all $F \in \mathcal{F}_n$, $n \ge 0$.
The sequence $\{\mathcal{F}_n\}_{n\ge 0}$ is said to have the 
\definition{Kolmogorov extension property}
if for each consistent sequence $\{\mu_n\}_{n\ge 0}$ there exists $\mu \in \Prob{Y}{\mathcal{F}}$
such that $\mu(F) = \mu_n(F)$ for all $F \in \mathcal{F}_n$, $n \ge 0$.
(This measure $\mu$ is then unique, since it is uniquely determined by the sequence $\{\mu_n\}_{n\ge 0}$
on the algebra $\mathcal{G} = \bigcup_{n\ge 0} \mathcal{F}_n$ and $\sigma(\mathcal{G}) = \mathcal{F}$.)
Finally, the $\sigma$-algebra $\mathcal{F}$ is called the
\definition{inverse limit} of the sequence $\{\mathcal{F}_n\}_{n\ge 0}$ if
$\bigcap_{n\ge 0} A_n \ne \varnothing$ holds 
whenever $\{A_n\}_{n\ge 0}$ is a decreasing sequence of atoms with $A_n \in \mathrm{A}(\mathcal{F}_n)$
for each $n \ge 0$.

\begin{theorem}
Let $(Y,\mathcal{F}_n)$ be a type $\mathcal{A}$ space for each $n \ge 0$ 
and let $\mathcal{F}$ be the inverse limit of the sequence $\{\mathcal{F}_n\}_{n\ge 0}$.
Then $(Y,\mathcal{F})$ is also a type $\mathcal{A}$ space and the
sequence $\{\mathcal{F}_n\}_{n\ge 0}$ has the Kolmogorov extension property.
\end{theorem}

\proof 
This follows immediately from Proposition~\ref{kolmogorov}.1 and Theorem~\ref{kolmogorov}.1 
with $(X,\mathcal{E}) = (Y,\mathcal{F})$,
$(X_n,\mathcal{E}_n) = (Y,\mathcal{F}_n)$ and $i_n = \tau_n = \id_Y$ for all $n \ge 0$.
\eop

Theorem~\ref{kolmogorov}.2 is the form of the Kolmogorov extension theorem occurring 
in Chapter V of Parthasarathy \cite{partha}.

We start the preparations for the proofs of Proposition~\ref{kolmogorov}.1 and Theorem~\ref{kolmogorov}.1;
they are based on the proof of Theorem~\ref{kolmogorov}.2 in \cite{partha}.

For each $n \ge 0$ let $f_n : (X_n,\mathcal{E}_n) \to (\specM,\mathcal{B})$ be an exactly measurable mapping and put
$q_n = f_n \circ \tau_n$; thus $q_n : (X,\mathcal{E}) \to (\specM,\mathcal{B})$ is measurable (but usually it will not be
exactly measurable).

Consider the product $\specM^{\Nat}$ as a compact metric space in the usual way; then the
product $\sigma$-algebra $\mathcal{ B}^{\Nat}$ is also the Borel $\sigma$-algebra. Define a mapping
$q : X \to \specM^\Nat$ by letting $q(x) = \{q_n(x)\}_{n \ge 0}$ for each $x \in X$.

\begin{lemma}
The mapping $q : (X,\mathcal{E}) \to (\specM^\Nat,\mathcal{B}^\Nat)$ is exactly measurable. 
\end{lemma}

\proof 
If $B_n \in \mathcal{B}$ for each $n \ge 0$ then $q^{-1}(\prod_{n\ge 0} B_n) = \bigcap_{n\ge 0} q_n^{-1}(B_n) \in \mathcal{E}$
and this implies that $q$ is measurable. Now fix $m \ge 0$ and let $E \in \mathcal{E}_m$; there thus exists $B \in \mathcal{B}$ with
$f_m^{-1}(B) = E$. Let $B' = \{ \{z_n\}_{n \ge 0} \in \specM^\Nat : z_m \in B \}$; then $B' \in \mathcal{B}^\Nat$ and
\[ q^{-1}(B') = q_m^{-1}(B) = \tau_m^{-1}(f_m^{-1}(B)) = \tau_m^{-1}(E)\;,\]
which shows that $\tau_m^{-1}(\mathcal{E}_m) \subset q^{-1}(\mathcal{B}^\Nat)$. Therefore
$\bigcup_{n\ge 0} \tau_n^{-1}(\mathcal{E}_n) \subset q^{-1}(\mathcal{B}^\Nat)$ and so $\mathcal{E} \subset q^{-1}(\mathcal{B}^\Nat)$.
Hence $q^{-1}(\mathcal{B}^\Nat) = \mathcal{E}$, i.e., $q$ is exactly measurable.
\eop

As usual let $\Delta : \specM \to \specM$ be the mapping given by
\[\Delta(\{z_n\}_{n\ge 0}) = \{z'_n\}_{n\ge 0}\;,\]
where $z'_n = z_{2n+1}$ for each $n \ge 0$; thus $\Delta$ is continuous and surjective. Also let
\[\specM_\Delta = \{ \{z_n\}_{n\ge 0} \in \specM^{\Nat} : z_n = \Delta(z_{n+1})\ \mbox{for all}\ n \ge 0 \}\;;\]
$\specM_\Delta$ is a closed (and thus compact) subset of $\specM^{\Nat}$.
The $\sigma$-algebra of Borel subsets of $\specM_\Delta$ will be denoted by $\mathcal{B}_{\Delta}$, thus
$\mathcal{B}_{\Delta}$ is also the trace $\sigma$-algebra $\mathcal{B}^\Nat_{\,|\specM_{\Delta}}$.
The next result is given in a form which will be needed later but we present it here because it also establishes that
$\specM_\Delta$ is non-empty.

\begin{lemma}
Let $m \ge 0$ and $w \in \specM$.
Then there exists $z = \{z_n\}_{n \ge 0} \in \specM_\Delta$ with $z_m = w$, and so in particular $\specM_\Delta$ is non-empty.
\end{lemma}

\proof 
Define an element $z = \{z_n\}_{n \ge 0}$ of $\specM^\Nat$
as follows: Set $z_m = w$ and let $z_{m-1},\,\ldots,\,z_0 \in \specM$ be given (uniquely) by the requirement that
$z_n = \Delta(z_{n+1})$ for $n = m - 1,\,\ldots,\,0$. Now choose $z_n$ for $n \ge m+1$ inductively so that
$\Delta(z_n) = z_{n-1}$  for $n = m+1,\,m+2,\,\ldots\,$. This can be done, since $\Delta$ is surjective, although at each stage
the choice is never unique. (It is possible, however, to make an explicit choice, for example by taking the even components of $z_n$
to be $0$.) Then $z \in \specM_\Delta$ and $z_m = w$.
\eop

Now let $\mathcal{D}$ be a classifying class closed under finite products,
and assume that $(X_n,\mathcal{E}_n)$ is a type $\mathcal{D}$ space for each $n \ge 0$.

\begin{proposition}
For each $n \ge 0$ there exists an exactly measurable mapping $f_n : (X_n,\mathcal{E}_n) \to (\specM,\mathcal{B})$ 
with $f_n(X_n) \in \mathcal{D}$ and such that $f_n \circ i_n = \Delta \circ f_{n+1}$ for all $n \ge 0$. 
\end{proposition}

\proof 
There exists an exactly measurable mapping $f_0 : (X_0,\mathcal{E}_0) \to (\specM,\mathcal{B})$ 
with $f_0(X_0) \in \mathcal{D}$, since $(X_0,\mathcal{E}_0)$ is a type $\mathcal{D}$ space.
Suppose for some $n \ge 0$ we have 
exactly measurable mappings $f_k : (X_k,\mathcal{E}_k) \to (\specM,\mathcal{B})$, $k = 0,\,\ldots,\,n$, 
with $f_k(X_k) \in \mathcal{D}$ and 
such that $f_k \circ i_k = \Delta \circ f_{k+1}$ for $k = 0,\,\ldots,\,n-1$.
Then, applying Proposition~\ref{cl-classes}.7 with $h = f_n$ and $f = i_n$, there exists an
exactly measurable mapping $f_{n+1} : (X_{n+1},\mathcal{E}_{n+1}) \to (\specM,\mathcal{B})$ 
with $f_{n+1}(X_{n+1}) \in \mathcal{D}$ and such that
$f_n \circ i_n = \Delta \circ f_{n+1}$. The result therefore follows by induction. \eop

From now on suppose that the mappings $f_n$, $n \ge 0$, have been chosen as in Proposition~\ref{kolmogorov}.2.

\begin{lemma}
$q(X) = \specM_\Delta \cap \prod_{n\ge 0} f_n(X_n)$.
\end{lemma}

\proof 
Note that $f_n(X_n) = q_n(X)$, since $\tau_n$ is surjective.
Now if $x \in X$ then $q(x) = \{q_n(x)\}_{n \ge 0}$ and $q_n(x) \in q_n(X)$; thus
$q(x) \in \prod_{n\ge 0} q_n(X) = \prod_{n\ge 0} f_n(X_n)$.
Moreover, $q_n(x) = \Delta(q_{n+1}(x))$ for all $n \ge 0$, since
\[ \Delta \circ q_{n+1} = \Delta \circ f_{n+1} \circ \tau_{n+1} = f_n \circ i_n \circ \tau_{n+1} 
     = f_n \circ \tau_n = q_n\;, \]
and hence $q(x) \in \specM_\Delta$. This shows that $q(X) \subset \specM_\Delta \cap \prod_{n\ge 0} f_n(X_n)$.

Conversely, consider $z = \{z_n\}_{n\ge 0} \in \specM_\Delta \cap \prod_{n\ge 0} f_n(X_n)$; then
$z_n \in f_n(X_n)$ and $z_n = \Delta(z_{n+1})$ for each $n \in 0$.
Put $A_n = f_n^{-1}(\{z_n\})$; by Lemma~\ref{cg-meas-spaces}.1 $A_n \in \mathrm{A}(\mathcal{E}_n)$ and, since $z_n = \Delta(z_{n+1})$,
\[A_{n+1} = f_{n+1}^{-1}(\{z_{n+1}\}) \subset  f_{n+1}^{-1}(\Delta^{-1}(\{z_{n+1}\})) 
= i_n^{-1}(f_n^{-1}(\{z_n\})) = i_n^{-1}(A_n)\]
for each $n \ge 0$. Therefore by assumption there exists an element $x \in X$ with $\tau_n(x) \in A_n$ for all $n \ge 0$
and then $q_n(x) = f_n(\tau_n) \in f_n(A_n) = \{z_n\}$, i.e., $q_n(x) = z_n$ for all $n \ge 0$.
Hence $q(x) = z$, which shows that
$q(X) \supset \specM_\Delta \cap \prod_{n\ge 0} f_n(X_n)$. \eop

\textit{Proof of Proposition~\ref{kolmogorov}.1:} 
(We are here assuming that $\mathcal{D}$ is closed under countable products.)
If $h : \specM^\Nat \to \specM$ is a homeomorphism then by Lemma~\ref{kolmogorov}.3
\[ h(q(X)) = h\Bigl(\specM_\Delta \cap \prod_{n\ge 0} f_n(X_n)\Bigr) = h(\specM_\Delta) \cap h\Bigl(\prod_{n\ge 0} f_n(X_n)\Bigr)\]
is an element of $\mathcal{D}$, since $h(\specM_\Delta)$ is a compact subset of $\specM$ and thus in $\mathcal{B}$.
Put $g = h \circ q$; then by Lemma~\ref{e-meas-maps}.1~(2) $g : (X,\mathcal{E}) \to (\specM,\mathcal{B})$ is exactly measurable
and $g(X) = h(q(X)) \in \mathcal{D}$. Hence $(X,\mathcal{E})$ is a type $\mathcal{D}$ space. \eop

For each $m \ge 0$ let $\pi_m : \specM^\Nat \to \specM$ be the projection mapping with
$\pi_m(z) = z_m$ for each $z = \{z_n\}_{n\ge 0} \in \specM^\Nat$ and let $\theta_m : \specM_\Delta \to \specM$ be the restriction
of $\pi_m$ to $\specM_\Delta$. Then $\theta_m$ is continuous (since $\pi_m$ is)
and Lemma~\ref{kolmogorov}.2 implies that $\theta_m$ is surjective.
Moreover, since $\Delta(z_{n+1}) = z_n$ for each $z = \{z_n\}_{n\ge 0} \in \specM_\Delta$ it follows that
$\Delta \circ \theta_{n+1} = \theta_n$ for each $n \ge 0$.

\begin{lemma}
$\{\theta_n^{-1}(\mathcal{B})\}_{n \ge 0}$ is an increasing sequence of sub-$\sigma$-algebras of $\mathcal{B}_\Delta$ with
\[  \mathcal{B}_\Delta = \sigma \Bigl(\bigcup_{n \ge 0} \theta_n^{-1}(\mathcal{B})\Bigr)\;. \]
\end{lemma}

\proof
If $B \in \mathcal{B}$ then $\theta_n^{-1}(B) = \pi_n^{-1}(B) \cap \specM_\Delta \in \mathcal{B}_\Delta$, since 
$\mathcal{B}_\Delta = \mathcal{B}^\Nat_{\,|\specM_\Delta}$, and thus $\theta_n^{-1}(\mathcal{B})$ is a sub-$\sigma$-algebra of
$\mathcal{B}_\Delta$.
Moreover, $\theta_n^{-1}(\mathcal{B}) = \theta_{n+1}^{-1}(\Delta^{-1}(\mathcal{B})) \subset \theta_{n+1}^{-1}(\mathcal{B})$
for each $n \ge 0$ (since $\Delta \circ \theta_{n+1} = \theta_n$) and therefore
$\{\theta_n^{-1}(\mathcal{B})\}_{n \ge 0}$ is an increasing sequence of sub-$\sigma$-algebras of $\mathcal{B}_\Delta$.
Now let $\mathcal{Z}$ denote the set of all subsets of $\specM^\Nat$ having the form $\pi_n^{-1}(B)$ for some $n \ge 0$ and some $B \in \mathcal{B}$.
Then $\mathcal{B}^\Nat = \sigma(\mathcal{Z})$, and since
$\theta_n^{-1}(B) = \pi_n^{-1}(B) \cap \specM_\Delta$ it follows that
$\bigcup_{n \ge 0} \theta_n^{-1}(\mathcal{B}) = \mathcal{Z}_{|\specM_\Delta}$ (with $\mathcal{Z}_{|\specM_\Delta}$  the set of all sets of the form
$Z \cap \specM_\Delta$ with $Z \in \mathcal{Z}$). Thus by Lemma~\ref{cl-classes}.5
\[ \sigma \Bigl(\bigcup_{n \ge 0} \theta_n^{-1}(\mathcal{B})\Bigr) 
= \sigma(\mathcal{Z}_{|\specM_\Delta}) = \sigma(\mathcal{Z})_{|\specM_\Delta} = \mathcal{B}^\Nat_{\,\specM_\Delta}
= \mathcal{B}_\Delta\;.\ \eop \]

\begin{lemma}
Let $\{\nu_n\}_{n\ge 0}$ be a sequence from $\Prob{\specM}{\mathcal{B}}$ with $\nu_n = \nu_{n+1} \Delta^{-1}$
for each $n \ge 0$. Then there exists $\nu \in \Prob{\specM_\Delta}{\mathcal{B}_\Delta}$
with $\nu_n = \nu \theta_n^{-1}$ for all $n \ge 0$.
\end{lemma}

\proof 
Let $m \ge 0$ and $w_0,\,\ldots,\,w_m \in \{0,1\}$; then
\[ \Delta^{-1}(\specM(w_0,\ldots,w_m)) = \bigcup_{z_0,\ldots,z_m \in \{0,1\}} \specM(z_0,w_0,\ldots,z_m,w_m) \]
and hence $\Delta^{-1}(\mathcal{C}_{\specM}) \subset \mathcal{C}_{\specM}$.
Therefore $\theta_n^{-1}(\mathcal{C}_{\specM}) = \theta_{n+1}^{-1}(\Delta^{-1}(\mathcal{C}_{\specM})) \subset \theta_{n+1}^{-1}(\mathcal{C}_{\specM})$
for each $n \ge 0$, and so $\{\theta_n^{-1}(\mathcal{C}_{\specM})\}_{n\ge 0}$ is an increasing sequence of countable algebras.
Put $\mathcal{C}_\Delta = \bigcup_{n\ge 0} \theta_n^{-1}(\mathcal{C}_{\specM})$;
then $\mathcal{C}_\Delta$ is a countable algebra and by Lemma~\ref{kolmogorov}.4 
\[ \sigma(\mathcal{C}_\Delta) = \sigma\Bigl(\bigcup_{n\ge 0} \theta_n^{-1}(\mathcal{C}_{\specM})\Bigr)
= \sigma\Bigl(\bigcup_{n\ge 1} \theta_n^{-1}(\sigma(\mathcal{C}_{\specM}))\Bigr)
= \sigma\Bigl(\bigcup_{n\ge 1} \theta_n^{-1}(\mathcal{B})\Bigr) = \mathcal{B}_\Delta\;.\]
Moreover, each element of $\mathcal{C}_\Delta$ is compact (since the mappings $\theta_n$ are continuous and
$\specM_\Delta$ is compact), and hence $\mathcal{C}_\Delta$ has the finite intersection property.

Let $n \ge 0$; by Lemma~\ref{kolmogorov}.3 $\theta_n$ is surjective and thus Proposition~\ref{image-measures}.1 
(applied to the exactly measurable mapping $\theta_n : (\specM_\Delta,\theta_n^{-1}(\mathcal{B})) \to (\specM,\mathcal{B})$)
implies there is a unique $\nu'_n \in \Prob{\specM_\Delta}{\theta_n^{-1}(\mathcal{B})}$ with $\nu_n = \nu'_n \theta_n^{-1}$,
and the sequence $\{\nu'_n\}_{n\ge 1}$ is consistent in that $\nu'_{n+1}(D) = \nu'_n(D)$ for all 
$D \in \theta_n^{-1}(\mathcal{B})$, $n \ge 0$. (Let $D \in \theta_n^{-1}(\mathcal{B})$
with $D = \theta_n^{-1}(B)$; then $D = \theta_{n+1}^{-1}(\Delta^{-1}(B))$, since $\theta_n = \Delta \circ \theta_{n+1}$
and hence $\nu'_{n+1}(D) = \nu_{n+1}(\Delta^{-1}(B)) = \nu_n(B) = \nu'_n(D)$.)
There is thus a unique mapping $\nu' : \bigcup_{n\ge 0} \theta_n^{-1}(\mathcal{B})  \to \Realpos$
such that $\nu'(D) = \nu'_n(D)$ for all $D \in \theta_n^{-1}(\mathcal{B})$, $n \ge 0$, and it is clear that
$\nu'$ is finitely additive. Now the restriction of $\nu'$ to $\mathcal{C}_\Delta$
is also finitely additive and hence (as in the proof of Proposition~\ref{image-measures}.2) 
there exists a unique $\nu \in \Prob{\specM_\Delta}{\mathcal{B}_\Delta}$
with $\nu(D) = \nu'(D)$ for all $D \in \mathcal{C}_\Delta$.
But then the restriction of $\nu$ to $\theta_n^{-1}(\mathcal{B})$ is a probability measure which is an extension of the
restriction of $\nu'_n$ to $\theta_n^{-1}(\mathcal{C}_{\specM})$, and $\theta_n^{-1}(\mathcal{C}_{\specM})$ is an algebra with
$\sigma(\theta_n^{-1}(\mathcal{C}_{\specM})) = \theta_n^{-1}(\mathcal{B})$.
This means that $\nu$ is an extension of $\nu'_n$ and from this it immediately follows that $\nu_n = \nu \theta_n^{-1}$
for each $n \ge 0$. \eop

\begin{lemma}
For each $n \ge 0$ let $\mu_n \in \Prob{X_n}{\mathcal{E}_n}$ and suppose the sequence $\{\mu_n\}_{n\ge 0}$
is consistent in that $\mu_{n+1} i_n^{-1} = \mu_n$ for each $n \ge 0$.
For each $n \ge 0$ let $\nu_n \in \Prob{\specM}{\mathcal{B}}$ be the image measure $\mu_n f_n^{-1}$. 
Then $\nu_{n+1}\Delta^{-1} = \nu_n \Delta^{-1}$ for all $n \ge 0$.
\end{lemma}

\proof
Let $n \ge 0$; then for all $B \in \mathcal{B}$
\begin{eqnarray*}  
\nu_{n+1}(\Delta^{-1}(B)) &=& \mu_{n+1}(f_{n+1}^{-1}(\Delta^{-1}(B)))\\ 
&=& \mu_{n+1}(i_n^{-1}(f_n^{-1}(B))) = \mu_n(f_n^{-1}(B)) = \nu_n(B)
\end{eqnarray*}
and therefore $\nu_{n+1}\Delta^{-1} = \nu_n$. \eop

\textit{Proof of Theorem~\ref{kolmogorov}.1:}
We are here assuming that each $(X_n,\mathcal{E}_n)$ is a type $\mathcal{A}$ space and hence the mappings $f_n$, $n \ge 0$,
in Proposition~\ref{kolmogorov}.2 can be chosen so that $f_n(X_n) \in \mathcal{A}$ for all $n \ge 0$.
This means that $f_n(X_n) \in \mathcal{B}_*$ for all $n \ge 0$, since by
Theorem~\ref{uni-meas}.1 $\mathcal{A} \subset \mathcal{B}_*$.

By Lemma~\ref{kolmogorov}.3 $q(X) \subset \specM_\Delta$ and so we will consider $q$ as a mapping from $X$ to $\specM_\Delta$;
by Lemma~\ref{kolmogorov}.1 and Proposition~\ref{e-meas-maps}.2 the mapping
$q : (X,\mathcal{E}) \to (\specM_\Delta,\mathcal{B}_\Delta)$ is exactly measurable; 
also $\theta_n \circ q = f_n \circ \tau_n$ for each $n \ge 0$.
Note that if $\{C_n\}_{n \ge 0}$ is a sequence of subsets of $\specM$ then 
\[\specM_\Delta \cap \prod_{n\ge 0} C_n  = \specM_\Delta \cap \bigcap_{n \ge 0} \pi_n^{-1}(C_n)
= \bigcap_{n \ge 0} \theta_n^{-1}(C_n)\]
and hence by Lemma~\ref{kolmogorov}.3 $q(X) = \bigcap_{n \ge 0} \theta_n^{-1}(f_n(X_n))$. 
Moreover, $q(X) \in (\mathcal{B}_\Delta)_*$, since by Lemma~\ref{uni-meas}.5
$\theta_n^{-1}(\mathcal{B}_*) \subset (\mathcal{B}_\Delta)_*$.

Now let $\{\mu_n\}_{n\ge 0}$ be a consistent sequence of measures (with $\mu_n \in \Prob{X_n}{\mathcal{E}_n}$ for each
$n \ge 0$) and for each $n \ge 0$ let $\nu_n = \mu_n f_n^{-1}  \in \Prob{\specM}{\mathcal{B}}$. Then by Lemma~\ref{kolmogorov}.6
$\nu_n = \nu_{n+1} \Delta^{-1}$ for each $n \ge 0$ and hence by Lemma~\ref{kolmogorov}.5 there exists a unique
measure $\nu \in \Prob{\specM}{\mathcal{B}_\Delta}$ such that $\nu_n = \nu \theta_n^{-1}$ for all $n \ge 0$. But
by Lemma~\ref{uni-meas}.5
$\bar{\nu} \theta_n^{-1} = \bar{\nu}_n = \bar{\mu}_n f_n^{-1}$ and therefore
 \[\bar{\nu}(\theta_n^{-1}(f_n(X_n))) = \bar{\nu}_n(f_n(X_n)) = \bar{\mu}_n(f_n^{-1}(f_n(X_n)))  = \bar{\mu}_n(X_n) = 1\] 
for each $n \ge 0$. This shows that 
$\bar{\nu}(q(X)) = \bar{\nu}\bigl(\bigcap_{n \ge 0} \theta_n^{-1}(f_n(X_n))\bigr) = 1$.  

Proposition~\ref{image-measures}.1 and Lemma~\ref{image-measures}.1 thus imply there exists a measure $\mu \in \Prob{X}{\mathcal{E}}$ with
$\nu = \mu q^{-1}$. Let $n \ge 0$ and $E \in \mathcal{E}_n$; 
then $E = f_n^{-1}(B)$ with $B \in \mathcal{B}$ and
\begin{eqnarray*}
 \mu_n(E) = \mu_n(f_n^{-1}(B)) = \nu_n(B) = \nu(\theta_n^{-1}(B)) &=& \mu(q^{-1}(\theta_n^{-1}(B)))\\ 
&=& \mu(\tau_n^{-1}(f_n^{-1}(B))) = \mu(\tau_n^{-1}(E)) 
\end{eqnarray*}
which shows that $\mu_n \tau_n^{-1} = \mu$ for all $n \ge 0$.
Finally, $\mu$ is the unique measure with this property, since by definition $\mu$ is determined by the sequence $\{\mu_n\}_{n\ge 0}$ 
on the algebra $\bigcup_{n \ge 0} \tau_n^{-1}(\mathcal{E}_n)$ and
$\mathcal{E} = \sigma\big(\bigcup_{n \ge 0} \tau_n^{-1}(\mathcal{E}_n)\bigr)$. \eop


\startsection{Finite point processes}
\label{pp}

The following is a fundamental construction in the theory of point processes:
For a measurable space $(X,\mathcal{E})$ let $\pp{X}{}$
denote the set of all measures on $(X,\mathcal{E})$ taking only values in the set $\Nat$
(and so each $p \in \pp{X}{}$ is a finite measure, since $p(X) \in \Nat$); put 
$\pp{\mathcal{E}}{} = \sigma(\mathcal{E}_\Diamond)$, where
$\mathcal{E}_\Diamond$ is the set of all subsets of $\pp{X}{}$ having the form
$\{ p \in \pp{X}{} : p(E) = k \}$ with $E \in \mathcal{E}$ and $k \in \Nat$.

In this section we give a proof of the following result (which is well-known
to those working in point processes):

\begin{proposition}
If $(X,\mathcal{E})$ is a type $\mathcal{B}$ space then so is $(\pp{X}{},\pp{\mathcal{E}}{})$.
\end{proposition}

A proof of this, or of results which are equivalent to it, can be found in
Matthes, Kerstan and Mecke \cite{mkm}, Kallenberg \cite{kall} and Bourbaki \cite{bourb}.

We start with some constructions which will be needed in our proof of this result.
Let $(X,\mathcal{E})$ and $(Y,\mathcal{F})$ be measurable spaces and let $f : (X,\mathcal{E}) \to (Y,\mathcal{F})$ be a measurable
mapping. If $p \in \pp{X}{}$ and $p f^{-1}$ is the image measure on $(Y,\mathcal{F})$ then
$(p f^{-1})(F) = p(f^{-1}(F)) \in \Nat$ for all $F \in \mathcal{F}$ and so $p f^{-1} \in \pp{Y}{}$.
Thus there is a mapping $\pp{f}{} : \pp{X}{} \to \pp{Y}{}$ given by
$\pp{f}{}(p) = p f^{-1}$ for each $p \in \pp{X}{}$.

\begin{lemma}
(1)\enskip The mapping $\pp{f}{} : (\pp{X}{},\pp{\mathcal{E}}{}) \to (\pp{Y}{},\pp{\mathcal{F}}{})$ is measurable. 

(2)\enskip If $f$ is exactly measurable then so is $\pp{f}{}$.

(3)\enskip If $f$ is exactly measurable and $f(X) \in \mathcal{F}$ then
$\pp{f}{}(\pp{X}{}) \in \pp{\mathcal{F}}{}$.
\end{lemma}

\proof
(1)\enskip Let $F \in \mathcal{F}$ and $k \in \Nat$; then
\begin{eqnarray*}
\lefteqn{ \pp{f}{-1}(\{ q \in \pp{Y}{} : q(F) = k \})}\hspace{60pt} \\
&=& \{ p \in \pp{X}{} : \pp{f}{}(p)(F) = k \} = \{ p \in \pp{X}{} : p(f^{-1}(F)) = k \}\;.
\end{eqnarray*}
Thus $\pp{f}{-1}(\mathcal{F}_\Diamond) \subset \mathcal{E}_\Diamond$ and therefore 
\[ \pp{f}{-1}(\pp{\mathcal{F}}{}) = \pp{f}{-1}(\sigma(\mathcal{F}_\Diamond)) 
= \sigma (\pp{f}{-1}(\mathcal{F}_\Diamond)) 
\subset \sigma(\mathcal{E}_\Diamond) = \pp{\mathcal{E}}{}\;.\]

(2)\enskip
Let $E \in \mathcal{E}$ and $k \in \Nat$; then there exists $F \in \mathcal{F}$ with $f^{-1}(F) = E$ and the calculation in (1) 
shows that
\[\pp{f}{-1}(\{ q \in \pp{Y}{} : q(F) = k \}) = \{ p \in \pp{X}{} : p(E) = k \}\;.\]
This implies 
$\pp{f}{-1}(\mathcal{F}_\Diamond) = \mathcal{E}_\Diamond$
(since in (1) we showed that $\pp{f}{-1}(\mathcal{F}_\Diamond) \subset \mathcal{E}_\Diamond$).
Therefore 
$\pp{f}{-1}(\pp{\mathcal{F}}{}) =  \pp{f}{-1}(\sigma(\mathcal{F}_\Diamond)) =  \sigma(\pp{f}{-1}(\mathcal{F}_\Diamond)) 
= \sigma(\mathcal{E}_\Diamond) = \pp{\mathcal{E}}{}$. 

(3)\enskip
Put $f(X) = D$ and so $D \in \mathcal{F}$. If $p \in \pp{X}{}$ then
\[(p f^{-1})(D) = p(f^{-1}(D)) = p(X) = p(f^{-1}(Y)) = (p f^{-1})(Y) \;.\] 
On the other hand, if $q \in \pp{Y}{}$ with $q(D) = q(Y)$ then by Proposition~\ref{image-measures}.1 there exists a measure $p$ on $(X,\mathcal{E})$
with $p f^{-1} = q$. Moreover,
$p(f^{-1}(F)) = q(F) \in \Nat$ for all $F \in \mathcal{F}$ 
and $f^{-1}(\mathcal{F}) = \mathcal{E}$ and so it follows that  $p \in \pp{X}{}$. Therefore
\begin{eqnarray*}
 \pp{f}{}(X) = \{ \pp{f}{}(p) : p \in \pp{X}{} \} &=& \{ p f^{-1} : p \in \pp{X}{} \} \\
&=& \bigcup_{n \in \Nat} \{ q \in \pp{Y}{} : q(D) = n \} \cap \{ q \in \pp{Y}{} : q(Y) = n \} 
\end{eqnarray*}
and hence $\pp{f}{}(X) \in \pp{\mathcal{F}}{}$. \eop

It is also useful to partition the space $\pp{X}{}$ into components consisting of those measures having the same total measure, and for this
we recall the definition of the $\sigma$-algebra occurring in the  disjoint union of measurable spaces.
Let $S$ be a non-empty set and for each $s \in S$ let $(Y_s,\mathcal{F}_s)$ be a measurable space.
Assume the sets $Y_s$, $s \in S$, are disjoint and put $Y = \bigcup_{s\in S} Y_s$. Then
\[\mathcal{F} = \{ A \subset Y : A \cap Y_s \in \mathcal{F}_s \ \mbox{for each}\ s \in S \}\]
is a $\sigma$-algebra of subsets of $Y$ and $(Y,\mathcal{F})$ is called the \definition{disjoint union} of the measurable spaces
$(Y_s,\mathcal{F}_s)$, $s \in S$.

Now for each $n \in \Nat$ let $\pp{X}{n}$
denote the set of all measures $p$ on $(X,\mathcal{E})$ taking only values in the set $\Nat_n = \{0,1,\ldots,n\}$
and with $p(X) = n$; put 
$\pp{\mathcal{E}}{n} = \sigma(\mathcal{E}^n_\Diamond)$, where
$\mathcal{E}^n_\Diamond$ is the set of all subsets of $\pp{X}{n}$ having the form
$\{ p \in \pp{X}{n} : p(E) = k \}$ with $E \in \mathcal{E}$ and $k \in \Nat_n$.
Thus $\pp{X}{}$ is the disjoint union of the sets
$\pp{X}{n}$, $n \in \Nat$. 

\begin{lemma}
$\pp{\mathcal{E}}{} = \{ A \subset \pp{X}{} : A \cap \pp{X}{n} \in \pp{\mathcal{E}}{n} \ \mbox{for each}\ n \in \Nat \}$
and thus the measurable space $(\pp{X}{},\pp{\mathcal{E}}{})$ is the disjoint union of the measurable spaces
$(\pp{X}{n},\pp{\mathcal{E}}{n})$, $n \in \Nat$.
\end{lemma}

\proof 
Put $\mathcal{D} = \{ A \subset \pp{X}{} : A \cap \pp{X}{n} \in \pp{\mathcal{E}}{n} \ \mbox{for each}\ n \in \Nat \}$, so $\mathcal{D}$
is the $\sigma$-algebra in the definition of the disjoint union.

Let $\pp{\mathcal{D}}{n} = \{ A \cap \pp{X}{n} : A \in \pp{\mathcal{E}}{} \}$;
then $\pp{\mathcal{D}}{n}$ is the trace $\sigma$-algebra of $\pp{\mathcal{E}}{}$ on $\pp{X}{n}$
and thus $\pp{\mathcal{D}}{n} = \sigma(\mathcal{D}^n_\Diamond)$, where
$\mathcal{D}^n_\Diamond = \{ A \cap \pp{X}{n} : A \in \mathcal{E}_\Diamond \}$.
But $\mathcal{D}^n_\Diamond = \mathcal{E}^n_\Diamond$ and hence
$\pp{\mathcal{D}}{n} = \pp{\mathcal{E}}{n}$, i.e.,
$\pp{\mathcal{E}}{n} = \{ A \cap \pp{X}{n} : A \in \pp{\mathcal{E}}{} \}$.
Therefore if $A \in \pp{\mathcal{E}}{}$ then 
$A \cap \pp{X}{n} \in \pp{\mathcal{E}}{n}$ for each $n \in \Nat$, which implies that $A \in \mathcal{D}$.
This shows $\pp{\mathcal{E}}{} \subset \mathcal{D}$. 

Conversely, let $A \in \mathcal{D}$; 
then $A \cap \pp{X}{n} \in \pp{\mathcal{E}}{n}$ and thus there exists $A_n \in \pp{\mathcal{E}}{}$
with $A \cap \pp{X}{n} = A_n \cap \pp{X}{n}$ and this implies that
$A \cap \pp{X}{n} \in \pp{\mathcal{E}}{}$ for each $n \in \Nat$, since $\pp{X}{n} \in \pp{\mathcal{E}}{}$.
Finally, we then have $A = \bigcup_{n\in \Nat} (A \cap \pp{X}{n}) \in \pp{\mathcal{E}}{}$, i.e., 
$\mathcal{D} \subset \pp{\mathcal{E}}{}$, and hence $\mathcal{D} = \pp{\mathcal{E}}{}$. \eop

We can now describe the main steps in the proof of Proposition~\ref{pp}.1, thus let $(X,\mathcal{E})$ be a type $\mathcal{B}$ space.
Then there exists an exactly measurable mapping $f : (X,\mathcal{E}) \to (\specM,\mathcal{B})$ with $f(X) \in \mathcal{B}$.
Therefore by Lemma~\ref{pp}.1 the mapping
$\pp{f}{} : (\pp{X}{},\pp{\mathcal{E}}{}) \to (\pp{\specM}{},\pp{\mathcal{B}}{})$ is exactly measurable
and $\pp{f}{}(\pp{X}{}) \in  \pp{\mathcal{B}}{}$, and so by Proposition~\ref{type-ba}.2~(2) it is enough to show that 
$(\pp{\specM}{},\pp{\mathcal{B}}{})$ is a type $\mathcal{B}$ space.
But by Lemma~\ref{pp}.2 $(\pp{\specM}{},\pp{\mathcal{B}}{})$ is the disjoint union of the measurable spaces
$(\pp{\specM}{n},\pp{\mathcal{B}}{n})$, $n \in \Nat$, and 
if $(Y,\mathcal{F})$ is the disjoint union of type $\mathcal{B}$ spaces $(Y_n,\mathcal{F}_n)$, $n \in \Nat$, then 
by Proposition~\ref{type-ba}.2~(4) $(Y,\mathcal{F})$
is also a type $\mathcal{B}$ space. It is thus enough to show that $(\pp{\specM}{n},\pp{\mathcal{B}}{n})$ is a type $\mathcal{B}$ space for each $n \in \Nat$.

Now fix $n \in \Nat$. We consider $\pp{\specM}{n}$ as a topological space:
Let $\pp{\mathcal{U}}{n}$ be the set of all non-empty subsets of $\pp{\specM}{n}$ having the form
\[\{ p \in \pp{\specM}{n}  : p(C) = v_C \ \mbox{for all}\ C \in N \}\]
with $N$ a finite subset of $\mathcal{C}_{\specM}$ and $\{v_C\}_{C \in N}$ a sequence from $\Nat_n$.
Clearly for each $p \in \pp{\specM}{n}$ there exists $U \in \pp{\mathcal{U}}{n}$ with $p \in U$ and if $U_1,\,U_2 \in \pp{\mathcal{U}}{n}$ and 
$p \in U_1 \cap U_2$
then there exists $U \in \pp{\mathcal{U}}{n}$ with $p \in U \subset U_1 \cap U_2$. Thus $\pp{\mathcal{U}}{n}$ is the base for a topology
$\pp{\mathcal{O}}{n}$ on $\pp{\specM}{n}$. This means that $U \in \pp{\mathcal{O}}{n}$ if and only if for each $p \in U$ there exists a finite
subset $N$ of $\mathcal{C}_{\specM}$ such that
\[\{ q \in \pp{\specM}{n}  : q(C) = p(C) \ \mbox{for all}\ C \in N \} \subset U\;.\]

\begin{lemma}
The topological space $\pp{\specM}{n}$ is compact and metrisable
and $\pp{\mathcal{B}}{n}$ is the Borel $\sigma$-algebra of $\pp{\specM}{n}$. In particular, 
$(\pp{\specM}{n},\pp{\mathcal{B}}{n})$ is a type $\mathcal{B}$ space.
\end{lemma}

\proof 
We start by showing that the topology $\pp{\mathcal{O}}{n}$ on $\pp{\specM}{n}$ is given by a metric.
Let $\{C_k\}_{k\ge 1}$ be an enumeration of the elements in the countable set $\mathcal{C}_{\specM}$
and define a mapping $\varrho : \pp{\specM}{n} \times \pp{\specM}{n} \to \Realpos$ by
\[ \varrho(p,q) = \sum_{k\ge 1} 2^{-k} |p(C_k) - q(C_k)|\;. \]
If $\varrho(p,q) = 0$ then $p(C) = q(C)$ for all $C \in \mathcal{C}_{\specM}$ and hence $p = q$ (since $\mathcal{C}_{\specM}$ is an algebra
with $\sigma(\mathcal{C}_{\specM}) = \mathcal{B}$).
Thus $\varrho$ is a metric since by definition it is symmetric and it is clear that the triangle inequality holds.
Moreover, if $p \in \pp{\specM}{n}$ then for each $\varepsilon > 0$ there exists a finite subset $N$ of $\mathcal{C}_{\specM}$ with
\[\{ q \in \pp{\specM}{n}  : q(C) = p(C) \ \mbox{for all}\ C \in N \} \subset \{ q \in \pp{\specM}{n} : \varrho(q,p) < \varepsilon \} \]
and for each finite subset $N$ of $\mathcal{C}_{\specM}$ there exists $\varepsilon > 0$ such that
\[\{ q \in \pp{\specM}{n} : \varrho(q,p) < \varepsilon \}  \subset \{ q \in \pp{\specM}{n}  : q(C) = p(C) \ \mbox{for all}\ C \in N \}\;. \]
This means that $\pp{\mathcal{O}}{n}$ is the topology given by the metric $\varrho$.
Note that if $\{p_k\}_{k\ge 1}$ is a sequence from $\pp{\specM}{n}$ and $p \in \pp{\specM}{n}$ then
$\lim_{k} p_k = p$ (i.e., $\lim_{k} \varrho(p_k,p) = 0$)
if and only if $\lim_{k} p_k(C) = p(C)$ for each $C \in \mathcal{C}_{\specM}$.

In order to show that $\pp{\specM}{n}$ is compact it is enough to show that the metric space $\pp{\specM}{n}$ is sequentially compact.
Let $\{p_k\}_{k\ge 1}$ be a sequence of elements of $\pp{\specM}{n}$.
By the usual diagonal argument there exists a subsequence $\{k_j\}_{j \ge 1}$ such that
$\lim_{j} p_{k_j}(C)$ exists for each $C \in \mathcal{C}_{\specM}$.
Define $p : \mathcal{C}_{\specM} \to \Realpos$ by $p(C) = \lim_{j} p_{k_j}(C)$.
Then $p$ is clearly finitely additive and $p(\specM) = m$ and so by Proposition~\ref{image-measures}.2
$p$ is a measure on $(\specM,\mathcal{C}_{\specM})$ which has a unique
extension to a measure (also denoted by $p$) on $(\specM,\mathcal{B})$. But
$\mathcal{D} = \{ B \in \pp{\mathcal{B}}{n} : p(B) \in \Nat_n \}$ is a monotone class containing the algebra
$\mathcal{C}_{\specM}$
and thus $p \in \pp{\specM}{n}$. Therefore $p \in \pp{\specM}{n}$ and $\lim_{j} \varrho(p_{k_j},p) = 0$ and
this shows that the metric space $\pp{\specM}{n}$ is sequentially compact.

It remains to show that $\pp{\mathcal{B}}{n}$ is the Borel $\sigma$-algebra of $\pp{\specM}{n}$.
First, the set $\pp{\mathcal{U}}{n}$ is countable and so each element of $\pp{\mathcal{O}}{n}$ can be written as a countable
union of elements from $\pp{\mathcal{U}}{n}$. Thus $\pp{\mathcal{O}}{n} \subset \sigma(\pp{\mathcal{U}}{n})$, which implies that
$\sigma(\pp{\mathcal{O}}{n}) = \sigma(\pp{\mathcal{U}}{n})$, since $\pp{\mathcal{U}}{n} \subset \pp{\mathcal{O}}{n}$. 
Second, each element of $\pp{\mathcal{U}}{n}$ is a finite intersection of elements from $\mathcal{B}^n_\Diamond$ and hence
$\pp{\mathcal{U}}{n} \subset \pp{\mathcal{B}}{n}$. This shows that
$\sigma(\pp{\mathcal{O}}{n}) = \sigma(\pp{\mathcal{U}}{n}) \subset \pp{\mathcal{B}}{n}$.
Finally, let $k \in \Nat_n$ and let 
$\mathcal{D}$ be the set of those $B \in \mathcal{B}$ for which
$\{ p \in \pp{\specM}{n} : p(B) = k \} \in \sigma(\pp{\mathcal{O}}{n})$.
Then $\mathcal{C}_{\specM} \subset \mathcal{D}$ and $\mathcal{D}$ is a monotone class, and so by the monotone class theorem
$\mathcal{D} = \mathcal{B}$, and this means that
$\{ p \in \pp{\specM}{n} : p(B) = k \} \in \sigma(\pp{\mathcal{O}}{n})$ for all $B \in \mathcal{B}$, $k \in \Nat_n$, i.e.,
$\mathcal{B}^n_\Diamond \subset \sigma(\pp{\mathcal{O}}{n})$. Thus
$\pp{\mathcal{B}}{n} = \sigma(\mathcal{B}^n_\Diamond) \subset \sigma(\pp{\mathcal{O}}{n})$, and this shows
$\pp{\mathcal{B}}{n} =  \sigma(\pp{\mathcal{O}}{n})$. \eop

This completes the proof of Proposition~\ref{pp}.1. \eop


\startsection{Existence of conditional distributions}
\label{cond-dist}

Let us say that 
\definition{conditional distributions} exist for measurable spaces $(X,\mathcal{E})$ and $(Y,\mathcal{F})$
if for each measure $\mu \in \Prob{X \times Y}{\mathcal{E}\times \mathcal{F}}$ there exists 
a probability kernel $\pi : X \times \mathcal{F} \to \Realpos$ such that
\[\mu(E \times F) =   \mu_1(I_E \pi(I_F))\]
 for all $E \in \mathcal{E}$, $F \in \mathcal{F}$, where
$\mu_1 = \mu \proj_1^{-1}$ is the image measure of $\mu$ under the
projection $\proj_1 : X \times Y \to X$ onto the first component.
(Beware that this definition is not symmetric in $(X,\mathcal{E})$ and $(Y,\mathcal{F})$.)
By a probability kernel we here mean a mapping $\pi : X \times \mathcal{F} \to \Realpos$ such that
$\pi(x,\cdot) \in \Prob{Y}{\mathcal{F}}$ for each $x \in X$ and such that $\pi(\cdot,F) : X \to \Realpos$ is
$\mathcal{E}$-measurable for each $F \in \mathcal{F}$.

Conditional distributions do not exist in general. However, they do exist if
$(X,\mathcal{E})$ is countably generated and $(Y,\mathcal{F})$ is a type $\mathcal{B}$ space.
Proofs of this fact can be found in Chapter~1 of Doob \cite{doob}, Chapter~V of Parthasarathy \cite{partha},
and also in Appendix~4 of Dynkin and Yushkevich \cite{dynkinyu}. We also give a proof:

\begin{theorem}
Conditional distributions exist for $(X,\mathcal{E})$ and $(Y,\mathcal{F})$ 
if $(X,\mathcal{E})$ is countably generated and $(Y,\mathcal{F})$ is a type $\mathcal{B}$ space.
\end{theorem}

\proof We reduce things to the case in which $(X,\mathcal{E}) = (Y,\mathcal{F}) = (\specM,\mathcal{B})$.

\begin{lemma}
Conditional distributions exist for $(\specM,\mathcal{B})$ and $(\specM,\mathcal{B})$.
\end{lemma}

\proof 
For $m \ge 1$ again let
$\mathcal{C}_m = \proq_m^{-1}(\mathcal{P}(\{0,1\}^{m}))$, where
$\proq_m : \specM \to \{0,1\}^{m}$ is given by $\proq_m(\{z_n\}_{n\ge 1}) = (z_1,\ldots,z_m)$.
Thus $\{\mathcal{C}_m\}_{m\ge 1}$ is an increasing sequence of finite algebras with 
$\mathcal{C}_\specM = \bigcup_{m\ge 0} \mathcal{C}_m$.
For each $z \in \specM$ and each $n \ge 1$ let $\mathrm{a}_n(z)$ be the atom of $\mathcal{C}_n$ containing $z$.
Let $\mathcal{N} \subset \mathcal{B}$ be the trivial $\sigma$-algebra with $\mathcal{N} = \{\varnothing,M\}$.

Let $\mu \in \Prob{\specM\times \specM}{\mathcal{B}\times \mathcal{B}}$ and  $\mu_1 = \mu\proj_1^{-1}$ with
$\proj_1 : \specM \times \specM \to \specM$ projecting onto the first component. For each $n \ge 1$ define
$\gamma_n : \specM \times \mathcal{B} \to \Realpos$ by
\[ \gamma_n(z,B) = \frac{\mu(\mathrm{a}_n(z) \times B)}{\mu_1(\mathrm{a}_n(z))}\]
with $0/0$ taken to be $0$. 
Then $\gamma_n(z,\cdot)$ is either $0$ or an element of $\Prob{\specM}{\mathcal{B}}$ for each $z \in \specM$,
$\gamma_n(\cdot,B) \in \Mappi(\mathcal{C}_n)$ for each $B \in \mathcal{B}$ and
$\mu(C \times B) =  \mu_1(I_C \gamma_n(I_B))$
for all $C \in \mathcal{C}_n$, $B \in \mathcal{B}$.
Consider the mapping $\gamma_n' : (\specM \times \specM) \times \mathcal{B} \to \Realpos$ with
$\gamma_n'((z_1,z_2),B) = \gamma(z_1,B)$; then $\gamma_n'(\cdot,B) \in \Mappi(\mathcal{C}_n\times \mathcal{N})$
and
\[ 
\mu(I_{C\times N} \gamma'_n(I_B)) = \mu((C \times N) \cap (\specM \times B)) =  \mu(I_{C\times N}I_{\specM\times B}) \]
for all $C \in \mathcal{C}_n$, $N \in \mathcal{N}$, $B \in \mathcal{B}$. 
Thus $\gamma'_n(I_ B)$ is a version of the conditional expectation of $I_{\specM\times B}$ with respect to
$\mathcal{C}_n\times \mathcal{N}$ for each $n \ge 1$ and it therefore follows from 
the martingale convergence theorem (see, for example, Breiman, \cite{breiman}, Theorem~5.24) that
\begin{eqnarray*}  
\lefteqn{\mu_1\bigl(\bigl\{ z \in \specM : \lim\limits_{n\to\infty} \gamma_n(z,B)
  \ \mbox{exists} \bigr\}\bigr)}\hspace{80pt} \\
&=& \mu\bigl( \bigl\{ (z_1,z_2) \in \specM\times \specM : \lim\limits_{n\to\infty} \gamma'_n((z_1,z_2),B) \ \mbox{exists} \bigr\}\bigr) = 1
\end{eqnarray*}
for each $B \in \mathcal{B}$. Put
\[\specM_{\mathcal{C}} = \bigl\{ z \in \specM : \lim\limits_{n\to\infty} \gamma_n(z,C)
  \ \mbox{exists for all}\ C \in \mathcal{C}_\specM \bigr\}\;;\]
since $\mathcal{C}_\specM$ is countable it follows that $\specM_{\mathcal{C}} \in \mathcal{B}$ and $\mu_1(\specM_{\mathcal{C}}) = 1$.
Choose $z_0 \in \specM_{\mathcal{C}}$ and
define a mapping $\gamma : \specM \times \mathcal{C}_\specM \to \Realpos$ by letting
\[ \gamma(z,C) =  \left\{ \begin{array}{cl}
                     \lim_{n} \gamma_n(z,C) &\ \mbox{if}\ z \in \specM_\mathcal{C}\;,\\
                     \lim_{n} \gamma_n(z_0,C) &\ \mbox{if}\ z \in \specM \setminus \specM_\mathcal{C}\;.
                \end{array} \right. \]
Then $\gamma(\cdot,C) \in \Mappi(\mathcal{B})$ for each $C \in \mathcal{C}_\specM$ 
and by the dominated convergence theorem
\[ \mu(C_1 \times C_2) = \lim_{n\to\infty} \mu_1(I_{C_1} \gamma_n(I_{C_2})) = \mu_1(I_{C_1} \gamma(I_{C_2}))\]
for all $C_1,\,C_2 \in \mathcal{C}_\specM$.
Hence by $\mu(B_1 \times C_2) = \mu_1(I_{B_1} \gamma(I_{C_2}))$
for all $B_1 \in \mathcal{B}$, $C_2 \in \mathcal{C}_\specM$.
Now it is clear that the mapping $\gamma(z,\cdot) : \mathcal{C}_\specM \to \Realpos$ is additive
with $\gamma(z,\specM) = 1$ for each $z \in \specM$
and so by Proposition~\ref{image-measures}.2 it has a unique extension to an element of $\Prob{\specM}{\mathcal{B}}$
which will also be denoted by $\gamma(z,\cdot)$. It follows that
$\gamma : \specM \times \mathcal{B} \to \Realpos$ is a probability kernel satisfying
$\mu(B_1 \times B_2) = \mu_1(I_{B_1} \gamma(I_{B_2}))$ for all $B_1,\,B_2 \in \mathcal{B}$.
\eop

\begin{lemma}
If $(Y,\mathcal{F})$ is a type $\mathcal{B}$ space then
conditional distributions exist for $(\specM,\mathcal{B})$ and $(Y,\mathcal{F})$. 
\end{lemma}

\proof 
Since $(Y,\mathcal{F})$ is a type $\mathcal{B}$ space there exists an exactly measurable mapping 
$f : (Y,\mathcal{F}) \to (\specM,\mathcal{B})$ with $f(Y) \in \mathcal{B}$.
Let $\mu \in \Prob{\specM \times Y}{\mathcal{B}\times \mathcal{F}}$ and let $\mu_1 = \mu \proj_1^{-1}$
with $\proj_1 : \specM \times Y \to \specM$ the projection onto the first component.
Put $\nu = \mu g^{-1}$, where $g = \id_\specM\times f : \specM \times Y \to \specM \times \specM$,
so $\nu \in \Prob{\specM\times\specM}{\mathcal{B}\times \mathcal{B}}$. 
Then by Lemma~\ref{cond-dist}.1 there exists a probability
kernel $\gamma : \specM \times \mathcal{B} \to \Realpos$ such that
\[\nu(B_1 \times B_2) = \nu_1(I_{B_1} \gamma(I_{B_2}))\]
for all $B_1,\,B_2 \in \mathcal{B}$, where
$\nu_1 = \nu \breve{\proj}_1^{-1}$ with $\breve{\proj}_1 : \specM \times \specM \to \specM$ projecting onto the first component,
and note that $\nu_1 = \mu_1$, since $\breve{\proj}_1 \circ g = \breve{\proj}_1 \circ (\id_\specM \times f) = \proj_1$.
Now consider $\specM_0 = \{ z \in \specM : \gamma(z,f(Y)) = 1 \}$; then $\specM_0 \in \mathcal{B}$ and $\nu_1(\specM_0) = 1$, since
\[ 1 = \mu(\specM \times Y) = \mu(g^{-1}(\specM \times f(Y))) 
= \nu(\specM \times f(Y)) = \nu_1( \gamma(f(Y))) \;.\]
Choose some point $z_0 \in \specM_0$ and define $\gamma_o : \specM \times \mathcal{B} \to \Realpi$  by
\[ \gamma_o(z,B) = \left\{ \begin{array}{cl}
        \gamma(z,B)  &\ \mbox{if}\ z \in \specM_0\;,\\
        \gamma(z_0,B)  &\ \mbox{if}\ z \in \specM \setminus \specM_0\;;
\end{array} \right. \]       
then $\gamma_o$ is a probability kernel with $\gamma_o(z,f(Y)) = 1$ for all $z \in \specM$ and
\[\nu(B_1 \times B_2) = \nu_1(I_{B_1} \gamma_o(I_{B_2})) \]
for all $B_1,\,B_2 \in \mathcal{B}$. Now by Proposition~\ref{image-measures}.1 there exists 
for each $z \in \specM$ a probability measure $\tau(z,\cdot) \in \Prob{Y}{\mathcal{F}}$ so that
$\tau(z,f^{-1}(B)) = \gamma_o(z,B)$ for all $B \in \mathcal{B}$, and then 
$\tau : \specM \times \mathcal{F} \to \Realpos$ is clearly a probability kernel.
Let $B \in \mathcal{B}$ and $F \in \mathcal{F}$; then $F = f^{-1}(B')$ for some $B' \in \mathcal{B}$ and so
\begin{eqnarray*}
\mu(B \times F) &=& \mu(B \times f^{-1}(B')) = \mu(g^{-1}(B \times F))
= \nu(B \times B')\\
&=& \nu_1(I_{B} \gamma_o(I_{B'}))  = \mu_1(I_{B} \tau(I_{f^{-1}(B')}))
= \mu_1(I_{B} \tau(I_F))
\end{eqnarray*}
and this shows that 
conditional distributions exist for $(\specM,\mathcal{B})$ and $(Y,\mathcal{F})$. \eop

\textit{Proof of Theorem~\ref{cond-dist}.1:\ }
By Proposition~\ref{cg-meas-spaces}.2 there exists an exactly measurable mapping $f : (X,\mathcal{E}) \to (\specM,\mathcal{B})$.
Let $\mu \in \Prob{X\times Y}{\mathcal{E}\times \mathcal{F}}$ and $\mu_1 = \mu \proj_1^{-1}$
with $\proj_1 : X \times Y \to X$ the projection onto the first component.
Put $\nu = \mu g^{-1}$, where $g = f\times\id_Y : X \times Y \to \specM \times Y$, so
$\nu \in \Prob{\specM\times Y}{\mathcal{B}\times \mathcal{F}}$. Then by Lemma~\ref{cond-dist}.2 there exists a probability
kernel $\tau : \specM \times \mathcal{F} \to \Realpos$ such that
$\nu(B \times E) = \nu_1(I_{B} \tau(I_F))$
for all $B \in \mathcal{B}$, $F \in \mathcal{F}$, where
$\nu_1 = \nu \breve{\proj}_1^{-1}$ with $\breve{\proj}_1 : \specM \times Y \to \specM$ the projection onto the first component,
and $\nu_1 = \mu_1 f^{-1}$, since
$\breve{\proj}_1 \circ g = \breve{\proj}_1 \circ (f \times \id_Y) = f \circ \proj_1$.
Now define $\pi : X \times \mathcal{F} \to \Realpos$ by letting $\pi(x,F) = \tau(f(x),F)$ for all
$x \in X$, $F \in \mathcal{F}$, thus $\pi$ is clearly a probability kernel.
Let $E \in \mathcal{E}$, $F \in \mathcal{F}$; then $E = f^{-1}(B)$ for some $B \in \mathcal{B}$ and so
\begin{eqnarray*}
\mu(B \times F) &=& \mu(f^{-1}(B) \times F) = \mu(g^{-1}(B \times F))
= \nu(B \times F)\\
&=& \nu_1(I_{B} \tau(I_F)) =  (\mu_1 f^{-1})(I_{B} \tau(I_F)) \\
&=&  \mu_1(I_{f^{-1}(B)} \tau(f(\cdot),F))
= \mu_1(I_{E} \pi(I_F))
\end{eqnarray*}
and therefore conditional distributions exist for $(X,\mathcal{E})$ and $(Y,\mathcal{F})$. 
This completes the proof of Theorem~\ref{cond-dist}.1. \eop


\startsection{The Dynkin extension property}
\label{dynkin}

In what follows let $(X,\mathcal{E})$ be a measurable space.
A mapping $\pi : X \times \mathcal{E} \to \Realpos$ is a \definition{quasi probability kernel} if
$\pi(x,\cdot) : \mathcal{E} \to \Realpos$ is a measure with
$\pi(x,X)$ either $0$ or $1$ for each $x \in X$ and $\pi(\cdot,E) : X \to \Realpos$ is $\mathcal{E}$-measurable
for each $E \in \mathcal{E}$. If 
$\pi(x,X) = 1$ (i.e., $\pi(x,\cdot) \in \Prob{X}{\mathcal{E}}$)
for each $x \in X$ then $\pi$ is called a \definition{probability kernel}.

If $\mathcal{E}'$ is a sub-$\sigma$-algebra of $\mathcal{E}$ then
a quasi probability kernel $\pi(x,\cdot) : \mathcal{E} \to \Realpos$ is said to be
\definition{$\mathcal{E}'$-measurable} if the mapping $\pi(\cdot,E) : X \to \Realpos$ is $\mathcal{E}'$-measurable for each 
$E \in \mathcal{E}$.

If $\mathcal{E}'$ is a sub-$\sigma$-algebra of $\mathcal{E}$
and $\pi : X \times \mathcal{E} \to \Realpos$ is an $\mathcal{E}'$-measurable quasi probability kernel then let
\[\mathcal{G}(\pi) = \Bigl\{ \mu \in \Prob{X}{\mathcal{E}} :
     \mu(E' \cap E) = \int_{E'} \pi(x,E)\,d\mu(x)\ \mbox{for all}\ E' \in \mathcal{E}',\,E \in \mathcal{E}\Bigr\}\;,\]
thus if $\Elem_\mu(I_E | \mathcal{E}')$ denotes the conditional expectation of $I_E$ 
with respect to the measure $\mu$ and the sub-$\sigma$-algebra $\mathcal{E}'$, then in fact
\[\mathcal{G}(\pi) = \{ \mu \in \Prob{X}{\mathcal{E}} :
    \Elem_\mu(I_E | \mathcal{E}') = \pi(\cdot,E)\ \mbox{$\mu$-a.e.\ for all}\ E \in \mathcal{E}\}\;.\]

Now let $\{\mathcal{E}_n\}_{n\ge 0}$ be a decreasing sequence of sub-$\sigma$-algebras of $\mathcal{E}$
and denote the \definition{tail field} $\bigcap_{n\ge 0} \mathcal{E}_n$ by  $\mathcal{E}_\infty$.
A sequence of kernels $\{\pi_n\}_{n\ge 0}$ is 
\definition{adapted to} $\{\mathcal{E}_n\}_{n\ge 0}$ if 
$\pi_n : X \times \mathcal{E} \to \Realpos$  is an 
$\mathcal{E}_n$-measurable quasi probability kernel for each $n \ge 0$.
The sequence $\{\mathcal{E}_n\}_{n\ge 0}$ has the \definition{Dynkin extension property} if
for each sequence  $\{\pi_n\}_{n\ge 0}$ adapted to $\{\mathcal{E}_n\}_{n\ge 0}$
there exists an $\mathcal{E}_\infty$-measurable quasi probability kernel $\pi : X \times \mathcal{E} \to \Realpos$ such that
$\bigcap_{n\ge 0} \mathcal{G}(\pi_n) \subset \mathcal{G}(\pi)$.
(Of course, in general the set $\bigcap_{n\ge 0} \mathcal{G}(\pi_n)$ will be empty,
since no consistency assumptions have been placed on the kernels $\{\pi_n\}_{n\ge 0}$.)

This property does not hold in general. However, 
if $(X,\mathcal{E})$ is a type $\mathcal{B}$ space then any decreasing sequence of sub-$\sigma$-algebras of $\mathcal{E}$
has the Dynkin extension property. This is proved in F\"ollmer \cite{foellmer} (based on ideas in  
Dynkin \cite{dynkin}); another proof can be found in Chapter~7 of Georgii \cite{georgii}. 

In Theorem~\ref{dynkin}.1 we establish that the Dynkin extension property holds for a type $\mathcal{B}$ space. In the second half of the
section we give F\"ollmer's construction in \cite{foellmer} (based on a technique from Dynkin \cite{dynkin})
which shows how the kernel occurring in Theorem~\ref{dynkin}.1 can be improved to obtain one which is much more suitable
for applications. This refinement does not depend on properties of type $\mathcal{B}$ spaces, except in that it needs
the kernel from Theorem~\ref{dynkin}.1 as a starting point.

\begin{theorem}
If $(X,\mathcal{E})$ is a type $\mathcal{B}$ space then any decreasing sequence of sub-$\sigma$-algebras of $\mathcal{E}$
has the Dynkin extension property.
\end{theorem}

\proof
Let $\{\mathcal{E}_n\}_{n\ge 0}$ be a decreasing sequence of sub-$\sigma$-algebras of $\mathcal{E}$
and let 
$\{\pi_n\}_{n\ge 0}$ be a sequence of kernels
adapted to $\{\mathcal{E}_n\}_{n\ge 0}$, thus  $\pi_n : X \times \mathcal{E} \to \Realpos$  is an 
$\mathcal{E}_n$-measurable quasi probability kernel for each $n \ge 0$.
We are looking for an $\mathcal{E}_\infty$-measurable quasi probability kernel $\pi$ such that
\[ \mu(G \cap E) = \int_G \pi(\cdot,E)\,d\mu \]
for all  $G \in \mathcal{E}_\infty$, $E \in \mathcal{E}$
and all $\mu \in \mathcal{G}$, where $\mathcal{G} = \bigcap_{n\ge 0} \mathcal{G}(\pi_n)$ and
\[\mathcal{G}(\pi_n) = \Bigl\{ \mu \in \Prob{X}{\mathcal{E}} :
     \mu(E' \cap E) = \int_{E'} \pi_n(\cdot,E)\,d\mu\ \mbox{for all}\ E' \in \mathcal{E}_n,\,E \in \mathcal{E}\Bigr\}\;.\]
(Note that if $\mathcal{G} = \varnothing$ then we can simply take $\pi(x,E) = 0$ for all $x \in X$, $E \in \mathcal{E}$). 

Now since $(X,\mathcal{E})$ is a type $\mathcal{B}$ space there exists an exactly measurable mapping 
$f : (X,\mathcal{E}) \to (\specM,\mathcal{B})$ with $f(X) \in \mathcal{B}$.
As before let $\mathcal{C}_\specM \subset \mathcal{B}$ be the countable algebra of cylinder sets. Let
\[X_{\mathcal{C}} = \bigl\{ x \in X: \lim\limits_{n\to\infty} \pi_n(x,f^{-1}(C))
  \ \mbox{exists for all}\ C \in \mathcal{C}_\specM \bigr\}\;.\]

\begin{lemma}
$X_{\mathcal{C}} \in \mathcal{E}_\infty$ and $\mu(X_{\mathcal{C}}) = 1$ for each $\mu \in \mathcal{G}$. 
\end{lemma}

\proof For each $C \in \mathcal{C}_\specM$ let $X_C$ denote the set of those elements $x \in X$ for which the limit
$\lim_{n\to\infty} \pi_n(x,f^{-1}(C))$ exists, thus $X_{\mathcal{C}} = \bigcap_{C \in \mathcal{C}_\specM} X_C$
and therefore, since $\mathcal{C}_\specM$ is countable, it is enough to show for each $C \in \mathcal{C}_\specM$ that 
$X_C \in \mathcal{E}_\infty$ and $\mu(X_C) = 1$ for each $\mu \in \mathcal{G}$. 
Now clearly $X_C \in \mathcal{E}_\infty$, and $\mu(X_C) = 1$ holds for each $\mu \in \mathcal{G}$ since
$\pi_n(\cdot,f^{-1}(C))$ is a version of $\Elem_\mu(I_{f^{-1}(C)}|\mathcal{E}_n)$
for each $n \ge 0$ and by the martingale convergence theorem (see, for example, Breiman, \cite{breiman}, Theorem~5.24) 
it follows that
$\lim_{n \to \infty} \Elem_\mu(I_{f^{-1}(C)}|\mathcal{E}_n) = \Elem_\mu(I_{f^{-1}(C)}|\mathcal{E}_\infty)$
$\mu$-a.e.  \eop

Define a mapping $\tau : X \times \mathcal{C}_\specM \to \Realpos$ by letting
\[ \tau(x,C) =  \left\{ \begin{array}{cl}
                     \lim\limits_{n\to \infty} \pi_n(x,f^{-1}(C)) &\ \mbox{if}\ x \in X_\mathcal{C}\;,\\
                                 0                 &\ \mbox{if}\ x \in X \setminus X_\mathcal{C}\;.
                \end{array} \right. \]
                            
\begin{lemma}
For each $C \in \mathcal{C}_\specM$ the mapping $\tau(\cdot,C) : X \to \Realpos$
is $\mathcal{E}_\infty$-measurable and
$\mu(G \cap f^{-1}(C)) = \int_G \tau(\cdot,C)\,d\mu$ for all $G \in \mathcal{E}_\infty$, $\mu \in \mathcal{G}$.
\end{lemma}

\proof 
It is clear that $\tau(\cdot,C)$ is $\mathcal{E}_\infty$-measurable, since by Lemma~\ref{dynkin}.1 
$X_{\mathcal{C}} \in \mathcal{E}_\infty$. Moreover, if $\mu \in \mathcal{G}$ and $G \in \mathcal{E}_\infty$ 
then for each $n \ge 0$
\[ \mu(G \cap f^{-1}(C)) = \int_G \pi_n(\cdot,f^{-1}(C))\,d\mu\]
and by Lemma~\ref{dynkin}.1 $\mu(X_{\mathcal{C}}) = 1$; thus by the dominated convergence theorem
\[ \mu(G \cap f^{-1}(C)) 
= \lim_{n\to\infty} \int_G I_{X_{\mathcal{C}}}\pi_n(\cdot,f^{-1}(C))\,d\mu
= \int_G \tau(\cdot,C)\,d\mu\;.\ \eop \]

\begin{lemma}
The mapping $\tau(x,\cdot) : \mathcal{C} \to \Realpos$ is additive
with $\tau(x,\specM) \in \{0,1\}$ for each $x \in X$. 
\end{lemma}

\proof 
This is trivially true if $x \in X \setminus X_{\mathcal{C}}$ so let
$x \in X_{\mathcal{C}}$. Clearly $\tau(x,M)$ must be either $0$ or $1$, since
$\pi_n(x,X)$ takes on only these values, and if $C_1,\,C_2 \in \mathcal{C}$
with $C_1 \cap C_2 = \varnothing$ then $f^{-1}(C_1) \cap f^{-1}(C_2) = \varnothing$ and
\begin{eqnarray*}
 \tau(x,C_1 \cup C_2)
 &=& \lim_{n\to \infty} \pi_n(x,f^{-1}(C_1\cup C_2))
= \lim_{n\to \infty} \pi_n(x,f^{-1}(C_1)\cup f^{-1}(C_2))\\
  &=& \lim_{n\to \infty} 
   (\pi_n(x,f^{-1}(C_1)) + \pi_n(x,f^{-1}(C_2)))
 = \tau(x,C_1) + \tau(x,C_2)\;.\ \eop
\end{eqnarray*}

Let $x \in X$; by Proposition~\ref{image-measures}.2 the additive mapping $\tau(x,\cdot) : \mathcal{C}_\specM \to \Realpos$ has a unique extension to a 
measure on $\mathcal{B}$ which will also be denoted by $\tau(x,\cdot)$.
Thus the measure $\tau(x,\cdot)$ is either $0$ or an element 
of $\Prob{M}{\mathcal{B}}$. This defines a mapping $\tau : X \times \mathcal{B} \to \Realpos$. 
 
\begin{lemma}
For each $B \in \mathcal{B}$ the mapping $\tau( \cdot, B) : X \to \Realpos$ is $\mathcal{E}_\infty$-measurable
and $\mu(G \cap f^{-1}(B)) = \int_G \tau(\cdot,B)\,d\mu$ for all $G \in \mathcal{E}_\infty$,
$\mu \in \mathcal{G}$.
\end{lemma}

\proof 
This follows from Lemma~\ref{dynkin}.2 using the monotone class theorem. \eop

Now for the first time the fact that $f(X) \in \mathcal{B}$ will be needed. Let
\[X_f = \{ x \in X : \tau(x,f(X)) = 1 \}\;.\]
Then $X_f \in \mathcal{E}_\infty$ and applying Lemma~\ref{dynkin}.4 with $B = f(X)$ and $G = X$ shows 
\[\int \tau(\cdot,f(X))\,d\mu  = \mu(f^{-1}(f(X))) = \mu(X) = 1\]
and hence that $\mu(X_f) = 1$ for each $\mu \in \mathcal{G}$.
Define a mapping $\eta : X \times \mathcal{B} \to \Realpos$ by 
$\eta(x,B) = I_{X_f}(x)\tau(x,B)$.
Then $\eta(x,\cdot)$ is either $0$ or an element of $\Prob{M}{\mathcal{B}}$ with $\eta(x,f(X)) = 1$
for each $x \in X$, the mapping
$\eta(\cdot,B) : X \to \Realpos$ is $\mathcal{E}_\infty$-measurable for each $B \in \mathcal{B}$ and
$\mu(G \cap f^{-1}(B)) = \int_G \eta(\cdot,B)\,d\mu$
for all $G \in \mathcal{E}_\infty$, $\mu \in \mathcal{G}$,
(since $\mu(X_f) = 1$ for each $\mu \in \mathcal{G}$). 

By Proposition~\ref{image-measures}.1 there now exists for each $x \in X$ a unique measure $\pi(x,\cdot)$ (either $0$ or an element of
$\Prob{X}{\mathcal{E}}$) so that $\eta(x,\cdot) = \pi(x,\cdot) f^{-1}$.
The resulting mapping $\pi : X \times \mathcal{E} \to \Realpos$ is then an
$\mathcal{E}_\infty$-measurable quasi probability kernel. If $E \in \mathcal{E}$ then $E = f^{-1}(B)$
for some $B \in \mathcal{B}$ and thus
\[ \mu(G \cap E) = \mu(G \cap f^{-1}(B))
= \int_G \eta(\cdot,B)\,d\mu = \int_G \pi(\cdot,f^{-1}(B))\,d\mu 
= \int_G \pi(\cdot,E)\,d\mu \]
for all $G \in \mathcal{E}_\infty$, $\mu \in \mathcal{G}$. This completes the proof of Theorem~\ref{dynkin}.1. \eop

We now give F\"ollmer's construction in \cite{foellmer} which shows how the kernel occurring in Theorem~\ref{dynkin}.1 can be improved 
to obtain a much better one. However, for this an additional assumption
(strictness) has to be placed on the sequence of quasi probability kernels.

Let $(X,\mathcal{E})$ be a measurable space and let $\mathcal{E}'$ be a sub-$\sigma$-algebra of $\mathcal{E}$.
Then an $\mathcal{E}'$-measurable quasi probability kernel
$\pi : X \times \mathcal{E} \to \Realpos$   will be called \definition{strict} if
\[\pi(x,E' \cap E) = I_{E'}(x)\pi(x,E)\]
for all $E' \in \mathcal{E}'$, $E \in \mathcal{E}$ and all $x \in X$.

\begin{lemma}
Let $\pi : X \times \mathcal{E} \to \Realpos$ be a strict $\mathcal{E}'$-measurable quasi probability kernel;
then $\mathcal{G}(\pi) = \{ \mu \in \Prob{X}{\mathcal{E}} : \mu = \mu\pi \}$,
where the measure $\mu\pi$ is defined by $(\mu\pi)(E) = \int \pi(\cdot, E)\,d\mu$
for all $E \in \mathcal{E}$.
\end{lemma}

\proof 
If $\mu \in \mathcal{G}(\pi)$ then 
$(\mu\pi)(E) = \int \pi(\cdot, E)\,d\mu = \mu(X \cap E) = \mu(E)$
for all $E \in \mathcal{E}$, i.e., $\mu = \mu\pi$.  Conversely, if $\mu = \mu\pi$ then 
\[ \mu(E' \cap E) = (\mu\pi)(E' \cap E) = \int \pi(\cdot,E'\cap E)\,d\mu = \int_{E'} \pi(\cdot,E)\,d\mu \]
for all $E' \in \mathcal{E}_o$, $E \in \mathcal{E}$, and hence $\mu \in \mathcal{G}(\pi)$. \eop

In what follows let $(X,\mathcal{E})$ be a type $\mathcal{B}$ space and
$\{\mathcal{E}_n\}_{n\ge 0}$ be a decreasing sequence of sub-$\sigma$-algebras of $\mathcal{E}$;
put $\mathcal{E}_\infty = \bigcap_{n\ge 0} \mathcal{E}_n$. Let $\{\pi_n\}_{n\ge 0}$ be a sequence of strict kernels
adapted to $\{\mathcal{E}_n\}_{n\ge 0}$, thus  $\pi_n : X \times \mathcal{E} \to \Realpos$  is  now a strict
$\mathcal{E}_n$-measurable quasi probability kernel for each $n \ge 0$, and so by Lemma~\ref{dynkin}.5
\[\mathcal{G}(\pi_n) = \{ \mu \in \Prob{X}{\mathcal{E}} : \mu\pi_n = \mu \}\;.\]

\begin{theorem}
Again let $\mathcal{G} = \bigcap_{n\ge 0} \mathcal{G}(\pi_n)$.
Then there exists an $\mathcal{E}_\infty$-measurable quasi probability
kernel $\pi : X \times \mathcal{E} \to \Realpos$ with
\[ \mu(G \cap E) = \int_G \pi(\cdot,E)\,d\mu \]
for all  $G \in \mathcal{E}_\infty$, $E \in \mathcal{E}$
and all $\mu \in \mathcal{G}$, such that the following hold:

(1)\enskip $\pi(x,E) \in \{0,1\}$ for all $E \in \mathcal{E}_\infty$, $x \in X$. 

(2)\enskip  $\pi(x,X \setminus \Delta_x) = 0$ for each $x \in X$,
where $\Delta_x = \{ y \in X : \pi(y,\cdot) = \pi(x,\cdot) \}$.

(3)\enskip $\mathcal{G} = \{ \mu \in \Prob{X}{\mathcal{E}} : \mu\pi = \mu \}$.
\end{theorem}

(Note that if $\mathcal{G} = \varnothing$ then there is really nothing to prove:
We can still take $\pi(x,E) = 0$ for all $x \in X$, $E \in \mathcal{E}$,
since here $\{ \mu \in \Prob{X}{\mathcal{E}} : \mu\pi = \mu \} = \varnothing$.)

\proof 
By Theorem~\ref{dynkin}.1 there exists an $\mathcal{E}_\infty$-measurable quasi probability kernel
$\pi' : X \times \mathcal{E} \to \Realpos$ such that
$\mu(G \cap E) = \int_G \pi'(\cdot,E)\,d\mu$
for all  $G \in \mathcal{E}_\infty$, $E \in \mathcal{E}$ and all $\mu \in \mathcal{G}$.
Let us fix a countable algebra $\mathcal{D}$ with $\mathcal{E} = \sigma(\mathcal{D})$.
(This exists since a type $\mathcal{B}$ space is countable generated.)

\begin{lemma}
Let $\mu \in \mathcal{G}$ and
$f : X \to \Realpos$ be bounded and $\mathcal{E}$-measurable; then
\[\int_G f\,d\mu = \int_G \pi'f \,d\mu \]
for all $G \in \mathcal{E}_\infty$ (where $\pi'f$ is defined by
$(\pi'f)(x) = \int f(y)\pi'(x,dy)$ for each $x \in X$).
\end{lemma}

\proof If $f = I_E$ then this is true by assumption, and so the result also holds 
for all simple mappings (i.e., $\mathcal{E}$-measurable mappings taking on only finitely many values).
It therefore holds for a general bounded $\mathcal{E}$-measurable $f$, since such a mapping can be 
uniformly approximated using simple ones. \eop

\begin{lemma}
Let $E' = \{ x \in X: \pi'(x,\cdot) \in \mathcal{G}\}$;
then $E' \in \mathcal{E}_\infty$ and $\mu(E') = 1$ for each
$\mu \in \mathcal{G}$. 
\end{lemma}

\proof By the monotone class theorem $x \in E'$ if and only if $\pi'(x,X) = 1$ and
\[ \pi'(x,E) = \int \pi_n(y,E)\pi'(x,dy)\]
for all $n \ge 0$ and all $E \in \mathcal{D}$.
But there are only countably many equations involved here and therefore 
$E' \in \mathcal{E}_\infty$. 
Let $\mu \in \mathcal{G}$, $E \in \mathcal{D}$ and put
$f = \pi_n(\cdot,E)$. By Lemma~\ref{dynkin}.6 (noting that $f$ is bounded) it then follows that
for each $G \in \mathcal{E}_\infty$ 
\begin{eqnarray*}
\lefteqn{ \int_G \int \pi_n(y,E)\pi'(x,dy)\,d\mu(x)
= \int_G \pi' f\,d\mu  = \int_G f\,d\mu}\hspace{40pt}\\
 &=& \int_G \pi_n(x,E)\,d\mu(x) = \int \pi_n(x,E\cap G)\,d\mu(x)\\ 
 &=& (\mu\pi_n)(E\cap G) = \mu(E\cap G)  = \int_G \pi'(x,E)\,d\mu(x)\;,
\end{eqnarray*}
and this shows that
$\pi'(x,E) = \int \pi_n(y,E)\pi'(x,dy)$ holds for $\mu$-a.e.\ $x \in X$. 
Finally, $\int \pi'(\cdot,X)\,d\mu = \mu(X) = 1$ and so $\pi'(\cdot,X) = 1$ $\mu$-a.e. 
Therefore $\mu(E') = 1$. \eop

For each $x \in X$ let $\Delta'_x = \{ y \in X : \pi'(y,\cdot) = \pi'(x,\cdot) \}$; then
\[\Delta'_x = \{ y \in X : \pi'(y,E) = \pi'(x,E)\ \mbox{for all}\ E \in \mathcal{D} \}\]
by the monotone class theorem, and so $\Delta'_x \in \mathcal{E}_\infty$.

\begin{lemma}
Let $E = \{ x \in E':  \pi'(x,\Delta'_x) = 1 \}$;
then $E \in \mathcal{E}_\infty$ and $\mu(E) = 1$ for all $\mu \in \mathcal{G}$. 
\end{lemma}

\proof Note that $E = \bigcap_{E \in \mathcal{D}} E_E$, where
\[ E_E = \Bigl\{ x \in E': \int ( \pi'(x,E) - \pi'(y,E))^2 \pi'(x,dy) = 0 \Bigr\}\;.\]
Let $x \in E'$;
then $\int \pi'(y,E)\pi'(x,dy) = \pi'(x,E)$
(since $\pi'(x,\cdot) \in \mathcal{G}$), and thus
\begin{eqnarray*} 
\lefteqn{\int ( \pi'(x,E) - \pi'(y,E))^2 \pi'(x,dy)}\hspace{40pt}\\
 & =& \int ((\pi'(x,E))^2 - 2\pi'(x,E)\pi'(y,E) + (\pi'(y,E))^2 ) \pi'(x,dy)\\
 &=& \int (\pi'(y,E))^2\pi'(x,dy) - (\pi'(x,E))^2\;.
\end{eqnarray*}
Therefore
$\int (\pi'(y,E))^2\pi'(x,dy) \ge (\pi'(x,E))^2$
for all $x \in E'$ and all $E \in \mathcal{E}$,
and $E_E$ consists exactly of those elements $x \in E'$ for which 
\[\int (\pi'(y,E))^2\pi'(x,dy) = (\pi'(x,E))^2\;.\]
In particular this implies that $E_E \in \mathcal{E}_\infty$. 
Now let $\mu \in \mathcal{G}$ and put
$g = (\pi'(\cdot,E))^2$; then $g$ is bounded and so by Lemma~\ref{dynkin}.6 
\begin{eqnarray*} 
\int\!\int (\pi'(y,E))^2\pi'(x,dy)\,d\mu(x)
 &=& \int\!\int g(y)\pi'(x,dy)\,d\mu(x)  = \int \pi'g\,d\mu\\
 &=&  \int g \,d\mu = \int (\pi'(x,E))^2 \,d\mu(x)\;;
\end{eqnarray*}
hence (since $\int (\pi'(y,E))^2\pi'(x,dy) \ge (\pi'(x,E))^2$ for all $x \in E'$)
\begin{eqnarray*}
\lefteqn{\int \Bigl| \int (\pi'(y,E))^2\pi'(x,dy) - (\pi'(x,E))^2 \Bigr|\,d\mu(x)}\hspace{40pt}\\
 &=& \int\!\int \bigl( (\pi'(y,E))^2\pi'(x,dy) - (\pi'(x,E))^2\bigr) \,d\mu(x)\\
 &=& \int\!\int (\pi'(y,E))^2\pi'(x,dy)\,d\mu(x) - \int (\pi'(x,E))^2 \,d\mu(x) =  0\;,
\end{eqnarray*}
and thus $\mu(E_E) = 1$. Since $\mathcal{D}$ is countable it then follows that both
$E \in \mathcal{E}_\infty$ and $\mu(E) = 1$. \eop

\begin{lemma}
If $x \in E$ then $\pi'(x,E) \in \{0,1\}$ for each $E \in \mathcal{E}_\infty$. 
\end{lemma}

\proof Let $x \in E$; then $\pi'(x,\cdot) \in \mathcal{G}$ and $\pi'(x,\Delta'_x) = 1$. Thus for each
$E \in \mathcal{E}_\infty$
\begin{eqnarray*}
 \pi'(x,E) &=& \pi'(x,E\cap E) = \int_E \pi'(y,E)\pi'(x,dy) = \int_E I_{\Delta'_x}(y)\pi'(y,E)\pi'(x,dy)\\
  &=& \int_E I_{\Delta'_x}(y)\pi'(x,E)\pi'(x,dy)  = \int_E \pi'(x,E)\pi'(x,dy)
  = (\pi'(x,E))^2
\end{eqnarray*}
and hence $\pi'(x,E) \in \{0,1\}$. \eop

Now define $\pi : X \times \mathcal{E} \to \Realpos$ by
\[ \pi(x,E) = \left\{ \begin{array}{cl}
    \pi' (x,E)  &\ \mbox{if}\ x \in E\;,\\
    0  &\ \mbox{if}\ x \notin E\;.
\end{array} \right. \]
Then $\pi$ is clearly an $\mathcal{E}_\infty$-measurable quasi probability kernel, and 
\[ \mu(G \cap E) = \int_G \pi'(\cdot,E)\,d\mu = \int_G I_E \pi'(\cdot,E)\,d\mu 
= \int_G \pi(\cdot,E)\,d\mu \]
for all  $G \in \mathcal{E}_\infty$, $E \in \mathcal{E}$
and all $\mu \in \mathcal{G}$ (since $\mu(E) = 1$). Moreover, (1) holds: This follows from Lemma~\ref{dynkin}.9 
if $x \in E$ and it is trivially true if $x \in X \setminus E$.

If $x \in E$ then $\Delta_x = \Delta'_x \cap E$ and by Lemma~\ref{dynkin}.8 $\pi'(x,E) = 1$ (since $\pi'(x,\cdot) \in \mathcal{G}$).
Therefore $\pi(x,\Delta_x) = \pi'(x,\Delta'_x \cap E) = \pi'(x,\Delta'_x) = 1$, and so
$\pi(x,X\setminus \Delta_x) = 0$. But this is trivially true if $x \in X \setminus E$, and hence (2) holds. 

Finally (3) also holds: Let $\mu \in \Prob{X}{\mathcal{E}}$ with $\mu = \mu\pi$; then in particular
\[ \mu(E) = \int I_E \,d\mu = \int \pi(\cdot,X)\,d\mu = (\mu\pi)(X) = \mu(X) = 1\;. \]
Moreover, if $x \in E$ then $\pi(x,\cdot) \in \mathcal{G}$ and so
\[ \pi(x,E) = (\pi(x,\cdot) \pi_n)(E) = \int \pi_n(y,E)\pi(x,dy)\]
for all $E \in \mathcal{E}$, $n \ge 0$. Therefore
\begin{eqnarray*}
 (\mu\pi_n)(E) &=& \int \pi_n(\cdot,E)\,d\mu = \int \pi_n(\cdot,E)\,d(\mu\pi)\\
&=& \int \int \pi_n(y,E)\pi(x,dy)\,d\mu(x)   = \int_E \int \pi_n(y,E)\pi(x,dy)\,d\mu(x) \\
&=& \int_E \pi(x,E)\,d\mu(x) = \int \pi(x,E)\,d\mu(x) = (\mu\pi)(E) = \mu(E) 
\end{eqnarray*}
for all $E \in \mathcal{E}$, $n \ge 0$, i.e., $\mu = \mu\pi_n$ for all $n \ge 0$, and so $\mu \in \mathcal{G}$.
Conversely, if $\mu \in \mathcal{G}$ then
$\mu(E) = \mu(X \cap E) = \int_X \pi(\cdot,E)\,d\mu = (\mu\pi)(E)$ for all $E \in \mathcal{E}$, i.e., 
$\mu = \mu\pi$. This completes the proof of Theorem~\ref{dynkin}.2. \eop


\sbox{\ttt}{\textsc{Bibliography}}
\thispagestyle{plain}

\bigskip
\bigskip

{\sc Fakult\"at f\"ur Mathematik, Universit\"at Bielefeld}\\
{\sc Postfach 100131, 33501 Bielefeld, Germany}\\
\textit{E-mail address:} \texttt{preston@math.uni-bielefeld.de}\\
\textit{URL:} \texttt{http://www.math.uni-bielefeld.de/\symbol{126}preston}


\end{document}